\documentclass[11 pt]{amsart}

\usepackage{amsxtra,comment}
\usepackage[all]{xy}
\usepackage{palatino}
\usepackage{mathtools}
\usepackage{latexsym, amsmath, amssymb, longtable, booktabs,amscd,microtype,booktabs, cases}
\usepackage[square]{natbib}
\usepackage{graphicx}
\usepackage[dvipsnames]{xcolor}
\usepackage{caption}
\usepackage{tcolorbox}
\usepackage[all]{xy}
\usepackage{enumerate}
\usepackage{multirow}
\usepackage{tikz}
\usepackage{amsthm}

\theoremstyle{plain}
\newtheorem{theorem}{Theorem}[section]
\newtheorem{thm}[theorem]{Theorem}
\newtheorem{proposition}{Proposition}[section]
\newtheorem{lemma}[theorem]{Lemma}
\newtheorem{corollary}{Corollary}[section]
\newtheorem{cor}[corollary]{Corollary}

\theoremstyle{definition}
\newtheorem{definition}{Definition}[section]
\newtheorem{remark}{Remark}[section]
\newtheorem{example}{Example}[section]

\usetikzlibrary{matrix,arrows,decorations.pathmorphing}
\usepackage{collectbox}

\usepackage{amsxtra}
\usepackage[colorinlistoftodos,prependcaption, textsize=tiny]{todonotes}
\reversemarginpar

\usepackage{enumitem}
\usepackage{mathabx}

\usepackage{latexsym, amsmath, amssymb, longtable, booktabs,amscd,microtype,booktabs}
\usepackage[square]{natbib}

\usepackage{color}
\usepackage{amsfonts}
\usepackage{amsmath}
\usepackage[all]{xy}
\usepackage{geometry}
\usepackage{amssymb,mathrsfs}
\usepackage{amsmath}
\usepackage{amsthm}
\usepackage{amsopn,xcolor}
\usepackage{amscd}
\usepackage{exscale, verbatim}
\usepackage{graphicx}
\usepackage{hyperref}
\usepackage{bbold, bbm}
\bibpunct{[}{]}{,}{n}{}{;} 

\numberwithin{equation}{section}
\begin{document}
\renewcommand{\baselinestretch}{1.25}

\newcommand\p{\mathfrak{p}}
\newcommand\m{\mathfrak{m}}
\newcommand\J{\mathrm{J}}
\newcommand\A{\mathbb{A}}
\newcommand{\CC}{\mathbb{C}}
\newcommand\G{\mathbb{G}}
\newcommand\N{\mathbb{N}}
\newcommand{\T}{\mathbb{T}}
\newcommand{\cT}{\mathcal{T}}
\newcommand{\cL}{\mathcal{L}}
\newcommand\sO{\mathcal{O}}
\newcommand\sE{{\mathcal{E}}}
\newcommand\tE{{\mathbb{E}}}
\newcommand\sF{{\mathcal{F}}}
\newcommand\sG{{\mathcal{G}}}
\newcommand\sC{{\mathcal{C}}}
\newcommand\GL{{\mathrm{GL}}}
\newcommand{\HH}{\mathrm H}
\newcommand\mM{{\mathrm M}}
\newcommand\fS{\mathfrak{S}}
\newcommand\fP{\mathfrak{P}}
\newcommand\fQ{\mathfrak{Q}}
\newcommand\Qbar{{\bar{\Q}}}
\newcommand\sQ{{\mathcal{Q}}}
\newcommand\sP{{\mathbb{P}}}
\newcommand{\Q}{\mathbb{Q}}
\newcommand{\tH}{\mathbb{H}}
\newcommand{\Z}{\mathbb{Z}}
\newcommand{\R}{\mathbb{R}}
\newcommand{\I}{\mathcal{I}}
\newcommand{\F}{\mathbb{F}}
\newcommand{\D}{\mathfrak{D}}
\newcommand{\Div}{\mathrm{Div}}
\newcommand\gP{\mathfrak{P}}
\newcommand\Gal{{\mathrm {Gal}}}
\newcommand\SL{{\mathrm {SL}}}
\newcommand\Spec{{\mathrm {Spec}}}
\newcommand\Num{{\mathrm {Num}}}

\newcommand\Hom{{\mathrm {Hom}}}
\newcommand{\legendre}[2] {\left(\frac{#1}{#2}\right)}
\newcommand\iso{{\> \simeq \>}}
\newcommand{\cX}{\mathcal{X}}
\newcommand{\cJ}{\mathcal{J}}

\theoremstyle{definition}
\newtheorem{dfn}{Definition}

\theoremstyle{remark}

\makeatletter
\def\imod#1{\allowbreak\mkern10mu({\operator@font mod}\,\,#1)}
\makeatother
\newcommand{\mcF}{\mathcal{F}}
\newcommand{\mcG}{\mathcal{G}}
\newcommand{\mcM}{\mathcal{M}}
\newcommand{\mcO}{\mathcal{O}}
\newcommand{\mcP}{\mathcal{P}}
\newcommand{\mcS}{\mathcal{S}}
\newcommand{\mcV}{\mathcal{V}}
\newcommand{\mcW}{\mathcal{W}}

\newcommand{\cO}{\mathcal{O}}
\newcommand{\cC}{\mathcal{C}}
\newcommand{\mfa}{\mathfrak{a}}
\newcommand{\mfb}{\mathfrak{b}}
\newcommand{\mfH}{\mathfrak{H}}
\newcommand{\mfh}{\mathfrak{h}}
\newcommand{\mfm}{\mathfrak{m}}
\newcommand{\mfn}{\mathfrak{n}}
\newcommand{\mfp}{\mathfrak{p}}
\newcommand{\mfq}{\mathfrak{q}}

\newcommand{\mrB}{\mathrm{B}}
\newcommand{\mrG}{\mathrm{G}}
\newcommand{\gG}{\mathcal{G}}
\newcommand{\mrH}{\mathrm{H}}
\newcommand{\mH}{\mathrm{H}}
\newcommand{\mrZ}{\mathrm{Z}}
\newcommand{\Ga}{\Gamma}
\newcommand{\cyc}{\mathrm{cyc}}
\newcommand{\Fil}{\mathrm{Fil}}

\newcommand{\mrp}{\mathrm{p}}
\newcommand{\PGL}{\mathrm{PGL}}
\newcommand{\x}{{\mathcal{X}}}
\newcommand{\Sp}{\textrm{Sp}}
\newcommand{\ab}{\textrm{ab}}
\newcommand{\Jac}{\mathrm{Jac}}
\newcommand{\Frob}{\mathrm{Frob}}
\newcommand{\lra}{\longrightarrow}
\newcommand{\ra}{\rightarrow}
\newcommand{\rai}{\hookrightarrow}
\newcommand{\ras}{\twoheadrightarrow}

\newcommand{\repr}{\rho_{f,\wp}|_{G_p}}
\newcommand{\GRF}{{\rho}_{f,\wp}}

\newcommand{\lan}{\langle}
\newcommand{\ran}{\rangle}
\newcommand{\ord}{\mathrm{ord}}
\newcommand{\crit}{\mathrm{crit}}
\newcommand{\old}{\mathrm{old}}
\newcommand{\new}{\mathrm{new}}

\newcommand{\mo}[1]{|#1|}

\newcommand{\hw}[1]{#1+\frac{1}{2}}
\newcommand{\mcal}[1]{\mathcal{#1}}
\newcommand{\trm}[1]{\textrm{#1}}
\newcommand{\mrm}[1]{\mathrm{#1}}
\newcommand{\car}[1]{|#1|}
\newcommand{\pmat}[4]{ \begin{pmatrix} #1 & #2 \\ #3 & #4 \end{pmatrix}}
\newcommand{\bmat}[4]{ \begin{bmatrix} #1 & #2 \\ #3 & #4 \end{bmatrix}}

\newcommand{\pbmat}[4]{\left \{ \begin{pmatrix} #1 & #2 \\ #3 & #4 \end{pmatrix} \right \}}
\newcommand{\psmat}[4]{\bigl( \begin{smallmatrix} #1 & #2 \\ #3 & #4 \end{smallmatrix} \bigr)}

\newcommand{\bsmat}[4]{\bigl[ \begin{smallmatrix} #1 & #2 \\ #3 & #4 \end{smallmatrix} \bigr]}
\newcommand{\bstwomat}[2]{\bigl[ \begin{smallmatrix} #1 \\ #2 \end{smallmatrix} \bigr]}
\title{Cuspidal subgroups associated with non-rational Eisenstein maximal ideals}
\author{Debargha Banerjee}
\email{debargha.banerjee@gmail.com}
\address{Department of Mathematics, 
Indian Institute of Science Education and Research Pune,
Dr. Homi Bhabha Road, Pune 411008, INDIA.}
\author{Narasimha Kumar}
\email{narasimha@math.iith.ac.in}
\address{
Department of Mathematics, Indian Institute of Technology Hyderabad,
Kandi, Sangareddy 502284, INDIA.}
\author{Dipramit Majumdar}
\email{dipramit@gmail.com}
\address{Department of Mathematics, Indian Institute of Technology Madras, Chennai 600036, INDIA,  and Institut für Mathematik, Universität Heidelberg, GERMANY.
}
\thanks{It is a pleasure to acknowledge e-mail communications, advice, and remarks of Professors Ken Ribet and Barry Mazur. 
The SERB grant ANRF/ARGM/$2025/000421$/MTR has partially supported the first author. The second author's research was partially supported by the SERB grant No. CRG$/2023/003136$. The third author would like to thank Jo\"{e}l Bella\"{i}che for his guidance and encouragement. He had numerous discussions with Jo\"{e}l on Eisenstein ideals during his Ph.D., for which he is forever grateful. The third author thanks Professor Otmar Venjakov for hosting him in the University of Heidelberg as a visiting researcher. We are sincerely grateful to the anonymous referee for carefully reading our paper and for comments that helped us to improve the mathematical content of the paper.}

\subjclass[2010]{Primary: 11F67, Secondary: 11F11, 11F30, 11F80}
\keywords{Eisenstein series, Eisenstein ideal, Special values of $L$-functions, Congruence of modular forms}

\begin{abstract}
 In this paper, we are interested in the generalization of Ramanujan-like Eisenstein congruences (congruences between cusp forms and Eisenstein series) for congruence subgroups of the form $\Ga_0(N)$ with $N \in \N$. 
  We determine the possible primes that can produce Eisenstein congruences.  We provide several examples of Eisenstein congruences to substantiate our method. Ribet conjectured (\cite[p. 360]{MR3540618}) about these congruences for the square-free level $N$. Yoo proved the conjecture. For general $N$, Yoo proved a generalization of the conjecture, under some hypotheses, provided that those ideals are {\it rational}. 
We show that the generalization of Ribet's conjecture for certain non-square-free levels $N$ is true even for {\it non-rational} Eisenstein maximal ideals.
\end{abstract}
\maketitle
\begin{center} {\it Dedicated to Jo\"{e}l Bella\"{i}che.} \end{center}
 \section{Introduction}
\label{Introduction}
Ramanujan was the first to discover congruences between cusp forms and Eisenstein series. Modulo the miraculous prime $691$, there is a mysterious congruence between the cusp forms and the Eisenstein series for level $N=1$ and weight $12$. 
 If $N$ is a prime number,  Ogg computed~\cite{MR337974} the order of the cuspidal subgroup. 
 The order is $\widetilde{n}:=\Num(\frac{N-1}{12})$ (traditionally, it is denoted by $n$). This is the order of the cuspidal subgroup associated with the unique Eisenstein series of level $N$.  Eisenstein congruences, that is, congruences between cusp forms and the Eisenstein series of weight two, are controlled by the number $\widetilde{n}=\Num(\frac{N-1}{12})$, in the sense that the primes dividing $\widetilde{n}$ control Ramanujan-like congruences between Eisenstein series and cusp forms. 
 In a landmark paper on Eisenstein ideals,  Mazur~\cite{MR488287} showed that the cuspidal subgroup is of the same order as the torsion subgroup of the Jacobian of the modular curve, famously proving Ogg's conjecture. 

It is natural to ask if it is possible to predict numbers like $\widetilde{n}=\Num(\frac{N-1}{12})$ for which we have Ramanujan-like congruences between the Eisenstein series and the cusp forms, and to give explicit descriptions of the maximal Eisenstein ideals for the general level $N$. Our paper provides an answer to this question. We follow the strategy of Mazur (similar to that of Ohta, Yoo, and Ren). 

We now state our results.  Let $\J_0(N):= \Jac(X_0(N))/\Q$ be the Jacobian abelian variety \cite[Definition 6.1.1]{MR2112196} associated with the modular curve $X_0(N)$ (cf. \S \ref{maximalideals}) and $C_{\Ga_0(N)}$ be the cuspidal subgroup within $J_0(N)$
generated by the equivalence classes of divisors of degree zero supported at $\partial(X_0(N))$. Consider the Hecke algebra $\T(N)$, generated by all Hecke operators $T_l$ for all primes $l$ (including $l \mid N$),  acting on the Jacobian $\J_0(N):=\Jac(X_0(N))$ of the modular curves $X_0(N)$ \cite[\S 6.3, p. 230]{MR2112196}. Let $\m$ be a maximal ideal of $\T(N)$ with residue characteristic $\ell$ and denote by $\rho_{\m}:\Gal(\overline{\Q}/\Q) \rightarrow \GL_2(\T/\m)$ the two-dimensional semi-simple Galois representation
associated to $\m$.  Up to an isomorphism, this Galois representation is completely determined by the following two properties \cite[p. 2430]{MR3988582}:
\begin{enumerate}
\item 
The Galois representation $\rho_{\m}$ is unramified outside $\ell N$, 
\item 
for all $p \nmid \ell \cdot N$, the characteristic polynomial of $\rho_{\m}(\Frob_p)$ is $X^2- (T_p \pmod \m)X+p$.
\end{enumerate}

We call $\m$ an Eisenstein maximal ideal if $\rho_{\m}$ is reducible. Following Yoo, we first classify Eisenstein maximal ideals as rational and non-rational. 
\begin{definition}
\label{Rationalnon-rational}
(Rational and non-rational Eisenstein maximal ideals)
 We call an Eisenstein maximal ideal $\m$ {\it rational} if 
 $\rho_{\m} \simeq  \mathbb{1} \bigoplus \overline{\chi}_\ell$,  where $\overline{\chi}_\ell$ denotes the mod $\ell$ cyclotomic character. Otherwise, we call $\m$ {\it non-rational}.
\end{definition}

Given a natural number $N$, we define the square full and  square-free part of $N$ as
\begin{equation}\label{decomp}
N_1=\prod\limits_{p; v_p(N) \geq 2} p^{v_p(N)} \text{  and  } N_2=\prod\limits_{p; v_p(N) =1}  p
\end{equation}
and therefore $N=N_1N_2$  with $N_1$ and $N_2$ coprime integers.
 We prove the following result:

\begin{theorem}[Theorem~\ref{classification} in the text]
\label{mainthmclassification}
Let $N=N_1N_2$ be the decomposition of positive integers $N$ with $N_1, N_2$ as in equation~\ref{decomp}. Let $\m$ be a non-rational Eisenstein maximal ideal of $\T(N)$ whose residual characteristic $\ell$ is coprime to $6N_1$. The ideal $\m$ determines a non-trivial Dirichlet character
$\overline{\epsilon}_\m$ of the conductor $f$ (with $f^2 \mid N_1$). Moreover, there exist square-free natural numbers $t$ (coprime to $f$ and dividing $N_1$) and $M$  (coprime to $ft$ and dividing $N$) and an ideal $I_{\overline{\epsilon}_\m,M, t}(N)$ such that $\m = (\ell, I_{\overline{\epsilon}_\m,M, t}(N))$ and $ \T/\m \cong \F_\ell[\overline{\epsilon}_\m]$.

\end{theorem}

 Faltings-Jordan \cite[p. 37, Conjecture $3.22$]{MR1336315} conjectured that every Eisenstein maximal ideal corresponds to a characteristic-zero Eisenstein series for all weights. 
 Our investigation confirms this for the non-rational Eisenstein maximal ideals of weight $2$. 
Let $E_2(\Ga_0(N), \CC)$ be the vector space of all Eisenstein series of weight two associated with the congruence subgroup $\Ga_0(N)$. 
In our study of the classification of non-rational Eisenstein maximal ideals of level $N=N_1 N_2$, we see that if a non-rational Eisenstein maximal ideal $\m$ has a residual characteristic coprime to $6N_1$, then $\m$ corresponds to an Eisenstein series $E_\m:= E_{\epsilon_\m, Mt, Lt} \in E_2(\Ga_0(N),\CC)$ (see Proposition \ref{Em}).

The novelty in our work lies in the use of the cuspidal subgroup $C_{\Ga_0(N)}(E)$ (cf. \S ~\ref{Cuspidal}) associated with an Eisenstein series for the congruence subgroup $\Ga_0(N)$ rather than for $\Ga_1(N)$ with dependence on the Shimura subgroup as in Ren~\cite{MR4621854}.  In~\cite{MR800251}, Stevens computed the order of the subgroup for $\Ga_1(N)$. As a by-product, it can compute the order of $C_{\Ga_0(N)}(E)$ {\it up to}  
the order of intersection of $C_{\Ga_1(N)}(E)$ with the Shimura subgroup, and therefore, it can only compute the order of certain components (~\cite[Example 4.9, p. 542]{MR800251}) of the cuspidal group associated with an Eisenstein series, rather than the entire group due to dependency on the Shimura subgroup. We modify the arguments of Stevens' Theorem
~\cite[Theorem 1.3]{MR800251} and prove a variant of it (cf. Theorem~\ref{threestep})  that holds for $\Ga_0(N)$ for $N \in \N$.

To understand Ramanujan-like congruences for non-rational Eisenstein maximal ideals, we need to base change our work to a ring $R$ (instead of $\Z$) inside a carefully chosen number field obtained from a Dirichlet character. We get this Dirichlet character from the Eisenstein series associated with this non-rational Eisenstein maximal ideal.
Instead of computing $C_{\Ga_0(N)}(E)$, we compute the order of the cyclic $R$ module $C_{\Ga_0(N)}(E)\otimes_\Z R$ associated with an Eisenstein series $E\in E_2(\Ga_0(N),\CC)$ directly, without dependence on the order of the Shimura subgroups. However, we need to put a certain restriction on $N$ to achieve this. 
Note that $C_{\Ga_0(N)}(E)$ and the $R$ module $C_{\Ga_0(N)}(E)\otimes_\Z R$ have the same set of primes in support. 
Now we use the same method as Mazur's seminal work on Eisenstein ideals, but after base changing to $R$, to conclude that they indeed control Ramanujan-like congruences. This part is a routine generalization using the same technique as Mazur. 

For the remainder of the paper, we need to work with certain restrictions on $N$. 
 Let $p$ be a prime, $k\ge 2$ an integer, and $N=p^k N_2$ an integer where $N_2$ is a square-free integer such that all prime divisors $q$ of $N_2$ satisfy $q \equiv \pm 1 \pmod{p^{\lfloor \frac{k}{2} \rfloor}}$.
 In other words, we assume $N_1=p^k$ with certain conditions on the prime divisors of $N_2$. In particular, integers of the form $N=p^k$, $N=16N'$ and $N=9N'$ for any square-free integer $N'$, can be written in the form $N=p^k N_2$ with $N_2$ square-free such that all prime divisors $q$ of $N_2$ satisfy $q \equiv \pm 1 \pmod{p^{\lfloor \frac{k}{2} \rfloor}}$.

For an integer $x$, by $\nu_l(x)$ we denote the highest power of $l$ that divides $x$.
Suppose $\varphi$ is a nontrivial primitive Dirichlet character of conductor $f:=p^r$ with $r \leq {\lfloor \frac{\nu_p(N)}{2} \rfloor}$ (this places some restriction on $N$, for example, if $p=2$ then it forces $r \ge 2$ and hence $\nu_2(N) \ge 4$). This condition ensures that $\varphi(q)=\pm 1$ for all $q \mid N_2$. 
Let $M$ be a natural number that is a factor of $N_2$ and let $L$ be a factor of $\frac{N_2}{M}$. 
At primes dividing $M$, we obtain new Eisenstein series through ordinary refinements,  while for primes dividing $L$, we obtain new Eisenstein series by critical refinements (cf.~\ref{StevenEisen}). In this way, we define an Eisenstein series $E_{\varphi, M, L} \in E_2(\Ga_0(N),\CC)$. We base change to a ring $R:=\Z[\frac{1}{6N}, \zeta_f, \varphi]$ and study the structure of the corresponding $R$-module. For any $\Z$ module $C$, we denote by $C \otimes R$ the corresponding $R$-module $C \otimes_\Z R$ (suppressing 
$\Z$ in the notation). 

\begin{thm}\label{theorem1}
\label{ordercuspidal}
 Consider the ring $R:=\Z[\frac{1}{6N}, \zeta_f, \varphi]$. Let $p$ be a prime and $N$ be a natural number of the form $N=p^kN_2$ for some $k \ge 2$, where $N_2$ is a square-free integer such that all prime divisors $q$ of $N_2$ satisfy the congruence $q \equiv \pm 1 \pmod{p^{\lfloor \frac{k}{2} \rfloor}}$.
 Let $\varphi$ be a nontrivial Dirichlet character of conductor $f=p^r >1$, with $ r\leq  \lfloor \frac{k}{2} \rfloor$.  For a divisor $M$ of $N_2$, consider the {\it non-rational} Eisenstein series $E_{\varphi, M,L} \in E_2(\Ga_0(N),\CC)$, with $L= \frac{N_2}{M}$  (cf.  Proposition~\ref{lvalue}). 
 The cuspidal subgroup associated with $E_{\varphi,M,L}$ is $\Q(\zeta_{f})$ rational, and therefore $C_{\Ga_0(N)}(E_{\varphi, M,L})$ is a subgroup of $J_0(N)(\Q(\zeta_{f}, \varphi))$.  The $R$-module $C_{\Ga_0(N)}(E_{\varphi, M,L}) \otimes R$ is cyclic:
\[
C_{\Ga_0(N)}(E_{\varphi, M,L}) \otimes R \simeq R/ \left(\beta_{\Ga_0(N), \varphi,M,L}\right),
\]
where $\beta_{\Ga_0(N), \varphi,M,L} \in \Z[ \frac{1}{6N}, \zeta_f, \varphi]$ is a constant given by equation \eqref{betaimp} and $(\beta_{\Ga_0(N), \varphi,M,L})$ is the principal ideal generated by $\beta_{\Ga_0(N), \varphi,M,L}$ inside $R$.
\end{thm}
\begin{cor}
\label{support}
A prime $l$ divides $|C_{\Ga_0(N)}(E_{\varphi, M,L})|$ if and only if $l \mid |  \Z[\frac{1}{6N}, \zeta_f, \varphi]/\left(\beta_{\Ga_0(N), \varphi,M,L}\right) |$.
\end{cor}
As expected, the number appearing on the right is analogous to $\widetilde{n}$ for general levels and non-rational Eisenstein series. 
We need the assumption on $N$ to prove Theorem~\ref{theorem1} because the validity of one of our important results, Theorem~\ref{boundary}, depends on this. 

Ribet conjectured (for square-free $N$) that if $\m$ is an Eisenstein maximal ideal, then $C_{\Ga_0(N)}[\m] \neq 0$ \cite{MR3540618}. This conjecture was settled by Yoo~\cite{MR3540618}. 
Yoo~\cite{MR3988582}  (see also ~\cite{MR4663195}) proved that for {\it the rational} Eisenstein maximal ideal $\m$ with odd residue characteristic $\ell$ and $\ell^2 \nmid N$ (for general $N$), we have  $C_{\Ga_0(N)}[\m] \neq 0$ controlling  the support of Ramanujan-like congruences.

\begin{theorem}
\label{mainthmRibet}(A generalization of Ribet's conjecture)
Let $N$ be an integer as in Theorem~\ref{theorem1}. Let $\m$ be an Eisenstein maximal ideal of the Hecke algebra $\T(N)$ with the residual characteristic $\ell$. 
If ($\ell$,$6N)=1$, then 
$$ C_{\Ga_0(N)}[\m] \neq 0. $$
 \end{theorem}

 Given this maximal ideal $\m$, we associate an Eisenstein ideal $I_{\overline{\epsilon}_\m, M}(N) \subset \T(N)$.  We study the support of $\T(N)/I_{\overline{\epsilon}_\m, M}(N)$, that is, the set of primes $\ell$ such that $(\ell, I_{\overline{\epsilon}_\m, M}(N))$ is an Eisenstein maximal ideal of $\T(N)$. It follows that the residual characteristic of a non-rational Eisenstein ideal $\m$ of $\T(N)$ must divide $6N \cdot |C_{\Ga_0(N)}(E_\m)|$ for the Eisenstein series $E_{\m}$ associated with the maximal ideal $\m$. For some of these primes, there exist congruences between $E_\m$ and a new form of weight $2$. We compute examples of Eisenstein congruences (cf.\S \ref{numerical examples}) using SAGE and LMFDB.

Note that most of the techniques used in the paper go back to the seminal work of Mazur on Eisenstein ideals. 
This has been generalized for non-rational Eisenstein maximal ideals by several authors, including Gross-Lubin
~\cite{MR818355} and Calegari~\cite{MR2196762}, Lang-Wake (cf.~\cite{MR4977337},~\cite{MR4503449}) and a series of papers by Ren (cf.~\cite{MR4621854} and the reference therein). To our knowledge, this is the first paper that provides explicit examples of Eisenstein congruences for non-rational Eisenstein series demonstrated for square full level. We can do that because of Theorem~\ref{ordercuspidal}. In this theorem, we compute the order of the cuspidal subgroups associated with Eisenstein series, independent of the Shimura subgroups. 
\section{Eisenstein series for \texorpdfstring{$\Ga_0(N)$}{Gamma0(N)}}

 For any natural number $N$, the congruence subgroup $\Ga_0(N)$ acts on the complex upper half-plane $\tH$.  Let $Y_0(N):=\Ga_0(N) \backslash \tH$ be the Riemann surface, and $X_0(N)$ be the compactified modular curve obtained by adding a finite set of cusps $\partial(X_0(N))$. Denote by $M_2(\Ga_0(N), \CC)$, the vector space over $\CC$ of the classical modular elliptic forms of weight $2$ and level $N$ as in \cite{MR2112196}. Inside this vector space, let $S_2(\Ga_0(N), \CC)$ (respectively, $E_2(\Ga_0(N), \CC)$) be the space of cusp forms (respectively, Eisenstein series). 

\subsection{Non-rational Eisenstein series}
\label{StevenEisen} 

In this section, we use a refined version of Stevens' Eisenstein series \cite{MR670070}, \cite{MR800251} to obtain a new Eisenstein series $E_{\varphi}$ of level $\Ga_0(f^2)$, which is due to Tang ~\cite{MR1376558}. Note that our Eisenstein series is different from that of Ren~\cite[Definition 4.2]{MR4621854}. The advantage of these Eisenstein series is that they are more amenable to computing subgroups of cuspidal groups, which help us compute the exact order of the cuspidal group; hence, we can predict the Eisenstein congruences.
By ~\cite[p. 850, Remark 2.4]{MR1376558}, the Fourier expansion of $E_\varphi$ at $\infty$ is given by 
$$E_\varphi(z):= E_{\varphi, \varphi^{-1}}(z) = \sum\limits_{c=1}^\infty \sum_{b=1}^\infty \varphi(c)  \varphi^{-1}(b) b q^{bc} = \sum_{n=1}^\infty b_n q^n,$$
 where $q = e^{2 \pi i z}$ and $b_n = \sum\limits_{bc=n} \varphi(c) \varphi^{-1}(b)b$ for all $n \geq 1$.

For an Eisenstein series $E$ of level $\Ga_0(N)$  and a Dirichlet character $\chi$ of conductor $m$ with $(m,N)=1$, we define the twisted $L$-function as the analytic continuation of the series 
$ L(E, \chi, s):=\sum_{n=1}^\infty \frac{a_n(E) \chi(n)}{n^s}. $
The above series is convergent in a half-plane $\mrm{Re}(s)>2$. If $\chi$ is trivial, we write it as $L(E,  s)$.
 Note that Stevens \cite[p. 523, (1.5)]{MR800251} uses the series 
 $\chi^{-1}(N) N^s \sum_{n=1}^\infty \frac{a_n(E) \chi(n)}{n^s}$
  to define $L(E,\chi,s)$ for the subgroups $\Ga(N), \Ga_0(N), \Ga_1(N)$. 
 In \cite{MR800251}, Stevens first considers the Eisenstein series of level $\Ga(N)$, and therefore $a_n(E)$ corresponds to the coefficient of $q_N^n = e^{\frac{2 \pi i nz}{N}}$. 
In {\it loc.\ cit.}, the  Eisenstein series for the subgroups $\Ga_0(N)$ or $\Ga_1(N)$ are obtained as a suitable linear combination of Eisenstein series of $\Ga(N)$  and $ \left(\begin{smallmatrix}
1 & N \\
0 & 1\\
\end{smallmatrix}\right)$ is the smallest parabolic matrix in $\Ga(N)$. The significant difference here is that we chose to follow 
~\cite{MR1376558} and start directly with the Eisenstein series for $\Ga_0(N)$.

\begin{proposition}
\label{Eisenimp}
\begin{enumerate}
\item 
The $L$-function $L(E_{\varphi}, s)$ associated with the Eisenstein series $E_{\varphi}$ is given by
$ L(E_\varphi, s)=L(\varphi, s)L(\varphi^{-1}, s-1).$
More generally, for a Dirichlet character $\chi$ of conductor coprime to $f$, we have
$$ L(E_\varphi, \chi,s)=L(\chi \varphi, s)L(\chi \varphi^{-1}, s-1). 
$$
\item 
The action of the Hecke operators on $E_\varphi$ is given by 
\begin{equation*}
T_l (E_{\varphi}) =
\begin{cases}
(l\varphi^{-1}(l)+ \varphi(l))E_{\varphi}& \  \text{if $l \nmid f$,} \\
0 &  \ \text{if $l \mid f$. } \\
\end{cases}
\end{equation*}
\end{enumerate}
\end{proposition}

\begin{proof}
Part (1) of the Proposition essentially follows
from Proposition~\cite[Proposition 2.3, p. 850]{MR1376558}. 
By~\cite[p. 850, Remark 2.4]{MR1376558}, we have 
\[
a_n(E_{\varphi}) \chi(n) =\chi(n) \sum\limits_{bc=n} \varphi(c) \varphi^{-1}(b)b = \sum\limits_{bc=n} \chi(c)\varphi(c) \chi(b)\varphi^{-1}(b)b =  a_n(E_{\chi \varphi, \chi \varphi^{-1}}).
\]
Now, the required claim follows from ~\cite[Proposition 2.3(b)]{MR1376558}.

Part (2) of the proposition is similar to the proof of ~\cite[Lemma 2.7, p. 852]{MR1376558}, which we sketch below. The action of $T_l$ on $E_\varphi$ is given by
$$ T_l(E_\varphi) = \begin{cases}
\sum\limits_{n=1}^\infty b_{l n}q^n + \sum\limits_{n=1}^\infty l b_{n} q^{l n} &\text{ for } l \nmid f,\\
 \sum\limits_{n=1}^\infty b_{l n}q^n &\text{ for } l \mid f.
 \end{cases}$$

If $l \nmid f$, then by proceeding as in~\cite[p. 853]{MR1376558}
to obtain
$$b_{l n} = \sum_{bc= l n}  \varphi(c)\varphi^{-1}(b)b = \varphi(l) \sum\limits_{bc' = n, l \nmid b}  \varphi^{-1}(b) \varphi(c') +l \varphi^{-1}(l) \sum_{b'c =n}  \varphi^{-1}(b') \varphi(c)b' ,$$
$$l b_{n} = l \sum_{bc=n}  \varphi^{-1}(b) \varphi(c)b = \varphi(l) \sum_{bc=n} \varphi^{-1}(bl) \varphi(c)l b = \varphi(l) \sum_{b'c= l n , l \mid b'} \varphi^{-1}(b') \varphi(c)b'.$$

If $l \mid f$, then $b_{l n}=0$ and the result follows.
\end{proof}

\subsubsection{Ordinary and critical refinement of $E_\varphi$}
\label{Ordcrit}

Recall the ordinary and critical refinement of the Eisenstein series following \cite[\S 1.3]{MR2929082} and \cite{MR3357177}. The Eisenstein series $E_{\varphi}$, as defined above, is a new form of level $\Ga_0(f^2)$. If $l$ is a prime that does not divide $f$, then $E_\varphi$ has $T_l$ eigenvalue $\varphi(l) + l \varphi^{-1}(l)$. The polynomial 
$$  X^2-(\varphi(l)+l \varphi^{-1}(l)) X+l  $$
 has two roots $\alpha:=\varphi(l)$ and $\beta:=l \varphi^{-1}(l)$. 
 Note that $\nu_l(\alpha)=0$ and $\nu_l(\beta)=1=2-1$. 

For any prime $l \nmid f$, consider the matrix  $\gamma_l=\left(\begin{smallmatrix}
l & 0\\
0 & 1\\
\end{smallmatrix}\right)$.  For $g \in M_2(\Ga_0(f^2))$,  define
 \begin{align*}
 [l]_{\varphi}^{+}g(z):= g(z)-\varphi(l)g(l z) =g\mid_{(1-\frac{\varphi(l)}{l}[\gamma_l])} (z),
 \end{align*}
 and 
 \begin{align*}
 [l]_{\varphi}^{-}g(z):= g(z) -l \varphi^{-1}(l)g(l z) =g\mid_{(1-\varphi^{-1}(l)[\gamma_l])}(z). 
 \end{align*}
 
Note that $[l]_\varphi^\pm$ commutes with Hecke operators $T_n$ if $(n,l)=1$, and hence they preserve the eigenspaces of these Hecke operators. Now $[l]_{\varphi}^{+}E_\varphi$ is an Eisenstein series of level $\Ga_0(f^2l)$, it is an eigenvector for the $U_l$ operator with $U_l([l]_\varphi^+ E_\varphi) = (l \varphi^{-1}(l)) [l]^+_\varphi E_\varphi$, and hence $[l]^+_\varphi(E_\varphi)$ is an eigenform in $E_2(\Ga_0(f^2l))$. In fact, $[l]_\varphi^+$ acts as a critical refinement on $E_{\varphi}$.  Similarly, $[l]^-_\varphi$ acts as an ordinary refinement on $E_{\varphi}$. Since $[l]_\varphi^- \circ [l]_\varphi^+$ commutes with $T_n$ for all $(n,l)=1$, $[l]_\varphi^- \circ [l]_\varphi^+$ preserves the eigenspaces of Hecke operators $T_n$ if $(n,l)=1$. An easy computation shows that $U_l([l]_\varphi^- \circ [l]_\varphi^+ E_\varphi) =0$, hence $[l]_\varphi^- \circ [l]_\varphi^+ E_\varphi$ is an eigenform in $E_2(\Ga_0(f^2l^2))$. 

Let $N=N_1N_2$ be the decomposition of positive integers $N$ with $N_1, N_2$ as in equation~\ref{decomp}. Assume that $f^2 \mid N_1$. Let $t$ be a square-free natural number such that $t^2 \mid \frac{N_1}{f^2}$,  $M$ be a square-free natural number coprime to $ft$, and $L$ be a square-free natural number coprime to $Mft$ such that $ML \mid N$. Define 
$$
E_{\varphi,Mt,Lt}:= \prod_{l \mid Mt} [l]_{\varphi}^{+} \circ \prod_{q \mid Lt} [q]_{\varphi}^{-} E_{\varphi} \in E_2(\Ga_0(f^2t^2ML),\CC) \subseteq E_2(\Ga_0(N),\CC). $$
The Eisenstein series $E_{\varphi, Mt , Lt}$ is an eigenform for all Hecke operators, and the eigenvalues are given by
\begin{align*}
T_{r}(E_{\varphi, Mt, Lt}) &=  (\varphi(r)+ r \varphi^{-1}(r)) E_{\varphi, Mt,Lt} \quad  \text{for}\ r \nmid N, \\
U_{l}(E_{\varphi, Mt, Lt}) &= l \varphi^{-1}(l) E_{\varphi, Mt, Lt} 
\qquad \qquad \quad \text{for}\  l \mid M, \\
U_{q}(E_{\varphi, Mt, Lt}) &=  \varphi(q) E_{\varphi, Mt, Lt}  
\qquad \qquad \qquad  \text{for}\  q \mid L ,  \\
U_{p}(E_{\varphi, Mt, Lt}) &=  U_{p}(E_\varphi) =0  
\qquad \qquad \qquad \quad \text{for}\   p \mid ft.  
\end{align*}

\begin{proposition}
\label{lvalue}
Let $E_{\varphi, Mt,Lt} \in E_2(\Ga_0(N), \CC)$ be the Eisenstein series as defined above, and $\chi$ be a Dirichlet character whose conductor is coprime to $N$. The twisted $L$-function  $L(E_{\varphi, Mt,Lt}, \chi, s)$ is given by 
$$ L(E_{\varphi,Mt,Lt}, \chi, s) = \Big( \prod_{l \mid Mt} (1 - \chi(l)\varphi(l)l^{-s}) \Big) \Big( \prod_{q \mid Lt} (1- \chi(q)\varphi^{-1}(q) q^{1-s}) \Big) L( \chi\varphi, s) L( \chi\varphi^{-1}, s-1). $$
\end{proposition}
\begin{proof}
For any modular form $f$, we have $L(f\mid_{\gamma_d}, \chi, s) = \chi(d) d^{1-s}L(f, \chi, s)$, where $\gamma_d = \psmat{d}{0}{0}{1}$. This essentially follows from the fact that if the Fourier expansion of $f(z)=\sum\limits_{k \geq 0}^{} a_k q^k$ then the Fourier expansion of $f\mid_{\gamma_d}(z)=\sum\limits_{n \geq 0 }^{} b_n q^n$, where $b_n = \begin{cases} a_{dk} \quad \text{ if } n=dk,\\
0 \qquad \text{ otherwise.}
\end{cases}$
Using this, the $L$-functions of the critical and ordinary refinements of $f$ are calculated as follows:
 $$L( [l]_\varphi^+ E_{\varphi}, \chi,  s) = L( E_\varphi, \chi, s) - \frac{\varphi(l)}{l}L( E_{\varphi} \mid_{\gamma_l}, \chi, s) = (1- \chi(l)\varphi(l)l^{-s}) L(E_{\varphi}, \chi, s),$$
and 
 $$L([q]_\varphi^- E_{\varphi}, \chi,  s ) = L(E_\varphi, \chi, s) - \varphi^{-1}(q)L(E_{\varphi} \mid_{\gamma_q}, \chi, s) = (1- \chi(q) \varphi^{-1}(q)q^{1-s}) L(E_{\varphi}, \chi,  s).$$
 Combining both of these, we get
 $$L( [l]_\varphi^+ \circ [q]_\varphi^-E_{\varphi}, \chi,  s) = L([q]_\varphi^- E_\varphi, \chi, s) - \frac{\varphi(l)}{l}L( ([q]_\varphi^-E_{\varphi}) \mid_{\gamma_l}, \chi, s) = (1- \chi(l)\varphi(l)l^{-s}) L([q]_\varphi^-E_{\varphi}, \chi, s).$$
 Since $L(E_\varphi, \chi, s)=L(\chi\varphi, s)L(\chi\varphi^{-1}, s-1)$, the proof of the proposition is complete by iterating these identities over all primes dividing $Mt$ and $Lt$.
\end{proof}

\section{Classification of non-rational Eisenstein maximal ideals}
\label{maximalideals}
Let $M^B_2(\Ga_0(N),\Z)$ denote the subset of elements $ M_2(\Ga_0(N), \CC)$ whose Fourier coefficients at the cusp 
 $\infty$  belong to  $\Z$.
 By  \cite[Cor. 12.3.12, Proposition 12.4.1]{MR1357209},  $M^B_2(\Ga_0(N),\Z)$ is stable under the action of Hecke operators and contains a basis of $M_2(\Ga_0(N), \CC)$.
 For any $\Z$-algebra  $S$, we define $M^B_2(\Ga_0(N), S) := M^B_2(\Ga_0(N), \Z) \otimes_\Z S$. We denote by $S^B_2(\Ga_0(N), S)$ the corresponding space of cusp forms (the forms for which the constant terms at all cusps are zero).

Let $\T:= \T(N) \subseteq \mathrm{End}_{\Z}(S_2(\Ga_0(N), \Z)$
be the Hecke algebra generated by Hecke operators $T_n$ for all $n \in \N$. It is well known that $\T$ is a free $\Z$-module of rank $d$, where $d= \mathrm{dim}_\CC S_2(\Ga_0(N), \CC)$. For any $\Z$ module $M$, denote
$\T_M:=\T_\Z \otimes M$. For any maximal ideal $\m$ of $\T$, let $\kappa(\m):=\T/\m$ be its residue field. This is a finite field of characteristic $\ell$. By \cite[Proposition 5.1]{MR1047143}, there exists a unique continuous semi-simple Galois representation associated with $\m$ such that
$$
\rho_{\m}: \Gal(\overline{{\Q}}/{\Q}) \rightarrow \GL_2(\T/\m)
$$
it is unramified outside $\ell N$, and for all primes $p \nmid \ell N$, we have $\mrm{Tr}(\rho_\m(\Frob_p)) \equiv T_p \pmod \m$ and $\mrm{det}(\rho_\m(\Frob_p)) \equiv  p \pmod \m$. By Chebotarev's density theorem, we have $ \mrm{det}(\rho_\m)=\overline{\chi}_\ell$, where $\overline{\chi}_\ell$ is the mod $\ell$ cyclotomic character.

\begin{definition}(Cusp form associated with the maximal ideal $\m$)
\label{modlcuspform}
 There is a nondegenerate bilinear pairing \cite[(10), p. 465]{MR1047143}: 
$$ S^B_2(\Ga_0(N), \kappa(\m)) \times  \left(\T \otimes_\Z\kappa(\m)\right)  \rightarrow \kappa(\m). $$
The natural ring homomorphism $T_{\m}: \T \rightarrow  \kappa(\m)$ gives rise to a cusp form $f_{\m}:= \sum_{n \ge 1} T_\m(T_n)q^n \in S_2(\Ga_0(N), \kappa(\m))$. We call $f_\m$  {\it the cusp form associated with the Eisenstein maximal ideal $\m$}.    The mod $\ell$ modular form  $f_{\m}$ may lift to several modular forms of different weights, levels, and nebentypus~\cite[Theorem 2.1]{MR3512875}. 
\end{definition}

 Recall that a maximal ideal $\m$ in $\T$ is called an Eisenstein maximal ideal if the associated Galois representation $\rho_\m: \Gal(\overline{{\Q}}/{\Q}) \rightarrow \GL_2(\kappa(\m))$, associated with $\m$, is reducible. Since $\rho_\m$ is semi-simple, 
there exist continuous characters
$\bar{\alpha}_{\m}, \bar{\beta}_{\m}: \Gal(\overline{{\Q}}/{\Q}) \rightarrow \kappa(\m)^{\times}$
such that $\rho_{\m} \simeq \bar{\alpha}_{\m} \bigoplus \bar{\beta}_{\m}$.
We now classify the Eisenstein maximal ideals of $\T$ in the following proposition, which generalizes ~\cite[Lemma 0.1]{EisenJ} (see also Ren~\cite{MR4621854}).

\begin{proposition}
\label{Eisensteincriteria}
Let $\m$ be an Eisenstein maximal ideal of $\T$ with residual characteristic $\ell$ such that $\ell \nmid 6N_1$.
Then there exists a unique character $\overline{\epsilon}_\m : \Gal(\overline{{\Q}}/{\Q}) \rightarrow \kappa(\m)^{\times} $ unramified outside $N_1$
such that
$$\rho_\m \simeq \overline{\epsilon}_\m \oplus (\overline{\epsilon}_\m)^{-1} \overline{\chi}_\ell.$$
\end{proposition}
\begin{proof}
 If $\overline{\chi}: \Gal(\overline{{\Q}}/{\Q}) \to \kappa(\m)^\times$ is a continuous character, by \cite[Prop 21.6.3, p 223] {MR1860042}, then $\overline{\chi} = \overline{\theta} \overline{\chi}_\ell^i$ for some integer $i$ and some   character $\overline{\theta}$ unramified at $\ell$. Since $\rho_\m = \bar{\alpha}_\m \oplus \bar{\beta}_\m$ for some continuous character $\bar{\alpha}_\m, \bar{\beta}_\m : \Gal(\overline{{\Q}}/{\Q}) \to \kappa(\m)^\times$, there exist characters $\overline{\epsilon}_1,\overline{ \epsilon}_2: \Gal(\overline{{\Q}}/{\Q}) \to \kappa(\m)^\times$ unramified at $\ell$ and some integers $i$ and $j$ such that 
$$
\rho_{\m} \simeq  \left(\begin{smallmatrix}
\overline{\epsilon}_1\overline{\chi}_\ell^i& 0\\
0& \overline{\epsilon}_2\overline{\chi}_\ell^j \\
\end{smallmatrix}\right).
$$
Since $\overline{\chi}_\ell  = \det(\rho_\m)  = \overline{\epsilon}_1 \overline{\epsilon}_2 \overline{\chi}_\ell^{i+j}$, the congruences $i+j \equiv 1 \pmod{\ell-1}$ and $\overline{\epsilon}_1 = (\overline{\epsilon}_2)^{-1}$ hold. Define $\overline{\epsilon}_\m:= \overline{\epsilon}_2$, which is unramified at $\ell$.
   
        Let $f_{\m} \in S_2(\Ga_0(N), \kappa(\m))$ be the cusp form associated with the Eisenstein maximal ideal $\m$ (as in Definition \ref{modlcuspform}). Let $\rho_{f_\m}$ be the unique two-dimensional semisimple odd continuous Galois representation 
  $ \rho_{f_{\m}}: \Gal(\overline{{\Q}}/{\Q}) \to \GL_2(\overline{\mathbb{F}}_\ell)$ associated to $f_\m$. 
  Note that the representations $\rho_{f_\m}$ and $\rho_{\m}$ have the same characteristic polynomial at all $\Frob_r$ for primes $r \nmid \ell N$, and hence by the Brauer-Nesbitt theorem, for any prime $q$, we have
  $(\rho_{f_{\m}})^{ss}\mid_{G_q} \simeq \rho_{\m} \mid_{G_q}$, where $G_q$ denotes the decomposition group at $q$.  Since $N_2$ is square-free, the representation $\rho_{f_{\m}}$ is semi-stable outside $\ell N_1$ (\cite[Proposition 14.1, p. 113]{MR488287} or \cite[Lemma 0.1]{EisenJ}) and therefore the characters $\overline{\epsilon}_\m\overline{\chi}_\ell^j$, $(\overline{\epsilon}_\m)^{-1}\overline{\chi}_\ell^i$ are unramified outside $\ell N_1$. Therefore, $\overline{\epsilon}_\m$ is unramified outside $N_1$ as $\ell \nmid N_1$.
  Let $a_\ell$ be the $T_\ell$-eigenvalue of $f_\m$.

 \begin{itemize}
  \item  If $\ell \nmid N_2$ then $\rho_{f_{\m}\mid_{G_\ell}}$ is 
reducible, and we get $a_\ell \neq 0$; otherwise, it contradicts Fontaine's theorem~\cite[Theorem 2.6]{MR1176206}. Hence, by Deligne's theorem~\cite[Theorem 2.5]{MR1176206}, we get
   $$ \rho_{f_{\m}\mid_{G_\ell}} =  \left(\begin{smallmatrix}
\lambda(a_\ell^{-1})\overline{\chi}_\ell & \star\\
0& \lambda(a_\ell) \\
\end{smallmatrix}\right), $$
where $\lambda(x)$ denotes the unramified character that sends $\Frob_\ell$ to $x$.
\item If $\ell \mid N_2$, we cannot directly apply~\cite[Theorem 2.5, Theorem 2.6]{MR1176206}. We use the following Ribet trick as described in \cite{EisenJ}. By \cite[Theorem 2.1, Step 4]{MR1265566}, we replace the mod $\ell$ cusp form $f_\m \in S_2(\Ga_0(N), \kappa(\m))$ with some $g_\m \in S_{\ell+1}(\Ga_0(N/\ell), \kappa(\m))$ such that $f_\m$ and $g_\m$ have the same mod $\ell$ Galois representation. Since $\ell \nmid N/\ell$, by~\cite[Theorem 2.5, Theorem 2.6]{MR1176206}, we obtain
   $$
    \rho_{g_{\m}\mid_{G_\ell}} =  \left(\begin{smallmatrix}
\lambda(a_\ell^{-1})\overline{\chi}_\ell^l & \star\\
0& \lambda(a_\ell) \\
\end{smallmatrix}\right).  $$
\end{itemize}
As a result, in either case, we may assume that $ \overline{\alpha}_{\m}=\overline{\epsilon}_\m^{-1} \overline{\chi}_\ell$ and  $\overline{\beta}_{\m}=\overline{\epsilon}_\m$.
  \end{proof}

 \begin{definition}[Character associated to Eisenstein maximal ideals]
 \label{CharEisen}
 Let $\m$ be a maximal ideal of $\T = \T(N)$, and let the residual characteristic of $\m$ be coprime to $6N_1$. The unique character $\overline{\epsilon}_\m$ (as in Prop. \ref{Eisensteincriteria}) is called the Dirichlet character associated with the Eisenstein maximal ideal $\m$.  As in~\cite[Prop 21.6.3, p. 223-224] {MR1860042}, we can consider $\overline{\epsilon}_\m$ as a Dirichlet character $(\Z/N_1\Z)^\times \to \kappa(\m)^\times$. 
An Eisenstein maximal ideal $\m$ is {\bf rational} (resp., {\bf non-rational}) if the associated Dirichlet character $\overline{\epsilon}_\m$ is trivial (resp., non-trivial).  
\end{definition}

 \begin{cor}\label{Trmodm}
  Let  $\m$ be an Eisenstein maximal ideal of $\T = \T(N)$
  with the residual characteristic $\ell$ such that $\ell \nmid 6N_1$.
  For any prime $r \nmid N$, we have
  $$T_r \equiv \overline{\epsilon}_\m(r) + r \cdot (\overline{\epsilon}_\m)^{-1}(r) \pmod{\m}.$$ 
  \end{cor}

  \begin{proof}
  The proof of this corollary lies in the proof of Proposition \ref{Eisensteincriteria}. If $r \nmid \ell N$, then the corollary follows from $\rho_{f_\m}$ and $\rho_{\m}$ having the same trace at $\Frob_r$. If $r=\ell \nmid 6N_1$, then $T_\ell \pmod{\m} = a_\ell = \overline{\epsilon}_m (\ell) =  \overline{\epsilon}_\m(\ell) + \ell (\overline{\epsilon}_\m)^{-1}(\ell)$ as shown in~\cite[Lemma 1.1]{EisenJ} 
  (which is a private note).  
  \end{proof}
Let $\m$ be a non-rational Eisenstein maximal ideal of level $N$ with residual characteristic $\ell \nmid 6N_1$. 
For a prime $r \nmid N$, let $\overline{f_r(X)} \in \F_\ell[X]$ be the minimal polynomial of $\overline{\epsilon}_\m(r) + r (\overline{\epsilon}_\m)^{-1}(r) \in \kappa(\m)$ over $\F_\ell$. Note that the polynomial $\overline{f_r(X)}$ does not depend on the choice of isomorphism $\psi'$. Let $f_r(X) \in \Z[X]$ be a monic polynomial (necessarily irreducible) such that $f_r(X) \equiv \overline{f_r(X)} \pmod{\ell}$.

\begin{cor}\label{thirdequivalence}
Let $\m$ be a non-rational Eisenstein maximal ideal of $N$ with
the residual characteristic $\ell$ such that $\ell \nmid 6N_1$. If 
$I_{\overline{\epsilon}_\m}(N) = (f_r(T_r), r \text{ primes such that } r \nmid N ),$
then $\m \supseteq (\ell, I_{\overline{\epsilon}_\m}(N) )$.
\end{cor}
\begin{proof}
The ideal $(\ell, I_{\overline{\epsilon}_\m}(N) )$ is independent of the choice of $f_r(X)$ as any lift of $\overline{f_r(X)}$ satisfies $f_r(X) \equiv \overline{f_r(X)} \pmod{\ell}$. Consider the quotient map
$$\T(N) \to \T(N)/\m \cong \kappa(\m) \text{ sending } T_r \mapsto \overline{\epsilon}_\m(r) + r (\overline{\epsilon}_\m)^{-1}(r) \text{ for primes } r \nmid N.$$
By Corollary~\ref{Trmodm}, the ideal $(\ell, I_{\overline{\epsilon}_\m}(N) )$ is in the kernel of the map (for any choice of isomorphism between $\T/\m$ and $\kappa(\m)$). Hence,
we get $\m \supseteq (\ell, I_{\overline{\epsilon}_\m}(N) )$.
\end{proof}
We now recall the definition of the $q$-old and $q$-new Eisenstein ideals ( ~\cite[\S 2]{MR482230},~\cite[Section 2]{MR3988582}, \cite[\S 2.2]{MR4621854}). For any prime $q \mid N$, there are two degeneracy coverings $\alpha_q(N/q), \beta_q(N/q): X_0(N) \to X_0(N/q)$ given by $\alpha_q(N)(z) = z \pmod{\Ga_0(N/q)}$ and $\beta_q(N)(z) = qz \pmod{\Ga_0(N/q)}$; here $z \in \mathbb{H} \cup \mathbb{P}^1(\Q)$, and we identify $X_0(N)$ (resp. $X_0(N/q)$) with $\Ga_0(N) \backslash \mathbb{H} \cup \mathbb{P}^1(\Q)$ (resp. $\Ga_0(N/q) \backslash \mathbb{H} \cup \mathbb{P}^1(\Q)$). 
These two degeneracy maps, in turn, define two maps $\alpha_q(N/q)*, \beta_q(N/q)* : J_0(N/q) \to J_0(N)$.

Consider the map $\gamma_q(N)^*: J_0(N/q) \times J_0(N/q) \to J_0(N)$ given by 
\[
\gamma_q(N)^* (x,y) = \alpha_q(N/q)^*(x) + \beta_{q}(N/q)^*(y).
\]
The image of $\gamma_q(N)^*$ is called the {\it $q$-old subvariety} of $J_0(N)$, denoted by $J_0(N)^{q-\old}$. The quotient of $J_0(N)$ by $J_0(N)^{q-\old}$ is called the {\it $q$-new subvariety} of $J_0(N)$, denoted by $J_0(N)^{q-\new}$. 
The image of $\T(N)$ in $\mathrm{End}(J_0(N)^{q-\old})$ (resp. in $\mathrm{End}(J_0(N)^{q-\new})$) is denoted by $\T(N)^{q-\old}$ (resp. $\T(N)^{q-\new}$). 

 For $q \mid N$, the $q$-new subvariety of $J_0(N)$ is usually defined as the connected component of the kernel of the degeneracy map \cite{MR1159117} $J_0(N) \longrightarrow J_0(N/q) \times J_0(N/q)$. In contrast, we define the $q$-new quotient as the quotient $J_{q\text{-new}} := J_0(N)/J_{q\text{-old}}$. These two constructions are dual to each other~\cite[p.3]{MR1085264}. 

A maximal ideal of $\T(N)$ is called $q$-old (resp. $q$-new) if its image in $\T(N)^{q-\old}$ (resp. $\T(N)^{q-\new}$) is still maximal. Thus, any maximal ideal $\m$ of $\T(N)$ is either $q$-old, $q$-new, or both.

\begin{remark}
\label{old-new}
For a prime $r \mid N$, let $U_r$ (resp. $\tau_r$) denote the Hecke operator at $r$ of level $N$ (resp. $N/r$).
Let $p$ be a prime that divides $N_1$, and $q$ be a prime that divides $N_2$.  Additionally, $w_{q}$ denotes the Atkin-Lehner involution. By~\cite[(2.7), (2.10), p. 2439]{MR3988582}, it follows that 
\begin{enumerate}
\item
$U_{p} = 0 \in \T(N)^{p-\new}, U_{p}^2 -\tau_{p} U_{p}=0 \in \T(N)^{p-\old}  \text{ if $r \geq 2 $ }$
\item
$U_{q}+ w_{q} = 0 \in \T(N)^{q-\new},  U_{q}^2 -\tau_{q} U_{q} + q = 0 \in \T(N)^{q-\old} \text{ if $r=1 $ }$.
\end{enumerate}
\end{remark}

For a prime $r \mid N$, let $\overline{g_{r}(X)} \in \F_\ell[X]$ (resp. $\overline{h_{r}(X)} \in \F_\ell[X]$) be the minimal polynomial of $ r (\overline{\epsilon}_\m)^{-1}(r) \in \kappa(\m)$ (resp. $ \overline{\epsilon}_\m(r) \in \kappa(\m)$) over $\F_\ell$. Let $g_{r}(X) \in \Z[X]$ (resp. $h_{r}(X) \in \Z[X]$ ) be a monic (necessarily irreducible) polynomial such that $g_{r}(X) \equiv \overline{g_{r}(X)} \pmod{\ell}$ (resp. $h_{r}(X) \equiv \overline{h_{r}(X)} \pmod{\ell}$).

\begin{lemma}\label{valueatpi}
Let $\m$ be a non-rational Eisenstein maximal ideal of $\T(N)$ with residual characteristic $\ell$ coprime to $6N_1$, and let $q$ be a prime such that $q \mid N_2$.  Then $$U_{q} \equiv  q (\overline{\epsilon}_\m)^{-1}(q) \text{ or } \overline{\epsilon}_\m(q) \pmod{\m}.$$
In particular, either $g_{q}(U_{q}) \in \m$ or $h_{q}(U_{q}) \in \m$.
\end{lemma}

\begin{proof}
We proceed with the proof as in \cite[Lemma 2.1]{MR3540618}.
Suppose that $\m$ is $q$-old. Let $R$ be the common subring of $\T(N/q)$ and $\T(N)^{q-\old}$. 
Let $\mathfrak{n}$ be the corresponding maximal ideal of $\T(N/q)$ to $\m$. Let $\tau_{q}$ denote the $q$-th Hecke operator in $\T(N/q)$. Then
\begin{equation}
\T(N/q) = R[\tau_{q}], \qquad \T(N)^{q-\old}=R[U_{q}]\  \mathrm{and}\ \T(N)/\m \cong \T(N/q)/\mathfrak{n} .
\end{equation}
Note that $\mathfrak{n}$ is an Eisenstein ideal of $\T(N/q)$, it follows that $I_{\m}(N/p) \subseteq \mathfrak{n}$. 
Since $q \mid\mid N$, we have $q \nmid N/q$. Consequently, the operator $\tau_{q} - \epsilon_\m(q) - q (\epsilon_\m(q))^{-1} \in \mathfrak{n}$. Moreover, by Remark~\ref{old-new}, the Hecke operator $U_{q}^2 -\tau_{q} U_{q} + q = 0$. Therefore, in the ring $\T(N)/\m \cong \T(N/q)/\mathfrak{n}$,  we have $U_{q}^2 - (\overline{\epsilon}_\m(q) + q (\overline{\epsilon}_\m)^{-1}(q) ) U_{q} + q = (U_{q} - \overline{\epsilon}_\m(q))(U_{q} - q (\overline{\epsilon}_\m)^{-1}(q))=0 $.

Now suppose that $\m$ is $q$-new. If $\ell= q \mid N_2$ and $\ell \nmid 6N_1$, then the proof of Corollary \ref{Trmodm} implies that $U_{q} \equiv \overline{\epsilon}_\m(q) \pmod{\m}$. We may therefore assume that $\ell \nmid 6N$. Since $U_{q} + w_{q}=0 \in \T^{q-new}$, we get $U_{q} - a_{q} \in \m$, where $a_{q} \in \{ 1, -1 \}$. Now, the proof reduces to showing that 
$$a_{q} \equiv \overline{\epsilon}_\m(q)\ \mathrm{or}\ q (\overline{\epsilon}_\m)^{-1}(q) \pmod{\m}.$$

Since the residual characteristic $\ell$ of $\m$ is odd with $\ell \neq q$ and $q \mid \mid N$ , by~\cite[Theorem 3.1(e)]{DDT}, it follows that $(\rho_{f_\m}\mid_{G_{q}})^{ss}= \lambda(a_{q}) \overline{\chi}_\ell \oplus \lambda(a_{q})$, here $\lambda(a_{q})$ is the unramified quadratic character which sends $\Frob_{q}$ to $a_{q}$. Since $\rho_\m =  (\overline{\epsilon}_\m)^{-1} \overline{\chi}_\ell \oplus  \overline{\epsilon}_\m$, comparing the trace at $\Frob_{q}$, we get
$$ a_{q} + q a_{q} \equiv \overline{\epsilon}_\m(q) + q (\overline{\epsilon}_\m)^{-1}(q) \pmod{\m}.$$
Using the fact $a_{q} = a_{q}^{-1}$, along with some minor algebraic manipulation, we obtain
$$ (a_{q} -  \overline{\epsilon}_\m(q)) ( 1 - a_{q}^{-1} q (\overline{\epsilon}_\m)^{-1}(q)) \equiv 0 \pmod{\m}. $$
\end{proof}


\begin{lemma}\label{valueatp}
Let $\m$ be a non-rational Eisenstein maximal ideal of $\T(N)$ with residual characteristic $\ell$ that is coprime to $6N_1$. Let $p$ be a prime that divides $N_1$. 
\begin{enumerate}
\item If $\m$ is $p$-new, then $U_p \in \m$.
\item If $\m$ is $p$-old, then $U_p \equiv 0 \text{ or } p (\overline{\epsilon}_\m)^{-1}(p) \text{ or } \overline{\epsilon}_\m(p) \pmod{\m}$. In particular, either $U_p \in \m$ or $g_{p}(U_p) \in \m$ or $h_{p}(U_p) \in \m$.
\item If $p \mid \mrm{cond}(\overline{\epsilon}_\m)$, then $U_p \in \m$.
\end{enumerate}
\end{lemma}

\begin{proof}
If $\m$ is $p$-new, then by Remark~\ref{old-new}, we have
$U_p \in \m$ as $U_p=0 \in \T(N)^{p-\new}$.

Suppose that $\m$ is $p$-old and $U_p \notin \m$. Since $U_p(U_p - \tau_p)=0 \in \T(N)^{p-\old}$, there exists a maximal ideal $\mathfrak{n}$ in $\T(N/p)$ corresponding to $\m$~\cite[Lemma 2.3 (2)]{MR3988582} such that $U_p \equiv \alpha \pmod{\m}$ if and only if $\tau_p \equiv \alpha \pmod{\mathfrak{n}}$.
We inductively continue this level lowering until $v_p(N)=1$. Once we have $v_p(N)=1$, we use Lemma \ref{valueatpi} and deduce that $\tau_p \equiv p (\overline{\epsilon}_\m)^{-1}(p) \text{ or } \overline{\epsilon}_\m(p) \pmod{\mathfrak{n}}$, which implies part $(2)$ of the Lemma. 

Finally, suppose that $p \mid \mrm{cond}(\overline{\epsilon}_\m)$ and $U_p
\not \in \m$. Then, we obtain $U_p \equiv p (\overline{\epsilon}_\m)^{-1}(p) \text{ or } \overline{\epsilon}_\m(p) \pmod{\m}$, which is a contradiction to the assumption.
\end{proof}

We collate all the information obtained so far as a Theorem, which is a classification of Eisenstein maximal ideals. 

\begin{theorem}(Classification of Eisenstein maximal ideals)\label{classification}
 Given an integer $N$, write $N=N_1N_2$ as in equation \ref{decomp} and let $\m$ be a non-rational Eisenstein maximal ideal of $\T(N)$ whose residual characteristic $\ell$ is coprime to $6N_1$.
The ideal $\m$ determines a non-trivial Dirichlet character
$\overline{\epsilon}_\m$ of the conductor $f$ (with $f^2 \mid N_1$). Moreover, there exist square-free natural numbers $t$ and $M$ such that
  $$T_r \equiv \overline{\epsilon}_\m(r) + r \cdot (\overline{\epsilon}_\m)^{-1}(r) \pmod{\m} \text{ \qquad\qquad for primes } r \nmid N,$$ 
    $$U_p \equiv 0 \pmod \m \text{ \qquad\qquad for primes } p \mid ft \text{ with }  t \mid N_1 \text{ and }  (t,f)=1,$$ 
    $$U_s \equiv  s \cdot (\overline{\epsilon}_\m)^{-1}(s) \pmod{\m} \text{ \qquad\qquad for primes } s \mid M \text{ with } M \mid N \text{ and } (M,ft)=1,$$
        $$U_q \equiv  \overline{\epsilon}_\m(q) \pmod{\m} \text{ \qquad\qquad for all other primes } q.$$ 
In particular, the maximal ideal $\m$ contains 
$$I_{\overline{\epsilon}_\m,M, t}(N)= \Big(I_{\overline{\epsilon}_\m}(N), U_p \text{ for primes } p \mid ft ,  g_s(U_s) \text{ for primes } s \mid M, h_q(U_q)   \text{ for all other primes } \Big), $$ 
and therefore $\m = (\ell, I_{\overline{\epsilon}_\m,M, t}(N))$ and $\kappa(\m) \cong \T/\m \cong \F_\ell[\overline{\epsilon}_\m]$.
\end{theorem}

\begin{proof}
The proof of this theorem follows from Corollary \ref{thirdequivalence}, Lemma \ref{valueatp}, and Lemma \ref{valueatpi}.
\end{proof}

\begin{remark}
\label{classificationremark}
For any prime $q \mid N$ such that $(q,Mft)=1$ and $q \equiv \overline{\epsilon}_\m^2(q)  \pmod \ell$, we have 
\[
I_{\overline{\epsilon}_\m,M, t}(N) = I_{\overline{\epsilon}_\m,Mq, t}(N).
\]
Thus, from now on, we assume that if a prime $p' \neq \ell$ divides $N$ and satisfies $p' \equiv \overline{\epsilon}_\m^2(p') \pmod{\ell}$, then $p' \mid M$. Without loss of generality, we assume that
\begin{itemize}
\item 
the integer $M$ is the smallest integer that meets these properties. This uniquely determines $M$. 
\item
By definition, $t$ is a square-free natural number with the property that $(f,t)=1$ and a prime $p \mid t$ if and only if $U_p \in \m$. This uniquely determines $t$.
\end{itemize}
As a consequence, we associate a unique ideal $I_{\overline{\epsilon}_\m,M, t}(N)$ with $\m$.
\end{remark}
\section{Boundary of Eisenstein series \texorpdfstring{$E_{\varphi, M,L}$}{Eϕ, M, L}}
We compute the boundary of the non-rational Eisenstein series $E_{\varphi, M,L}$ in this section.
\subsection{Cuspidal divisor groups}
For any congruence subgroup $\Ga \subset \mathrm{SL}_2(\Z)$, let $\cX_{\Ga}$ be the model of the modular curve defined over $\Q$ \cite[p. 521]{MR800251} and let $\partial(X_{\Ga}):=\Ga \backslash \sP^1(\Q) \subset \cX_{\Ga}(\overline{\Q})$ be the set of all cusps of $\Ga$. In particular,  $\partial(X_0(N))$ (resp., $\partial(X(N))$) is the set of all cusps for the congruence subgroup $\Ga_0(N)$  (resp., $\Ga(N)$).

Let $\D_{N}:=\Div^0(\partial(X_0(N));\Z)$ be the group of divisors supported on the cusps $\partial(X_0(N))$. 
 For any ring $R$, consider the $R$-module $\Div^0(\partial(X_0(N));R) =\Div^0(\partial(X_0(N));\Z) \otimes_{\Z} R$. In particular, we have an isomorphism $\Div^0(\partial(X_0(N));\CC):=\Div^0(\partial(X_0(N));\Z) \otimes_\Z \CC$.
Recall that the  Galois group $G_{\Q}:=\Gal(\overline{\Q}/\Q)$ and the Hecke algebra $\T$  act on $\D_{N}$ \cite[pp.~521--522]{MR800251}.

For the congruence subgroup $\Ga(N)$, following Shimura's notation, we denote by $\bstwomat{a}{b}_{\Ga(N) } \in \partial(X(N)):=\Ga(N) \backslash  \sP^1(\Q) $ the cusp represented by $\frac{\alpha}{\beta} \in \Q$ with integers $\alpha,\beta$  such that $(\alpha, \beta)=1$, $\alpha \equiv a \pmod N$, and $\beta \equiv b \pmod N$. The group $\GL_2(\Z/N\Z)$ acts on $\partial(X(N))$ by matrix multiplication \cite[p. 539]{MR800251}.  The natural map $\Ga_0(N) \rightarrow \GL_2(\Z/N\Z)$ determines an action of $\Ga_0(N)$ on $ \partial(X(N))$.

For any congruence subgroup $\Ga$,  we denote a cusp in $\partial(X_{\Ga})$ by $\bstwomat{a}{b}_{\Ga}$.  Now, the cusps of $\partial(X_0(N)):=\Ga_0(N) \backslash \sP^1(\Q)$ are in bijection with a pair of coprime integers $(a,b)$ represented by $\bstwomat{a}{b}_{\Ga_0(N)}$.  These cusps are the $\Ga_0(N)$ orbits of $\bstwomat{a}{b}_{\Ga(N)}$ with the relations: 
\begin{enumerate}
\label{cuspGa0}
   \item $\bstwomat{a}{b}_{\Ga_0(N)} = \bstwomat{a^\prime}{b^\prime}_         
         {\Ga_0(N)}$ if $a \equiv a^\prime, b \equiv b^\prime \pmod N$. 
   \item $\bstwomat{a}{b}_{\Ga_0(N)} = \bstwomat{a^\prime}{b}_         
         {\Ga_0(N)}$ if $a \equiv a^\prime \pmod b$. 
   \item $\bstwomat{ra}{b}_{\Ga_0(N)} = \bstwomat{a}{rb}_         
         {\Ga_0(N)}$,
         for all $r \in \Z$ with $(r, N)=1$.             \end{enumerate}

For congruence subgroups of the form $\Ga_0(N)$ with $N\in \N$, the integer $d=\gcd(b,N)$ depends only on the cusp $\bstwomat{a}{b}_{\Ga_0(N)}$. Following the terminology of Stevens, we call $d$ the divisor of $\bstwomat{a}{b}_{\Ga_0(N)}$. The ramification indices of the cusp $x = \bstwomat{a}{b}_{\Ga_0(N)} \in \partial(X_0(N))$ over $X(1):=\SL_2(\Z) \backslash \bigg( \tH \cup \sP^1(\Q) \bigg)$ are denoted by $e_{\Ga_0(N)}(x)$.
 By {\it loc.\ cit.}, $e_{\Ga_0(N)}(\bstwomat{a}{b}_{\Ga_0(N)})=\frac{N}{d t}$ with $t=\gcd(d, \frac{N}{d})$. Moreover, the cusp $\bstwomat{a}{b}_{\Ga_0(N)}$ of divisor $d$ is $\Q(\zeta_t)$-rational (cf. \cite[(4.9)]{MR800251}).

For $d \mid N$, let $\D_{N, d} \subset \D_{N} $ be the subgroup of divisors supported on the cusps with divisor $d$. 
The group $\GL_2(\Z/N\Z)$ acts on the set of cusps $\partial(X_0(N))$ and therefore on $\D_{N}$. 
The group $T(N)$ acts on $\D_N$ via the natural action of $\GL_2(\mathbb{Z}/N\mathbb{Z})$ on the set of cusps $\partial(X_0(N))$. By \cite[p. 539]{MR800251}, the action of
$T(N)$ on $\D_{N}$ preserves the subspace $\D_{N,d}$. 
 Let $\D_{N,d}(\psi)$ be the $\Z[\psi]$ submodule of $\D_{N,d} \otimes \Z[\psi]$ in which $T(N)$ acts by the character $\psi$. By \cite[p. 540]{MR800251}, define $ \D_{N}(\psi):=\sum\limits_{d \mid N} \D_{N,d}(\psi)$.

As a group $T(N)=(\Z/N\Z)^{\times} \times (\Z/N\Z)^{\times}$ and hence the Pontryagin 
dual $\widehat{T(N)}=\widehat{(\Z/N\Z)^{\times}} \times \widehat{(\Z/N\Z)^{\times}}$. 
Consequently, for any (but fixed) character $\psi: T(N) \ra \CC^{\times}$, there are two primitive Dirichlet characters 
$\epsilon_1, \epsilon_2$ of conductors $M_1,M_2$ divisors of $N$ such that
$$\psi(\psmat{r}{0}{0}{s}) =\epsilon_1^{-1}(s) \epsilon_2(r)$$
for $r,s \in (\Z/N\Z)^{\times}$.

For $r \in (\Z/N\Z)^{\times}$, $\tau_r \in \Gal(\Q[\zeta_N]/\Q)$ acts on $\D_{N}(\psi)$ as $\epsilon_1(r)$. For $l \nmid N$, the Hecke 
operator $T_l$ acts as $\epsilon_1(l)+ l \epsilon_2(l)^{-1}$. 
Following Stevens \cite[p. 540]{MR800251}, define
 $$
      D_{N,d}(\psi) :=\sum\limits_{\substack{ \bstwomat{a}{db}_{\Ga_0(N)}\\ \text{with divisor } d}}  \epsilon_1(b) \cdot \epsilon_2^{-1}(a)
                                  \bstwomat{a}{db}_{\Ga_0(N)} \in \D_{N,d} \otimes \Z[\psi]. $$
           We now recall the following Proposition \cite[Proposition 4.5]{MR800251} from loc.\ cit.

\begin{proposition}
\label{stevendamn}
\begin{enumerate}
\item  The divisor group $\D_{N,d}(\psi) \neq 0$ if and only if 
       $M_2 \mid d, M_1 \mid \frac{N}{d}$ and $\epsilon_1=\epsilon_2^{-1}=\epsilon$.
\item If the divisor group $\D_{N,d}(\psi)$ is nonzero, then it is generated as a $\Z[\epsilon]$-module by the single element $D_{N,d}(\psi)$.
\end{enumerate}
\end{proposition}
\begin{remark}
The coefficients of $D_{N,d}(\psi) $ are the roots of unity. 
\end{remark}

\subsection{Divisor associated to an Eisenstein Series } 
 Recall that $E_2(\Ga_0(N),\CC)
$ is the vector space of all Eisenstein series of weight two for the congruence subgroup $\Ga_0(N)$.
Let $E$ be an Eisenstein series in $E_2(\Ga_0(N),\CC)$. Then $E$ corresponds to the divisor 
$$\delta_{\Ga_0(N)}(E)=2 \pi i \Bigg(\displaystyle\sum_{x \in \partial(X_0(N))} r_{N,E}(x) \{x\}\Bigg).$$
By \cite[Theorem 1.3(a)]{MR800251}, $r_{N,E}(x) = e_{\Ga_0(N)}(x) \cdot a_0(E\mid_{\theta_x})$ and $\theta_x \in \mathrm{SL}_2(\Z)$ are chosen arbitrarily so that $\theta_x \cdot i\infty$ is in the $\Gamma_0(N)$-orbit of $\mathbb{P}^1(\Q)$ representing $x$.  This defines a boundary map:
\[
\delta_{\Ga_0(N)}: E_2(\Ga_0(N),\CC) \rightarrow \Div^0(\partial(X_0(N));\CC).
\]

From the definition, it is clear that the boundary map is linear.   Now, the Hecke operators act on $E_2(\Ga_0(N),\CC)$ and also on $\partial(X_0(N))$. 
For the Hecke operators $T_l$ acting on the Jacobian of the modular curve $X_0(N)$, let $T_l^{t}$ be the dual of the endomorphism of the Jacobian that $T_l$ induces by the functoriality of the Picard group.

We now list some properties of these boundary maps.
\begin{itemize}
\item 
  $\delta_{\Ga_0(N)}(T_l (E)) = T^{t}_l \delta_{\Ga_0(N)}(E)$ (cf. \cite[p. 110]{MR670070}, 
see also~\cite[Lemma 1.11]{MR1376558}).
\item 
If $l \nmid N$, we have $T_l^t=T_l$ \cite[Corollary 2.6, p. 7]{MR4621854}.  Hence, if $l \nmid N$, we have $\delta_{\Ga_0(N)}(T_l (E)) = T_l \delta_{\Ga_0(N)}(E)$.
\end{itemize}

\subsection{Computation of the boundary of the Eisenstein series \texorpdfstring{$E_\varphi$}{Eϕ}}
\label{appendix}
 Let $\varphi$ be a non-trivial Dirichlet character of conductor $f$ 
and $\xi$ be the primitive Dirichlet character associated with $\varphi^2$ of conductor $n$. Let $K_\varphi:=\Q(\zeta_f, \varphi)$ be a suitable finite extension of $\Q$.  
\begin{enumerate}
\item 
 The Bernoulli functions $\mathbb{R} \to \mathbb{R}$ are defined by $B_2(x) = (x - \lfloor x \rfloor)^2 - (x - \lfloor x \rfloor) + \tfrac{1}{6}$,
where $\lfloor x \rfloor$ is the greatest integer less than or equal to $x$. If $\chi$ is a primitive Dirichlet character of conductor $m$, the generalized Bernoulli functions associated with $\chi$ are defined by
\[
B_{2,\chi}(x) := m\left( \sum_{a=0}^{m-1} \chi(a)\, B_2\!\left(\frac{x+a}{m}\right)\right).
\]
The second generalized Bernoulli number associated with $\chi$ is  $B_{2}(\chi):=B_{2,\chi}(0)$.
Further, for a primitive Dirichlet character $\chi$, we have \cite[p.~520]{MR800251}
\[
B_2(\chi) = -2 L(\chi,-1 ).
\]
\item 
 The Gauss sum  $\tau(\chi)$ of $\chi$ is defined by
$\tau(\chi) := \sum\limits_{a=0}^{m-1} \chi(a)\, e\!\left(\frac{a}{m}\right)$, 
where $e(x) = e^{2\pi i x}$.
\item 
For a primitive Dirichlet character $\chi$ of conductor $m$, $R(\chi)$ is defined as 
\[
R(\chi) :=\sum\limits_{a=0}^{m-1} \chi^{-1}(a) \cdot 
\psmat{1}{\frac{a}{m}}{0}{1}.
\]
\item
For any primitive Dirichlet character $\chi$ of conductor $d$, the completed $L$-function
$\Lambda(\chi, s) := (\frac{d}{\pi})^{\frac{s+ \delta }{2}} \Gamma \left( \frac{s+ \delta}{2} \right) L(\chi, s),$
with $\delta = \frac{1 - \chi(-1)}{2}$. 
Then \(\Lambda(\chi,s)\) satisfies the functional equation 
\[
\Lambda(\chi,s)
=
\frac{\tau(\chi)}{i^\delta\sqrt f}
\Lambda(\chi^{-1},1-s).
\]
\end{enumerate}
Equivalently,
\begin{align}
\label{functionalequation}
L(\chi,s)
=
\frac{\tau(\chi)}{i^\delta\sqrt d}
\left(\frac{d}{\pi}\right)^{\frac12-s}
\frac{
\Gamma\!\left(\frac{1-s+\delta}{2}\right)
}{
\Gamma\!\left(\frac{s+\delta}{2}\right)
}
L(\chi^{-1},1-s).
\end{align} 

For any $F \in M_2(\Ga_0(N), \CC)$, the completed $L$-function of $F$ is defined as $D(F,s):= i \Ga(s) (2 \pi  )^{-s} L( F,s)$ and $D(F,s)$ satisfies the functional equation $D(F \mid \bsmat{0}{-1}{1}{0}, s) = -D(F, 2-s)$(~\cite[Lemma 5.1 (a)]{MR800251}).
\begin{lemma}
\label{SpecialvalueMZ}
 For the Eisenstein series 
 \[
 E^\prime := \sum\limits_{a=0}^{f-1} \varphi^{-1}(a) E_{\varphi}\mid  \bsmat{1}{0}{af}{1};
 \]
we have 
\[
L(E', 0)=\frac{f^2}{n^2} \tau(\varphi^{-1}) \tau(\xi) L(\xi^{-1}, -1) \zeta(0) \prod\limits_{p \mid f} (1 - \xi(p)p^{-2})(1-p^{-1})
\]

\end{lemma}
\begin{proof}
For an Eisenstein series $E$ and $d \in \Z$, write down the $L$-series using \cite[p. 544, equation 5.1]{MR800251}:
\[
L(E \mid 
\psmat{d}{0}{0}{1}, s)=d^{1-s} L(E,s).
\]
Define 
\[
E_1 = E_\varphi \mid 
\psmat{0}{-1}{f^2}{0},\qquad  E_2 = E_1 \mid R(\varphi),\qquad  E_3 = E_2| \psmat{0}{-1}{n}{0}.
\]
By~\cite[p. 546, Eq. (5.8)]{MR800251}, we have  $E' = \varphi(-1)E_3 | 
\psmat{\frac{f^2}{n}}{0}{0}{1}$. By $\psmat{0}{-1}{f^2}{0}
=\psmat{0}{-1}{1}{0}
\psmat{f^2}{0}{0}{1}$, we get $L(E_1,s) 
=
f^{2(1-s)}
L(E_\varphi \mid \psmat{0}{-1}{1}{0},s)$.
We have $D(E_{\varphi} \mid \psmat{0}{-1}{1}{0},s)=-D(E_\varphi, 2-s)$ \cite[p. 546]{MR800251}.  This implies 
$$ L(E_1,s) = -\frac{\Gamma(2-s)}{\Gamma(s) } (\frac{2\pi}{f})^{2s-2} L(E_{\varphi},2-s)= -\frac{\Gamma(2-s)}{\Gamma(s) } (\frac{2\pi}{f})^{2s-2}L(\varphi, 2-s) L(\varphi^{-1}, 1-s).$$

Using the functional equation~\ref{functionalequation} for $\varphi$ and $\varphi^{-1}$, we deduce that 
$$L(E_1,s) =  -\frac{\Gamma(2-s)}{\Gamma(s) }(\frac{2\pi}{f})^{2s-2} \Bigg( \frac{\Gamma(\frac{s+ \delta}{2}) \Gamma(\frac{s+ \delta -1}{2} ) }{\Gamma(\frac{2-s+ \delta}{2} )  \Gamma(\frac{1-s+ \delta}{2} ) } \Big( \frac{f}{\pi} \Big)^{2s-2} \frac{\tau(\varphi) \tau(\varphi^{-1})}{i^{2\delta}f}  \Bigg) L(\varphi, s) L(\varphi^{-1}, s-1).   $$

Since $\tau(\varphi) \tau(\varphi^{-1}) = i^{2\delta}f $, we obtain
$$L(E_1,s) =  -\frac{\Gamma(2-s)}{\Gamma(s) }2^{2s-2} \Bigg( \frac{\Gamma(\frac{s+ \delta}{2}) \Gamma(\frac{s+ \delta -1}{2} ) }{\Gamma(\frac{2-s+ \delta}{2} )  \Gamma(\frac{1-s+ \delta}{2} ) }  \Bigg) L(\varphi, s) L(\varphi^{-1}, s-1).   $$

By~\cite[Lemma 5.1 (b)]{MR800251}, we have 
\begin{align*}
L(E_2,s) & = \tau(\varphi^{-1}) L(E_1, \varphi, s)\\
& = - \tau(\varphi^{-1})\frac{\Gamma(2-s)2^{2s-2}}{\Gamma(s)} \Bigg( \frac{\Gamma(\frac{s+ \delta}{2}) \Gamma(\frac{s+ \delta -1}{2} ) }{\Gamma(\frac{2-s+ \delta}{2} )  \Gamma(\frac{1-s+ \delta}{2} ) }  \Bigg)  L(\varphi^2, s)L(\varphi\varphi^{-1}, s-1)\\
&=-C_1 \cdot G_1(s) \cdot H_1(s) L(\varphi^2, s)L(\varphi\varphi^{-1}, s-1)
\end{align*}
for a suitably defined function:
\[
G_1(s):=\frac{\Gamma(2-s)}{\Gamma(s)} \cdot  \Bigg( \frac{\Gamma(\frac{s+ \delta}{2}) \Gamma(\frac{s+ \delta -1}{2} ) }{\Gamma(\frac{2-s+ \delta}{2} )  \Gamma(\frac{1-s+ \delta}{2} ) } \Bigg),  \qquad H_1(s):=2^{2s-2} 
\]
and the constant $C_1= \tau(\varphi^{-1}) $.

Recall that for an imprimitive Dirichlet character $\chi$ of modulus $N$, we have 
\[
L(\chi, s) = L(\chi^*,s) \prod_{p \mid N} (1- \chi^*(p)p^{-s}),
\] 

where $\chi^*$ is the primitive Dirichlet character associated with $\chi$ of conductor $n$.  Let $\xi$ is the primitive character associated with $\varphi^2$ (of conductor $n$). Then we have
\[
L(\varphi^2, s) = L(\xi,s) \prod_{p \mid f} (1- \xi(p)p^{-s}),
\] 
\[
L(\varphi \varphi^{-1}, s) = \zeta(s) \prod_{p \mid f} (1- p^{-s}).
\] 

Thus,
\begin{align*}
  D(E_2,s) &:= i \Ga(s) (2 \pi)^{-s} L(E_2,s)  \\
   &= -(iC_1)(\Gamma(s)G_1(s))(H_1(s)(2 \pi)^{-s})L(\varphi^2, s) L(\varphi \varphi^{-1}, s-1)\\
   &=  -(iC_1)(\Gamma(s)G_1(s))(H_1(s)(2 \pi)^{-s}) [L(\xi,s) \prod_{p \mid f} (1- \xi(p)p^{-s})][\zeta(s-1) \prod_{p \mid f} (1- p^{1-s})]\\
   &= -C_2G_2(s)H_2(s) L(\xi, s)\zeta(s-1),
\end{align*}

with $C_2=i C_1$, $G_2(s)= \Ga(s) \cdot G_1(s)$ and 
\[
H_2(s):= \frac{H_1(s)}{(2 \pi)^s} \prod_{p \mid f} (1- \xi(p)p^{-s})(1- p^{1-s})= \frac{2^{s-2}}{ \pi^s} \prod_{p \mid f} (1- \xi(p)p^{-s})(1- p^{1-s}).
\]
As a consequence, 
$$H_2(2-s) = 2^{-s} \pi^{s-2} \prod_{ p \mid f} [( 1- \xi(p)p^{s-2})(1-p^{s-1}) ]. $$

Now,
\begin{align*}
D(E_3,s)  = n^{1-s} D(E_2 \mid \psmat{0}{-1}{1}{0}, s) = -n^{1-s} D(E_2, 2-s)\\
=C_2 n^{1-s} \cdot G_2(2-s) \cdot H_2(2-s)L(\xi, 2-s) \zeta(1-s)\\
=C_2 \cdot G_2(2-s) \cdot H_3(s)L(\xi, 2-s) \zeta(1-s)\\
\end{align*}
with $H_3(s)=n^{1-s} H_2(2-s)$.\\

Applying this with \(s\mapsto 2-s\), and noting that $\xi(-1)=1$,  we obtain

\[
L(\xi,2-s)
=
\frac{\tau(\xi)}{\sqrt n}
\Big(\frac{n}{\pi}\Big)^{s - \frac{3}{2}}
\frac{
\Gamma\!\left(\frac{s-1}{2}\right)
}{
\Gamma\!\left(\frac{2-s}{2}\right)
}
L(\xi^{-1},s-1)=C_3 G_3(s) H_4(s) L(\xi^{-1},s-1);
\]

with $C_3=\frac{\tau(\xi)}{\sqrt{n}}$, $G_3(s)= \frac{
\Gamma\!\left(\frac{s-1}{2}\right)
}{
\Gamma\!\left(\frac{2-s}{2}\right)
}$, and $H_4(s)=\Big(\frac{n}{\pi}\Big)^{s- \frac{3}{2}}$, for
 the Riemann zeta function, we use
\[
\zeta(1-s)
=
2^{1-s}\pi^{-s}
\cos\!\left(\frac{\pi s}{2}\right)
\Gamma(s)\zeta(s)=  H_5(s) \zeta(s).
\]
Consequently,
\[
D(E_3,s)
= C \cdot G(s) \cdot H(s) L(\xi^{-1},s-1)\zeta(s),
\]

where
\[
C:=C_2 \cdot C_3 = i \tau(\varphi^{-1}) \tau(\xi) n^{-\frac{1}{2}}, G(s):= G_2(2-s) \cdot G_3(s), H(s):=  H_3(s) \cdot H_4(s) \cdot H_5(s) .
\]

Note that $G_2(s) = \Gamma(2-s) \frac{\Gamma(\frac{s+ \delta}{2}) \Gamma(\frac{s + \delta -1}{2})}{ \Gamma(\frac{2-s+ \delta}{2}) \Gamma(\frac{1-s+ \delta}{2}) }$. Consequently, 
$$G_2(2-s) = \Gamma(s) \frac{ \Gamma(\frac{2-s+ \delta}{2}) \Gamma( \frac{1-s + \delta}{2}) }{\Gamma(\frac{s + \delta}{2}) \Gamma( \frac{s-1 + \delta}{2}) } $$
and hence
$$ G(s) = \Gamma(s) \frac{ \Gamma(\frac{2-s+ \delta}{2}) \Gamma( \frac{1-s + \delta}{2}) }{\Gamma(\frac{s + \delta}{2}) \Gamma( \frac{s-1 + \delta}{2}) } \frac{\Gamma(\frac{s-1}{2})}{\Gamma(\frac{2-s}{2})}.$$
Thus, if $\varphi(-1)=1$, that is $\delta =0$, we obtain
$$G(s) = \frac{\Gamma(s) \Gamma(\frac{1-s}{2})}{\Gamma(\frac{s}{2})} = \varphi(-1)\frac{\Gamma(s) \Gamma(\frac{1-s}{2})}{\Gamma(\frac{s}{2})}.$$
On the other hand, if $\varphi(-1)=-1$, that is $\delta =1$, we obtain
$$G(s) = \Gamma(s) \frac{ \Gamma(\frac{1-s}{2} +1)}{\Gamma(\frac{s}{2}) \Gamma( \frac{s-1}{2} +1)} \Gamma(\frac{s-1}{2}).$$
Using the relation $\Gamma(s) = (s-1)\Gamma(s-1)$, we obtain
$$G(s) = \Gamma(s) \frac{\frac{(1-s)}{2} \Gamma(\frac{1-s}{2} ) }{\Gamma(\frac{s}{2}) \frac{(s-1)}{2} \Gamma(\frac{s-1}{2} ) } \Gamma(\frac{s-1}{2}) = - \frac{\Gamma(s) \Gamma(\frac{1-s}{2})}{\Gamma(\frac{s}{2})} = \varphi(-1)\frac{\Gamma(s) \Gamma(\frac{1-s}{2})}{\Gamma(\frac{s}{2})}.$$
Thus, in both cases we obtain,
$$G(s) = \varphi(-1) \frac{\Gamma(s)\Gamma(\frac{1-s}{2})}{\Gamma(\frac{s}{2})} = \varphi(-1) \frac{2^{s-1}\sqrt{\pi}}{\cos(\pi s/2)}.$$
Next, we simplify $H(s)$. First note that,
$$H_2(s) = 2^{s-2} \pi^{-s} \prod_{ p \mid f} [(1 - \xi(p)p^{-s})( 1- p^{1-s}) ],$$
consequently,
$$H_2(2-s) = 2^{-s} \pi^{s-2} \prod_{ p \mid f} [( 1- \xi(p)p^{s-2})(1-p^{s-1}) ]. $$
Hence 
\begin{align*}
    H(s) & =H_3(s)\cdot H_4(s) \cdot H_5(s)\\
    & = n^{1-s} H_2(2-s) (\frac{n}{ \pi})^{s - \frac{3}{2}} 2^{1-s} \pi^{-s} \cos(\frac{\pi s}{2}) \Gamma(s)\\
    & = n^{- \frac{1}{2}} \pi^{- \frac{1}{2} - s} 2^{1-2s} cos(\frac{\pi s}{2}) \Gamma(s) \prod_{ p \mid f} [( 1- \xi(p)p^{s-2})(1-p^{s-1}) ].
\end{align*}

Putting it all together, we obtain
$$D(E_3,s) = i \varphi(-1) \frac{ \tau(\varphi^{-1}) \tau(\xi)}{n} \frac{\Gamma(s) 2^{-s} }{ \pi^s} \prod_{ p \mid f} [( 1- \xi(p)p^{s-2})(1-p^{s-1}) ] L(\xi^{-1},s-1)\zeta(s). $$

Thus,
$$L(E_3,s) := \frac{D(E_3,s) (2\pi)^s}{i \Gamma(s)} = \varphi(-1) \frac{ \tau(\varphi^{-1}) \tau(\xi)}{n}  \prod_{ p \mid f} [( 1- \xi(p)p^{s-2})(1-p^{s-1}) ] L(\xi^{-1},s-1)\zeta(s). $$

Finally, we have 
\begin{align*}
L(E',s) &= \varphi(-1) (\frac{f^2}{n})^{1-s} L(E_3,s) \\
         &= \tau(\varphi^{-1}) \tau(\xi) \frac{ f^{2-2s}}{n^{2-s}} \prod_{ p \mid f} [( 1- \xi(p)p^{s-2})(1-p^{s-1}) ] L(\xi^{-1},s-1)\zeta(s).
\end{align*}

\end{proof}
\begin{proposition}
\label{PR}
 Define
\[
\beta_\varphi := -\frac{f}{4n}\frac{ \tau(\varphi^{-1})}{ \tau(\xi^{-1})} B_2(\xi^{-1}) \prod_{p \mid f} (1 - \xi(p)p^{-2}) \in \Q(\zeta_f, \varphi).
\]

 Then, the divisor associated with the Eisenstein series $E_\varphi$ is given by 
$$ \delta_{\Ga_0(f^2)}(E_{\varphi})= \beta_\varphi D_{f^2,f}(\varphi), $$
where $D_{f^2,f}(\varphi)=\sum  \varphi(ab) \bstwomat{a}{fb}_{\Ga_0(f^2)}$ (cf. Proposition~\ref{stevendamn}). The cusps $\bstwomat{a}{fb}_{\Ga_0(f^2)}$ are $\Q(\zeta_f)$-rational but not $\Q$-rational, and
the divisor $D_{f^2,f}(\varphi) \in  \D_{N}\otimes K_\varphi$.     
\end{proposition}
\begin{proof}
 Recall that $ \D_{f^2}(\psi)=\displaystyle\sum_{d \mid f^2} \D_{f^2,d}(\psi)$ and $ E_\varphi$ form a Hecke eigensystem~\cite[p. 540]{MR800251}:
\[
T_l E_\varphi = (\varphi(l) + l \varphi^{-1}(l)) E_\varphi
\quad (l \nmid f).
\]
Since $T_l ( \delta_{\Ga_0(f^2)}(E_{\varphi}))=(\varphi(l)+l \varphi^{-1}(l))\delta_{\Ga_0(f^2)}(E_{\varphi})$ for 
 $l \nmid f$, the divisor $\delta_{\Ga_0(f^2)}(E_{\varphi}) \in \D_{f^2}(\psi)$ with $\psi(\psmat{r}{0}{0}{s}) =\varphi^{-1}(rs)$.
This implies
$ \delta_{\Ga_0(f^2)}(E_{\varphi}) \in \D_{f^2} (\psi)$. Since $E_\varphi$ is nonzero, the constant term at some cusp is nonzero; hence, the divisor is nonzero . By
Proposition~\ref{stevendamn}, part $(1)$, we deduce that $f\mid d, f \mid \frac{f^2}{d}$ and hence $f=d$. By part $(2)$ of the same proposition, we get $\delta_{\Ga_0(f^2)}(E_{\varphi}) = c(f) \cdot D_{f^2,f} (\varphi)$.  On the other hand,
$$ \delta_{\Ga_0(f^2)}(E_{\varphi})= \displaystyle\sum_{x \in \partial(X_0(f^2))} r_{f^2,E_{\varphi}}(x) \{x\} . $$

 Now, arguing as in the proof of~\cite[Proposition 4.7(d), p. 545]{MR800251}, we conclude that 
 \[
 \beta_\varphi:= c(f) = r_{f^2, E_{\varphi}}(\bstwomat{1}{f}_{\Ga_0(f^2)})  =  e_{\Ga_{0}(f^2)}(\bstwomat{1}{f}_{\Ga_0(f^2)}) a_{0}(E_{\varphi} \mid 
 \bsmat{1}{0}{f}{1}).
 \]
 Note that $e_{\Ga_{0}(f^2)}(\bstwomat{1}{f}_{\Ga_0(f^2)}) = \frac{f^2}{dt}$, where $d=\gcd(f,f^2)=f$, $t = \gcd(f,f)=f$, and hence $e_{\Ga_{0}(f^2)}(\bstwomat{1}{f}_{\Ga_0(f^2)}) = 1$. Therefore,
\begin{equation}\label{coeff}
\delta_{\Ga_0(f^2)}(E_{\varphi}) =  a_{0}(E_{\varphi} \mid \bsmat{1}{0}{f}{1}) \cdot D_{f^2,f} (\varphi).
\end{equation}
For $a \in \Z$, by comparing the coefficients of~\eqref{coeff} at the cusp
$\bstwomat{1}{af}_{\Ga_0(f^2)}$, we obtain 
\[
 a_{0}(E_{\varphi} \mid \bsmat{1}{0}{af}{1}) = \varphi(a) a_{0}(E_{\varphi} \mid \bsmat{1}{0}{f}{1}).
 \]
 Hence, following Stevens \cite[(5.5)]{MR800251}, we define $E^\prime := \sum\limits_{a=0}^{f-1} \varphi^{-1}(a) E_{\varphi}\mid  \bsmat{1}{0}{af}{1}$. By computing the constant term of the Eisenstein series as a sum of the constant terms, we have 
 \[
 a_{0}(E_{\varphi} \mid  \bsmat{1}{0}{f}{1}) = \frac{1}{|(\Z/f\Z)^{\times}|} a_{0}(E').
 \]
Finally, following Stevens (cf. \cite[p. 545-547]{MR800251}), we compute
$a_{0}(E') = -L(0, E')$. 

 \begin{align*}
\beta_\varphi = - \frac{L(E',0)}{|(\Z/f\Z)^{\times}|} & = - \frac{\frac{f^2}{n^2} \tau(\varphi^{-1}) \tau(\xi) L(\xi^{-1}, -1) \zeta(0) \prod\limits_{p \mid f} (1 - \xi(p)p^{-2})(1-p^{-1}) }{ f \prod_{p \mid f} (1-p^{-1})} \\
&=-\frac{f}{4n^2} \tau(\varphi^{-1}) \tau(\xi) B_2(\xi^{-1}) \prod\limits_{p \mid f} (1 - \xi(p)p^{-2}).
\end{align*}
In the last equality, we used $\zeta(0) = - \frac{1}{2}$,  $L( \chi, -1) = -B_{2}(\chi)/2$. The final formula for $\beta_\varphi$ is obtained by using $\tau(\xi)\tau(\xi^{-1}) = \xi(-1)n = n$.
 \end{proof}

\subsection{Computation of boundaries of Eisenstein series \texorpdfstring{$E_{\varphi,M,L}$}{Ephi,M,L}}

In this subsection, we compute the boundary of the Eisenstein series $E_{\varphi,M,L}$
, and it involves the calculation of the ramification indices. 

The following propositions are analogous to~\cite[Lemma 5.2]{MR800251}. 
We now fix a prime number $l$. We denote by $\Ga= \Ga_{0}(A)$, $\Ga^{'}= \Ga_{0}(Al)$. Consider the degeneracy maps $\pi:=\pi_{(l)}: X_{0}(Al) \to X_0(A)$ the map induced by the inclusion of $\Ga_0(Al) \subseteq  \Ga_0(A)$. We also have the map $\pi_{l}: X_{0}(Al) \to X_0(A)$ the map induced from $\gamma \mapsto  \bsmat{l}{0}{0}{1} \gamma \bsmat{l^{-1}}{0}{0}{1}$.  For any map $f:X \rightarrow Y$ between two Riemann surfaces, if $f(x)=y$, then we let $e(f;x,y)$ denote the ramification index of $f$ at the cusp $x$.

\begin{proposition}
\label{pullbackpi}
Let $d$ be a divisor of $A$. Then, 
 \begin{enumerate}
 \item
\begin{align*}
\pi_{(l)}^\ast( D_{A,d}(\varphi)) &= \begin{cases}
e(\pi_{(l)}; \bstwomat{1}{d}_{\Ga_0(Al)}, 
\bstwomat{1}{d}_{\Ga_0(A)}) D_{Al,d}(\varphi) \quad \text{ if 
 } l \mid \frac{A}{d},\\
e(\pi_{(l)}; \bstwomat{1}{d}_{\Ga_0(Al)}, 
\bstwomat{1}{d}_{\Ga_0(A)}) D_{Al,d}(\varphi) + 
\varphi(l) e(\pi_{(l)}; 
\bstwomat{1}{d}_{\Ga_0(Al)}, 
\bstwomat{l}{dl}_{\Ga_0(A)}) D_{Al,dl}(\varphi) \quad  \text{ if  } l \nmid \frac{A}{d}.
\end{cases}\\
\end{align*}
\item
\[
\pi_l^\ast( D_{A,d}(\varphi)) =  \begin{cases}
e(\pi_l; \bstwomat{1}{dl}_{\Ga_0(Al)}, 
\bstwomat{1}{d}_{\Ga_0(A)}) D_{Al,dl}(\varphi) \qquad \qquad \text{  if } l \mid d,\\
e(\pi_l; \bstwomat{1}{dl}_{\Ga_0(Al)}, 
\bstwomat{1}{d}_{\Ga_0(A)}) D_{Al,dl}(\varphi) + 
\quad \varphi(l)  e( \pi_l, \bstwomat{1}{d}_{\Ga_0(Al)}, 
\bstwomat{l}{d}_{\Ga_0(A)})  D_{Al,d}(\varphi) \qquad \text{ if } l \nmid d.
\end{cases}
\]

\end{enumerate}
\end{proposition}
\begin{remark}
The cusps appearing in $D_{Al,dl}(\varphi)$ and $D_{Al,d}(\varphi)$ are disjoint by definition. 
\end{remark}
\begin{proof}
By the definition of the pullback maps on the Picard groups \cite[p. 228]{MR2112196}, 
\[
\pi^{*} D_{A,d}(\psi)=\sum \varphi(b) \pi^{*}\{\frac{1}{db}\}=\sum \varphi(b) (\sum_{y \in  \pi^{-1}\begin{bmatrix}
           1 \\
          db\\
       \end{bmatrix}_{\Ga_0(A)} }e(\pi;x,y) y).
\]

Observe that
\begin{align}
\label{primedivisor}
\pi^{*}(D_{A,d}(\psi)) = \sum_{k} c(k) D_{Al,k}(\psi),
\end{align}
with the sum being over $k \mid Al$.  Equation~\ref{primedivisor} is the equality of two divisors
in the divisor group $\Div^0\bigg(\partial(X_{0}(Al));\Z\bigg)\otimes_\Z \CC$. Observe that $c(k)=0$ if $k \not \in \{d,dl\}$  by looking at Equation~\ref{primedivisor} (see also ~\cite[p. 548]{MR800251}). Let $r(\begin{bmatrix}
           1 \\
          k\\
       \end{bmatrix}_{\Ga_0(A)})$ be the coefficient of $\begin{bmatrix}
           1 \\
          k\\
       \end{bmatrix}_{\Ga_0(A)}$ in $D_{A,d}(\psi)$. We have an equality from Equation~\ref{primedivisor}:
       \begin{align}
       \label{jhattjalano}
   r(\begin{bmatrix}
           1 \\
          k\\
       \end{bmatrix})e (\pi; \begin{bmatrix}
           1 \\
          k\\
       \end{bmatrix}_{\Ga_0(Al)}, \begin{bmatrix}
           1 \\
          k\\
       \end{bmatrix}_{\Ga_0(A)})= c(k). 
       \end{align}

Case 1: 
If $l \mid \frac{A}{d}$, we compute that 
$\gcd(dl,A)=\gcd(dl, \frac{A}{d} d)=d \cdot \gcd(l, \frac{A}{d} )=d \cdot l$.
 Hence, we deduce that $\pi(\begin{bmatrix}
           1 \\
          dl\\
       \end{bmatrix}_{\Ga_0(Al)}) \neq \begin{bmatrix}
           1 \\
        d\\
       \end{bmatrix}_{\Ga_0(A)}$ and $r(\begin{bmatrix}
           1 \\
          dl\\
       \end{bmatrix}_{\Ga_0(A)})=0$. In this case, $\pi^{*}D_{A,d}(\psi)=c(d)  D_{Al,d}(\psi)$. 
Comparing the coefficients for $D_{A,d}(\psi)$ and $D_{Al,d}(\psi)$, we get 
\[
\varphi(d)e (\pi; \begin{bmatrix}
           1 \\
          d\\
       \end{bmatrix}_{\Ga_0(Al)}, \begin{bmatrix}
           1 \\
          d\\
       \end{bmatrix}_{\Ga_0(A)})=c(d) \varphi(d)
\]
and therefore $c(d)=e (\pi; \begin{bmatrix}
           1 \\
          d\\
       \end{bmatrix}_{\Ga_0(Al)}, \begin{bmatrix}
           1 \\
          d\\
       \end{bmatrix}_{\Ga_0(A)})$.

Case 2:-
For $l \nmid \frac{A}{d}$,  we compute
$\gcd(dl,A)=\gcd(dl, \frac{A}{d} d)=d \cdot \gcd(l, \frac{A}{d} )=d$.
Hence,  $\pi(\begin{bmatrix}
           1 \\
          dl\\
       \end{bmatrix}_{\Ga_0(Al)})=\pi(\begin{bmatrix}
           1 \\
          d\\
       \end{bmatrix}_{\Ga_0(Al)})= \begin{bmatrix}
           1 \\
        d\\
       \end{bmatrix}_{\Ga_0(A)}$. In this case, we deduce that
       \[
       \pi^{*} D_{A,d}(\psi)=c(d)  D_{Al,d}(\psi)+c(dl)  D_{Al,dl}(\psi).
       \]
       We compute $c(d)$ and $c(dl)$ using Equation~\ref{jhattjalano}. Comparing the coefficients of cusps of level $d$ in  $D_{A,d}(\psi)$  and cusps of level $ld$ in $D_{Al,dl}(\psi)$, we get 
\[
\varphi(l)\varphi(b)e (\pi; \begin{bmatrix}
           1 \\
          dl\\
       \end{bmatrix}_{\Ga_0(Al)}, \begin{bmatrix}
           1 \\
          dl\\
       \end{bmatrix}_{\Ga_0(A)})=c(dl)\varphi(b)
\]       
and therefore $c(dl)=\varphi(l)e (\pi; \begin{bmatrix}
           1 \\
          dl\\
       \end{bmatrix}_{\Ga_0(Al)}, \begin{bmatrix}
           1 \\
          dl\\
       \end{bmatrix}_{\Ga_0(A)})$. This completes the proof for $\pi = \pi_{(l)}$.

The proof for $\pi_l$ is similar. Observe again that
\begin{align}
\label{primedivisorpil}
\pi_l^{*}(D_{A,d}(\psi)) = \sum_{k} c_l(k) D_{Al,k}(\psi),
\end{align}
with the sum being over $k \mid Al$. Now again $c_l(k)=0$ if $k \not \in \{d,dl\}$ by loc.\ cit.

If $l \mid d$, we again compare the coefficients and get $c_l(dl)=  e(\pi_l; \bstwomat{1}{dl}_{\Ga_0(Al)}, \bstwomat{1}{d}_{\Ga_0(A)})$. Hence, we get 
\[
\pi_l^{*}(D_{A,d}(\psi))=
e(\pi_l; \bstwomat{1}{dl}_{\Ga_0(Al)}, 
\bstwomat{1}{d}_{\Ga_0(A)}) D_{Al,dl}(\varphi).
\]
On the other hand, if $l \nmid d$, the usual comparison of coefficients gives us 
\[
\varphi(lab) e(\pi_l; \bstwomat{1}{dl}_{\Ga_0(Al)}, 
\bstwomat{1}{d}_{\Ga_0(A)})=c_l(d) \varphi(ab)
\]
and $c_l(d)=\varphi(l)$ if $l \nmid d$.  We can also compute $c_l(dl)$ in the same way. Hence, we get the desired result. 
\end{proof}

Next, we compute the required ramification indices for the maps $\pi_{(l)}$ and $\pi_l$. Recall that for an integer $x$, by $\nu_l(x)$ we denote the highest power of $l$ that divides $x$.

\begin{proposition}
\label{ramindexpi}
Let $d$ be a divisor of $A$. 
\begin{enumerate}
\item The ramification indices of $\pi_{(l)}$ are given as follows:
\begin{enumerate}
\item For any prime $l$, we have
$$
e(\pi_{(l)}; \bstwomat{1}{d}_{\Ga^\prime},
\bstwomat{1}{d}_{\Ga})= \begin{cases}
l \qquad \text{ if } \nu_l(d) \le \nu_l(\frac{A}{d}),\\
1 \qquad \text{ otherwise.}
\end{cases}
$$

\item  For any prime $l \nmid \frac{A}{d}$, we have
$ e(\pi_{(l)}; \bstwomat{1}{dl}_{\Ga^\prime}, \bstwomat{1}{dl}_{\Ga})=1
$
\end{enumerate}
\item The ramification indices of $\pi_{l}$ are given as follows:
\begin{enumerate}
\item For any prime $l$, we have
$$
e(\pi_l; \bstwomat{1}{dl}_{\Gamma^\prime},
         \bstwomat{1}{d}_{\Gamma})= 
         \begin{cases}
l \qquad \text{ if }  \nu_l(\frac{A}{d}) \le \nu_l(d),\\
1 \qquad \text{ otherwise.}
\end{cases}
$$

\item  For any prime $l \nmid d$, we have
$ e(\pi_l; \bstwomat{1}{d}_{\Gamma^\prime},
         \bstwomat{l}{d}_{\Gamma})=1$.
\end{enumerate}
\end{enumerate}
\end{proposition}
\begin{proof}
The computation in part (1) follows easily from 
\[
 e(\pi_{(l)}; \bstwomat{a}{db}_{\Gamma^\prime},
\bstwomat{a}{db}_{\Gamma}) = 
\frac{e_{\Gamma_0(Al)}(\bstwomat{a}{db}_{\Gamma^\prime})}{
e_{\Gamma_0(A)}( \bstwomat{a}{db}_{\Gamma}
 )}.
 \]
We now give a proof of part (2). Since the map $\pi_l: X_0(Al) \rightarrow X_0(A)$ is given by $z \rightarrow lz$, 
We now compute the ramification indices for the map $\pi_l$. Recall the description of local charts at the cusps of modular curves~\cite[p. 62]{MR2112196}. 

Consider the following subsets:
\begin{itemize}
\item 
$V, V' \subseteq \tH \cup \sP^1(\Q)$ around the cusps $x$, $y$, resp.,
\item
 let
$U,\widetilde{U}$ be the corresponding neighborhoods on the modular curves $X_0(Al)$, $X_0(A)$, resp.,
\item
 $U', \widetilde{U'}$ be the neighborhoods of  the point $0$ of the unit disc. 
\end{itemize}
We have a commutative diagram: 
$$
\xymatrix{
V \ar[d]^{\pi} \ar[r]^{\pi_l} & V' \ar[d]^{\widetilde{\pi}}\\
U\ar[d]^{\phi} \ar[r]^{\widetilde{\pi_l}} &\widetilde{U} \ar[d]^{\widetilde{\phi}}\\
U'\ar[r]^{\pi^0_l} & \widetilde{U'}\\
}.
$$
We will show that the map $\pi^0_l:U' \rightarrow  \widetilde{U'}$ is $z \rightarrow z^l$.

For a cusp $s$, let $\delta_s \in \SL_2(\Z)$ be such that $\delta_s(s)=\infty$ and $h_s$ be the width of the cusp $s$. Consider the wrapping map $\rho_s(z)=
e^{\frac{2 \pi iz}{h_s}}$ and hence $\rho^{-1}_s(z)=\frac{h_s}{2 \pi i} \log(z)$. In particular, let $h_{\frac{1}{d}}$ (resp., $h_{\frac{1}{dl}}$) be the width of the cusp $\bstwomat{1}{d}_{\Ga_0(A)},$ (resp., $ \bstwomat{1}{dl}_{\Ga_0(Al)}$) for $X_0(A)$  (resp., $X_0(Al)$).

Observe 
that $\pi_l$ is achieved by the matrix $\tau_l=\left(\begin{smallmatrix}
\sqrt{l} & 0\\
0 & \frac{1}{\sqrt{l}} \\
\end{smallmatrix}\right)\in \SL_2(\R)$. Note that $\delta_{ld}=\left(\begin{smallmatrix}
1 & 0\\
-ld & 1 \\
\end{smallmatrix}\right) \in \SL_2(\Z)$  (respectively  $\delta_d=\left(\begin{smallmatrix}
1 & 0\\
-d & 1 \\
\end{smallmatrix}\right) \in \SL_2(\Z)$)  is such that $\delta_{ld}(\frac{1}{ld})=\infty$ and 
$\delta_{d}(\frac{1}{d})=\infty$. A small computation using matrices shows that 
\[
\delta_{d} \circ  \tau_l \circ  \delta_{ld}^{-1}=\tau_l.
\]

By \cite[Section 1.4.1]{MR2289048}, it follows that the stabilizers of the cusps are conjugate and hence
$h_{1/d} = h_{1/(dl)}$.

The induced map on local parameters is
\begin{align*}
\pi_l^0(z)
&= \rho_{1/d} \circ \delta_d \circ \tau_l \circ (\rho_{1/(dl)} \circ \delta_{dl})^{-1}(z) \\
&= \rho_{1/d} \circ \tau_l \circ (\rho_{1/(dl)})^{-1}(z) \\
&= \exp\!\left( \frac{h_{1/(dl)}}{h_{1/d}} \, l \log z \right) \\
&= z^l.
\end{align*}
Thus, the ramification index is $l$. Using the standard formula for cusp widths,
\[
h_{1/d} = \frac{A}{d\, \gcd(d, A/d)}, \qquad
h_{1/(dl)} = \frac{A}{d\, \gcd(dl, A/d)} \qquad t= \gcd(d, \frac{A}{d}), t''= \gcd(dl, \frac{A}{d}).
\]

and writing $d = l^\alpha k_1$, $\frac{A}{d} = l^\beta k_2$, we obtain
\[
\gcd(dl, A/d) =
\begin{cases}
\gcd(d, A/d) & \text{if } \beta \le \alpha,\\
l\cdot \gcd(d, A/d) & \text{otherwise}.
\end{cases}
\]
The ramification index is now given by \[ e(\pi_l; \bstwomat{1}{dl}_{\Gamma^\prime}, \bstwomat{1}{d}_{\Gamma})=\frac{lt}{t''}=l. \] This concludes the proof of part (a).

We now show part (b). Now assume $l \nmid d$. Consider the cusps $\frac{1}{d}$ on $X_0(Al)$ and $\frac{l}{d}$ on $X_0(A)$.
Their widths are
\[
h_{1/d} = \frac{Al}{d\, \gcd(d, A/d)}, 
\qquad
h_{l/d} = \frac{A}{d\, \gcd(d, A/d)}.
\]
Hence
\[
\frac{h_{1/d}}{h_{l/d}} = l.
\]
Since $l \nmid d$, there exist integers $m,n$ such that $ml+dn=1$. Consider the matrix $\delta_{\frac{l}{d}}=\left(\begin{smallmatrix}
m& n\\
-d & l\\
\end{smallmatrix}\right) \in \SL_2(\Z)$. A direct computation shows that
$ \delta_{\frac{l}{d}}  \cdot  \tau_l \cdot \delta_{\frac{1}{d}}^{-1}=\left(\begin{smallmatrix}
\frac{1}{\sqrt{l}}& 0\\
0 & \sqrt{l}\\
\end{smallmatrix}\right) \cdot \left(\begin{smallmatrix}
1& n\\
0 & 1\\
\end{smallmatrix}\right)$. Therefore,
\begin{align*}
\pi_l^0(z)
&= \rho_{l/d} \circ \delta_{l/d} \circ \tau_l \circ (\rho_{1/d} \circ \delta_{1/d})^{-1}(z) \\
&= \exp\!\left( \frac{h_{1/d}}{l\, h_{l/d}} (\log z + n) \right) \\
&= \exp(\log z + n) \\
&= z.
\end{align*}
Thus, the ramification index is $1$.

\end{proof}
Combining Propositions \ref{pullbackpi}, \ref{ramindexpi}, we obtain:
\begin{proposition}\label{induction}
Let $A \in \N$, let $d \mid A$, and let $\varphi$ be a Dirichlet character. 
For each prime $l$, the pullbacks of the divisors $D_{A,d}(\varphi)$ under the degeneracy maps 
$\pi_{(l)}, \pi_l : X_0(Al) \to X_0(A)$ are given as follows. 
(Note that the fiber of a cusp under these maps consists of at most two cusps.)

\begin{enumerate}
\item
\[
\pi_{(l)}^{*}(D_{A,d}(\varphi)) =
\begin{cases} 
lD_{Al,d}(\varphi) + \varphi(l)\, D_{Al,dl}(\varphi) 
& \text{if } l \nmid \frac{A}{d},\; l \nmid d,\\[6pt]

D_{Al,d}(\varphi) + \varphi(l)\, D_{Al,dl}(\varphi) 
& \text{if } l \nmid \frac{A}{d},\; l \mid d,\\[6pt]

l\, D_{Al,d}(\varphi)  
& \text{if } l \mid \frac{A}{d},\; \nu_l(d) \le \nu_l\!\left(\frac{A}{d}\right),\\[6pt]

D_{Al,d}(\varphi) 
& \text{if } l \mid \frac{A}{d},\; \nu_l(d) > \nu_l\!\left(\frac{A}{d}\right).
\end{cases} 
\]

\item
\[
\pi_l^{*}(D_{A,d}(\varphi)) =
\begin{cases}
l\, D_{Al,dl}(\varphi) + \varphi(l)\, D_{Al,d}(\varphi)
& \text{if } l \nmid d,\; l \nmid \frac{A}{d},\\[6pt]

D_{Al,dl}(\varphi) + \varphi(l)\, D_{Al,d}(\varphi) 
& \text{if } l \nmid d,\; l \mid \frac{A}{d},\\[6pt]

l\, D_{Al,dl}(\varphi) 
& \text{if } l \mid d,\; \nu_l\!\left(\frac{A}{d}\right) \le \nu_l(d),\\[6pt]

D_{Al,dl}(\varphi)  
& \text{if } l \mid d,\; \nu_l\!\left(\frac{A}{d}\right) > \nu_l(d).
\end{cases} 
\]
\end{enumerate}
\end{proposition}
Recall that $\gamma_l=\left(\begin{smallmatrix}
l & 0\\
0 & 1\\
\end{smallmatrix}\right) \in M_2(\Z)$.
The map $\pi_l: X_0(Al) \rightarrow X_0(A)$ is given by $z \rightarrow lz$; hence, the pull-back map on the Picard group \cite[p. 228]{MR2112196} satisfies the equality $\pi^*_l(\delta_{\Ga_0(A)}(E))=\delta_{\Ga_0(Al)}(E\mid_{[\gamma_l]}(z)).$ Similarly, $
\pi_{(l)}^*(\delta_{\Ga_0(A)}(E))=\delta_{\Ga_0(Al)}(E(z)).$  Define $[l]_\varphi^{\pm} : \D_N \otimes \CC  \rightarrow  \D_N \otimes \CC $ similar to \S~\ref{Ordcrit} by 
$$
[l]_\varphi^+(z)= \pi_{(l)}^* (z) - \frac{\varphi(l)}{l} \pi_l^*z, [l]_\varphi^- (z)= \pi_{(l)}^* (z) - \varphi^{-1}(l) \pi_l^*z.
$$

By Proposition~\ref{induction}, we obtain the following lemma, which is analogous to Stevens~\cite[Lemma 5.3]{MR800251}.

\begin{lemma}
\label{induction2}
\begin{enumerate}
\item Suppose $l \nmid A$ and $d \mid A$. Then,
$[l]_\varphi^+ D_{A,d}(\varphi) = l (1 - \frac{\xi(l)}{l^2}) D_{Al,d}(\varphi).$

\item Suppose $l \nmid A$, $d \mid A$. Let $S_{\varphi} =\{ l \mid  \varphi(l) = \pm 1 \}$. Then, 
 $$[l]_\varphi^- D_{A,d}(\varphi)= \begin{cases}
(l-1)(\beta_{l,0} D_{Al,d}(\varphi) + \beta_{l,1} D_{Al, dl}(\varphi)) \quad \text{ if } l \in S_{\varphi},\\
\beta_{l,0} D_{Al,d}(\varphi) + \beta_{l,1} D_{Al, dl}(\varphi) \qquad \qquad \quad \text{if } l \notin S_{\varphi},
\end{cases}$$ 
where 
$$ \beta_{l,0} = \begin{cases}
                    1 \text{\qquad \quad \   if } l \in S_{\varphi},\\
                    l-1 \text{\qquad if } l \notin S_{\varphi}
                 \end{cases} 
   \text{ and  } \beta_{l,1} 
                     = \begin{cases} 
                         -\varphi(l) \text{\qquad \qquad \quad \  if } l \in S_{\varphi},\\                            
                          \varphi(l) - l \varphi^{-1}(l) 
                          \text{\qquad if } l \notin S_{\varphi}. \end{cases}$$

\item Suppose $l \mid\mid A$, $l \nmid d$ and $d \mid A$. Then,
$[l]^- D_{A,d}(\varphi)=
(l-1) D_{Al,d}(\varphi) - \varphi^{-1}(l) D_{Al, dl}(\varphi).$
\end{enumerate}
\end{lemma}

\begin{proof}
All these formulas are an easy consequence of Proposition \ref{induction}. We only show part (1) for illustration purposes.
The calculations in the remaining parts are similar. Observe that
$$
[l]_\varphi^+ D_{A,d}(\varphi) = \pi_{(l)}^* D_{A, d} (\varphi) - \frac{\varphi(l)}{l} \pi_l^*D_{A,d}(\varphi).
$$
Since $l \nmid A$ implies $l \nmid \frac{A}{d}$ and $l \nmid d$, applying Proposition \ref{induction} yields the result. \qedhere

\end{proof}

\begin{lemma}\label{boundarypm}
Let $E \in E_2(\Gamma_0(M),\mathbb{C})$ be an Eisenstein series. Then
\[
\delta_{\Gamma_0(Ml)}\bigl([l]_\varphi^\pm E\bigr)
=
[l]_\varphi^\pm \bigl(\delta_{\Gamma_0(M)}(E)\bigr).
\]
\end{lemma}

\begin{proof}
We prove the statement for $[l]_\varphi^+$; the case of $[l]_\varphi^-$ is similar.

By definition,
\[
[l]_\varphi^+ E
=
\pi_{(l)}^* E - \frac{\varphi(l)}{l}\, \pi_l^* E.
\]
Applying $\delta_{\Gamma_0(Ml)}$ and using the compatibility of $\delta$ with pullbacks, we obtain
\begin{align*}
\delta_{\Gamma_0(Ml)}\bigl([l]_\varphi^+ E\bigr)
&=
\delta_{\Gamma_0(Ml)}(\pi_{(l)}^* E)
-
\frac{\varphi(l)}{l}\, \delta_{\Gamma_0(Ml)}(\pi_l^* E) \\
&=
\pi_{(l)}^* \delta_{\Gamma_0(M)}(E)
-
\frac{\varphi(l)}{l}\, \pi_l^* \delta_{\Gamma_0(M)}(E) \\
&=
[l]_\varphi^+ \bigl(\delta_{\Gamma_0(M)}(E)\bigr).
\end{align*}
\end{proof}
From now on, let us assume that $N = p^k N_2$ for some $k \ge 2$, where $N_2$ is a square-free integer such that every prime divisor $q \mid N_2$ satisfies $q \equiv \pm 1 \pmod{p^{\lfloor k/2 \rfloor}}$.
For a real number $r$, we write $\lfloor r \rfloor$ for the greatest integer less than or equal to $r$. Let $\varphi$ be a non-trivial Dirichlet character of conductor $p^r$ with $r \le \lfloor k/2 \rfloor$, let $M \mid N_2$, and let $L \mid \frac{N_2}{M}$. The following theorem determines the divisor associated with the non-rational Eisenstein series $E_{\varphi,M,L} \in M_2(\Gamma_0(N))$. The general outline of the proof is similar to that of \cite[Proposition 4.7]{MR800251}.

 For $n \in \mathbb{N}$, let $\phi(n)$ denote the cardinality of $(\mathbb{Z}/n\mathbb{Z})^\times$.
\begin{theorem}\label{boundary}
Let $p$ be a prime and let $N = p^k N_2$ with $k \ge 2$, where $N_2$ is square-free and every prime divisor $q \mid N_2$ satisfies
\[
q \equiv \pm 1 \pmod{p^{\lfloor k/2 \rfloor}}.
\]
Let $\varphi$ be a non-trivial Dirichlet character of conductor $f=p^r>1$ with $r \le \lfloor k/2 \rfloor$. Let $M \mid N_2$ and $L = N_2/M$, and consider the Eisenstein series $E_{\varphi,M,L} \in E_2(\Gamma_0(N),\mathbb{C})$.
Let $\xi$ denote the primitive Dirichlet character associated with $\varphi^2$, of conductor $n$. Then
\[
\delta_{\Gamma_0(N)}(E_{\varphi,M,L})
=
\beta_{\Gamma_0(N),\varphi,M,L}\, D_{N,M,L}(\varphi),
\]
where

\begin{equation}\label{betaimp}
\beta_{\Gamma_0(N),\varphi,M,L}
= -
\frac{p^{k-r} M \phi(L)}{4 n}
\cdot
\frac{\tau(\varphi^{-1})}{\tau(\xi^{-1})}
\cdot
B_2(\xi^{-1})
\cdot
\prod_{l \mid pM}
\left(1 - \frac{\xi(l)}{l^2}\right)
\in \mathbb{Q}(\zeta_f,\varphi),
\end{equation}

and
\[
D_{N,M,L}(\varphi) \in \Div^0\bigl(X_0(N)\bigr)\bigl(\mathbb{Q}(\zeta_f,\varphi)\bigr).
\]
Moreover, the coefficient of each cusp appearing in the support of $D_{N,M,L}(\varphi)$ is a unit in the ring of integers of $\mathbb{Q}(\zeta_f,\varphi)$.
\end{theorem}

\begin{proof}
By Proposition~\ref{PR}, we have for $f=p^r$
\[
\delta_{\Gamma_0(f^2)}(E_\varphi)
=
\beta_\varphi\, D_{f^2,f}(\varphi),
\]
where
\[
\beta_\varphi
=
-\frac{f}{4 n}
\cdot
\frac{\tau(\varphi^{-1})}{\tau(\xi^{-1})}
\cdot
B_2(\xi^{-1})
\cdot
\left(1 - \frac{\xi(p)}{p^2}\right).
\]

We construct $E_{\varphi,M,L}$ from $E_\varphi$ by applying the operators $[l]_\varphi^+$ for primes $l \mid M$ and $[q]_\varphi^-$ for primes $q \mid L$. Since these operators commute and the boundary map $\delta$ is compatible with pullbacks (Lemma~\ref{boundarypm}), we may apply them successively to the divisor.

First, applying $[l]_\varphi^+$ for all $l \mid M$ and using Lemma~\ref{induction2}(1), we obtain
\begin{align*}
\delta_{\Gamma_0(f^2M)}(E_{\varphi,M,1})
&=
\prod_{l \mid M} [l]_\varphi^+ \, \delta_{\Gamma_0(f^2)}(E_\varphi) \\
&=
\beta_\varphi \prod_{l \mid M} [l]_\varphi^+ D_{f^2,f}(\varphi) \\
&=
\beta_\varphi \prod_{l \mid M} l \left(1 - \frac{\xi(l)}{l^2}\right) D_{f^2M,f}(\varphi) \\
&=
\beta_\varphi \, M \prod_{l \mid M} \left(1 - \frac{\xi(l)}{l^2}\right) D_{f^2M,f}(\varphi).
\end{align*}

Next, let $L = q_1 \cdots q_s$ be the prime factorization. By assumption, each $q_i$ satisfies $\varphi(q_i)=\pm 1$. Applying $[q_i]_\varphi^-$ successively and using Lemma~\ref{induction2}(2), we obtain
\[
\delta_{\Gamma_0(f^2ML)}(E_{\varphi,M,L})
=
\beta_\varphi \, M \, \phi(L)
\prod_{l \mid M} \left(1 - \frac{\xi(l)}{l^2}\right)
\widetilde{D},
\]
where
\[
\widetilde{D}
=
\sum_{j_1,\dots,j_s \in \{0,1\}}
\beta_{q_1,j_1} \cdots \beta_{q_s,j_s}
\, D_{f^2ML,\, f q_1^{j_1}\cdots q_s^{j_s}}(\varphi),
\]
with coefficients $\beta_{q_i,j_i} \in \{\pm 1\}$.

Finally, we pass from level $f^2ML$ to level $N = p^k N_2$ by applying $(\pi_{(p)}^*)^{k-2r}$. Since $\nu_p(f^2ML)=2r$ and $\nu_p(f q_1^{j_1}\cdots q_s^{j_s})=r$, Lemma~\ref{induction2}(3) shows that each application multiplies the divisor by $p$, hence
\[
\delta_{\Gamma_0(N)}(E_{\varphi,M,L})
=
\beta_\varphi \, M \, \phi(L) \, p^{k-2r}
\prod_{l \mid M} \left(1 - \frac{\xi(l)}{l^2}\right)
D_{N,M,L}(\varphi).
\]

Substituting the expression for $\beta_\varphi$ yields \eqref{betaimp}.

Finally, the cusps appearing in the divisors $D_{N,d}$ for distinct $d$ are disjoint; hence, the coefficients in $D_{N,M,L}(\varphi)$ are units in the ring of integers of $\mathbb{Q}(\zeta_f,\varphi)$.
\end{proof}
\section{Cuspidal subgroup associated to {\it non-rational} Eisenstein series \texorpdfstring{$E_{\varphi,M,L}$}{Ephi,M,L}}
\label{Cuspidal}
Let $\J_0(N):= \Jac(X_0(N))$ be the Jacobian variety of the modular curve $X_0(N)$  and $C_{\Ga_0(N)}(\overline{\Q})$ be the subgroup of $J_0(N)(\overline{\Q})$
consisting of equivalence classes of degree $0$ divisors of $X_0(N)(\overline{\Q})$ supported on the $\overline{\Q}$-rational cusps of $X_0(N)(\overline{\Q})$. 
Let $E \in E_2(\Ga_0(N),\CC)$ and $\delta_{\Ga_0(N)}(E)$ be the associated degree $0$ divisor of $E$ that is supported at the cusps of $X_0(N)$. Let $R_N(E)$ be the sub-$\Z$
module of $\CC$ generated by the coefficients of $\delta_{\Ga_0(N)}(E)$. Let 
\[
\widecheck{R_N(E)}:=\Hom_{\Z}(R_N(E),\Z)
\]
 be its dual $\Z$ module. 
 \begin{definition}
 The cuspidal subgroup $C_{\Ga_0(N)}(E)$ associated with $E$ is defined to be the image of the composition:
 \[
 \widecheck{R_N(E)} \rightarrow{} \Div^0(\partial(X_0(N))) \rightarrow  C_{\Ga_0(N)}(\overline{\Q}).
 \]
 The first map is given by $\phi \rightarrow \phi(\delta_{\Ga_0(N)}(E))$. Here, $\phi(\delta_{\Ga_0(N)}(E))$ is the divisor obtained by applying $\phi$ to the coefficients of $\delta_{\Ga_0(N)}(E)$.   This is a subgroup of $C_{\Ga_0(N)} (\overline{\Q})\subset J_0(N)(\overline{\Q})$ generated by $\{\phi \circ \delta_{\Ga_0(N)}(E)\}_{\phi \in \widecheck{R_N(E)}}$ ~\cite[p.41]{MR670070}. 
  \end{definition}
 By a classical theorem due to Manin-Drinfeld, $C_{\Ga_0(N)}$ (and hence $C_{\Ga_0(N)}(E)$) is a finite group. Let $w_N$ be the Atkin-Lehner involution of level $N$. 
 Then, the cuspidal subgroup $C^{Ren}_{\Ga_0(N)}(E)$ associated with $E$ is defined to be the subgroup of  $C_{\Ga_0(N)} (\overline{\Q})\subset J_0(N)(\overline{\Q})$  generated by $\{w_N \circ \phi \circ \delta_{\Ga_0(N)}(E)\}_{\phi \in \widecheck{R_N(E)}}$ ( \cite[Definition 4.6]{MR4621854}, ~\cite[p.3]{MR3805452} ). Note that this definition is slightly different from that of Stevens that we use in this paper, but the order of the cuspidal subgroup remains unchanged since $w_N$ is an isomorphism.  This small change ensures the validity of Lemma~\ref{annihilates}, but the order of the cuspidal subgroup remains unchanged since the Atkin-Lehner involutions are isomorphisms. We have 
   \[
   |C^{Ren}_{\Ga_0(N)}(E)|=|C_{\Ga_0(N)}(E)|.
   \]

Let $E \in E_2(\Gamma_0(N),\CC)$ be an Eisenstein series for the congruence subgroup $\Ga_0(N)$ with coefficients in $\CC$.  Let $\omega_E$ be the differential in $X_0(N)(\CC)$ whose pullback to $\tH$ is the differential form $E(z)dz$. The differential $\omega_E$ is holomorphic in $Y_0(N)(\CC)$.
\begin{definition}(The group of periods of Eisenstein series $E$)
Let $[c]$ be the homology class represented by the $1$-cocycle $c$ on $Y_0(N)(\CC)$. We define a homomorphism $\pi_E: \HH_1(Y_0(N)(\CC),\Z) \rightarrow \CC$
by 
\[
[c] \rightarrow \int_{c} \omega_E. 
\]
The homomorphism $\pi_E:\HH_1(Y_0(N),\Z) \rightarrow \CC$ ~\cite[p. 52]{MR670070}
is the ``period'' homomorphism of $E$.
Let $P_N(E)$ be the image of the homomorphism $\pi_E$ inside $\CC$. This is the group of periods of the Eisenstein series $E$.
\end{definition}
The residue is obtained by integrating along a small circle around a cusp.  Hence, we have $R_N(E) \subseteq P_N(E)$. 
For more details on the group $P_N(E)$, the reader can refer to~\cite[p.3]{MR3805452}. 
\begin{definition}
\label{eisensub}
(The dual cuspidal subgroup associated with an Eisenstein series $E$)
For any Eisenstein series $E \in E_2(\Ga_0(N),\CC)$, we associate an abelian group with $E$ as
$$ A_N(E):=P_N(E)/R_N(E).$$
\end{definition}
We list some properties of $ A_N(E)$ and $C_{\Ga_0(N)}(E)$
that we crucially use in the rest of the article~\cite[Proposition 1.1]{MR800251}, 
~\cite[Theorem 1.2]{MR800251}:
\begin{itemize}
\item 
The abelian group $A_N(E)$ is {\it finite}. 
\item
There exists a
non-degenerate pairing $C_{\Ga_0(N)}(E)\times A_N(E) \rightarrow \Q/ \Z$. Hence, these two groups are Pontryagin duals to each other. 
 \end{itemize}
Now, $C_{\Ga_0(N)}(E)$ and $A_N(E)$ are finite abelian groups. 
Assume that 
\[
 A_{N}(E) \simeq \bigoplus\limits_{i=1}^n \Z/p^{n_i}\Z.
\]

We find the structure of the finite abelian group $C_{\Ga_0(N)}(E)$:
\begin{align*}
C_{\Ga_0(N)}(E) &\simeq& \Hom_{\Z} (A_{N}(E), \Q/ \Z) 
&\simeq& \Hom_{\Z } (\bigoplus\limits_{i=1}^n \Z/p^{n_i}\Z, \Q/\Z) \\
&\simeq & \bigoplus\limits_{i=1}^n \Hom_{\Z } ( \Z/p^{n_i}\Z, \Q/\Z)
&\simeq & \bigoplus\limits_{i=1}^n \Z/p^{n_i}\Z.\\
\end{align*}
\subsection{Structure of the dual cuspidal subgroup associated to an Eisenstein series}
We generalize~\cite[Theorem 1.3(b)]{MR800251} to the case of $\Gamma_0(N)$ and use that to study the finite abelian groups $A_N(E_{\varphi, M,L})$ and $C_{\Ga_0(N)}(E_{\varphi,M,L})$.
\begin{definition}\label{setS}
Following~\cite[(1.8)]{MR800251} (see also \cite[p. 30]{MR4621854}), let  $S_N$ be a set of primes such that 
\begin{itemize}
\item 
every element of $S_N$ is $3$ modulo $4$, and 
\item
 for all $k \in \Z$, there exists an $m \in S_N$ such that $m \equiv -1 \pmod{Nk}$. 
\end{itemize}
\end{definition}
Let $\x_{N}$ be the set of non-quadratic Dirichlet characters whose conductors are in $S_N$. We denote by $\x_{N}^{+}$ (resp., $\x_{N}^{-}$) the set of even (resp., odd) characters of $\x_{N}$.  Let $\Z[\x_{N}]$ be the ring generated by the values of the characters in $ \x_N$. 

For a nontrivial primitive Dirichlet character $\chi$ of conductor $t$  with $(t,N)=1$ and $t \equiv -1 \pmod{N^2}$, Stevens ~\cite[Equations (2.1) and (2.2)]{MR800251} considers some special values as 
$$ \Lambda(\chi):=\sum_{\substack{a=0,\\ (a,t)=1}}^{t-1} \chi^{-1}(N a) \left \{0, \frac{Na}{t} \right \}_{\Ga_0(N)}  \in \HH_1(X_0(N),\Z[\chi]), $$
where $\left \{0, \frac{Na}{t} \right \}_{\Ga_0(N)} \in \HH_1(X_0(N),\Z)$ is the homology class represented by an oriented geodesic in $\mathbb{H}$ joining $0$ and $\frac{Na}{t}$. If $t=m^n$ for some $m \in S_N$, then let
$$ \Lambda_{\pm}(\chi):=\sum_{\substack{a=0,\\ \chi_m(Na)= \pm 1}}^{t-1} \chi^{-1}(N a) \left \{0, \frac{Na}{t} \right \}_{\Ga_0(N)}, $$
where $\chi_{m}(\cdot)=(\frac{\cdot}{m})$ is the quadratic residue symbol modulo $m$. Finally, for any $\Z[\chi]$ module $A$ and $\Phi: \HH_1(X_0(N),\Z) \to A$, he defines the twisted special values ~\cite[(2.4)]{MR800251} as
$$ \Lambda_{\pm}(\Phi, \chi) := (\Phi \otimes 1)(\Lambda_{\pm}(\chi)).
$$
The following theorem is a generalization of~\cite[Theorem 2.1]{MR800251} for $\Ga_0(N)$. 
\begin{theorem}
\label{threestep}
Let $S_N$ be a set of primes as in Definition~\ref{setS} 
and $A \subseteq \CC$ be a $\Z[\x_{N}]$-module. 
For a cohomology class $\Phi: \HH_1(X_0(N), \Z) \rightarrow A $ such that  $\Lambda_{\pm}(\Phi, \chi)=0$  for all Dirichlet characters $\chi \in \x_N^+$, we have that $\Phi$ is zero. 
\end{theorem}

\begin{proof}
Choose a prime ideal $\gP \subseteq \Z[\x_{N}]$ and
consider the localization $A_{\gP}$ 
of $A$ at $\gP$.  Consider now the cohomology class $\Phi_{\gP}$ obtained by composing $\Phi$ 
with the inclusion $A \rightarrow A_{\gP}$.
We define $\Phi_{\gP}: \Gamma_{0}(N) \to A_{\gP}$ by defining $\Phi_{\gP}(\gamma) := \Phi_{\gP}(\{x, \gamma x\})$. 
We show $\Phi_{\gP}=0$ in the following 4 steps.

\underline{Step 1}:
If $m \in S_N, m \equiv -1 \pmod{N^{2}}$ and $ \gP \nmid \left(\frac{m-1}{2} \right)\Z[\x_{N}]$, then $\Phi_{\gP}(\{0, \frac{Na}{m}\}_{\Ga_0(N)})=0$ whenever $(m, Na)=1$. We remark that the proof of this step is essentially the same as in Step $1$ of ~\cite[Theorem 2.1]{MR800251}.

\underline{Step 2}:
In this step, we show that 
$\Phi_{\gP}$ vanishes on the subgroup:
$$\Gamma^{'}(N) = \Big\{\begin{pmatrix}
                                    a & bN \\
                                  cN & d
                                 \end{pmatrix} \in SL_{2}(\Z)| a,d \equiv \pm 1 \pmod{N^2} \Big\}.$$

The proof of this step is the same as in Step $2$ of ~\cite[Theorem 2.1]{MR800251}.

\underline{Step 3}:
Let 
$$\Gamma_{0}^{'}(N) = \Big\{\begin{pmatrix}
                                    a & bN \\
                                  cN & d
                                 \end{pmatrix}\Big\}  \subset \Ga_{0}(N),$$
then $\Phi_{\gP}$ vanishes on $\Gamma_{0}^{'}(N)$. Since we have $\gcd(a,N)=\gcd(d,N)=1$, which is equivalent to saying $\gcd(a,N^2) = \gcd(d,N^2)=1$, we see that $$a^{\frac{\phi(N^2)}{2}} = d^{\frac{\phi(N^2)}{2}}= \pm 1 \pmod{N^2},$$ here $\phi$ denotes Euler's totient function. Let $\gamma=\begin{pmatrix}a & bN \\cN & d \end{pmatrix} \equiv \begin{pmatrix} a & 0 \\ 0 & d \end{pmatrix} \pmod{N} \in \Gamma_{0}^{'}(N)$.
Then 
\[
\gamma^{\frac{\phi(N^2)}{2}} \equiv \begin{pmatrix} a^{\frac{\phi(N^2)}{2}} & 0 \\ 0 & d^{\frac{\phi(N^2)}{2}}\end{pmatrix} \pmod{N}.
\]
Since $a^{\frac{\phi(N^2)}{2}} = d^{\frac{\phi(N^2)}{2}}= \pm 1 \pmod{N^2}$, we see that $\gamma^{\frac{\phi(N^2)}{2}} \in \Gamma^{'}(N)$.
Thus $\Phi_{\gP}(\gamma^{\frac{\phi(N^2)}{2}})=0$ by Step 2. Since $\Phi_{\gP}$ is a homomorphism, we have,
$$\Phi_{\gP}(\gamma^{\frac{\phi(N^2)}{2}}) = \frac{\phi(N^2)}{2} \Phi_{\gP}(\gamma)=0.$$

Since $A$ (and therefore $A_{\gP}$) is torsion-free, we have $\Phi_{\gP}(\gamma)=0$.

\underline{Step 4}:

In this step, we show that $\Phi_{\gP}=0$.  Let $\alpha \in \Ga_0(N)$ be 
the matrix $\alpha= \left(\begin{array}{cc}
a & b\\
Nc & d\\
\end{array}\right)$, now consider a parabolic matrix $\beta= \left(\begin{array}{cc}
1 & x\\
0 & 1\\
\end{array}\right)$, here $x$ is chosen such that $ax \equiv b \pmod N$.
 Clearly,  we then have $\alpha \beta^{-1} \in  \Gamma_{0}^{'}(N)$. By Step 3, $\Phi_{\gP}(\alpha \beta^{-1} )=0$ and by a well-known theorem due to Manin~\cite{MR0314846}, we have $\Phi_{\gP}( \beta)=0$. Hence, we conclude that $\Phi_{\gP}(\alpha)=\Phi_{\gP}(\alpha \beta^{-1} )+\Phi_{\gP}( \beta ) =0$.

Since the homomorphism $\Gamma_{0}(N) \to H_{1}(X_{0}(N), \Z)$ is surjective, we see $\Phi_{\gP}=0$. 
Since $\Phi_{\gP}=0$ for all primes $\gP$, we get $\Phi=0$.
\end{proof}

The {\it ``twisted special values"} of the $L$-function of an Eisenstein series $E$ are defined as follows:
$$\Lambda(E, \chi, 1) := \frac{\tau(\chi^{-1})}{2 \pi i} L(E, \chi, 1) \text{ \qquad and \qquad}  \Lambda_{\pm}(E, \chi, 1) :=\frac{1}{2}
[\Lambda(E, \chi, 1) \pm \Lambda(E, \chi \chi_m, 1)]. $$

\begin{theorem}
\label{main}(Stevens)
Let $\widetilde{R} \subset \CC$ be a Dedekind domain. Let $M \subseteq \CC$ be a finitely generated $\widetilde{R}$-module inside $\CC$. Let $S_N$ be a set of primes as in Definition~\ref{setS}. Then,  the following are equivalent
\begin{enumerate}
\item 
 $P_N(E) \otimes_\Z \widetilde{R} \subseteq M$;
\item 
$R_N(E) \otimes_\Z \widetilde{R}\subseteq M$ and $\Lambda_{\pm}(E,\chi,1) \in M[ \chi, \frac{1}{m_{\chi}}]$ for every $\chi \in \x_N$ with conductor $m_{\chi}$.
\end{enumerate}
\end{theorem}
\begin{proof}
Suppose that the set $P_N(E)  \otimes_\Z \widetilde{R}$ is contained in $M$.
Since $R_N(E)  \otimes_\Z \widetilde{R} \subset P_N(E)  \otimes_\Z \widetilde{R}$, we have $R_N(E)  \otimes_\Z \widetilde{R} \subset M$.  Part (2) follows from~\cite[Lemma 2.2(1)]{MR800251}. Note that $\widetilde{R}$ is a torsion-free ring and hence a flat $\Z$ module.   

For the converse part, the proof is the same as in the proof of \cite[Theorem 1.3(b), p. 528]{MR800251} with the
following modifications:
\begin{enumerate}
\item
Use Theorem \ref{threestep} instead of~\cite[Theorem 2.1]{MR800251}.
\item 
Work with the prime ideal $\gP$ of $\widetilde{R}$ instead of the prime ideal $p$ of $\Z$. 
\end{enumerate}
\end{proof}

In preparation for proving Theorem~\ref{ordercuspidal}, we apply the general theory outlined above to the case of Eisenstein series $E_{\varphi, M, L}$ as in Theorem~\ref{boundary}. In the following proposition, we calculate the twisted special values of $E_{\varphi, M, L}$.
\begin{proposition}
\label{specialvalues}
Let $\chi$ be a primitive Dirichlet character with conductor $m_{\chi}$ and $(m_{\chi},N)=1$.  Let $\x_N$ be the set of Dirichlet characters defined just after Definition~\ref{setS}. Then, we have 
\begin{enumerate}
\item $
\Lambda(E_{\varphi, M,L}, \chi, 1) =  \frac{\varphi(m_{\chi}) \chi(f) }{2 \tau(\varphi^{-1})}  \prod\limits_{l \mid M } \Big(1- \frac{\chi\varphi(l)}{l} \Big) \prod\limits_{q \mid L}(1- \chi \varphi^{-1}(q)) B_{1}(\chi^{-1}\varphi^{-1}) B_{1}(\chi \varphi^{-1}).
$
\item 
$\Lambda_{\pm}(E_{\varphi, M,L}, \chi, 1) =
 \frac{\varphi(-m_{\chi})\tau(\varphi) \chi(f)}{f}  \prod\limits_{l \mid M } \Big(1- \frac{\chi\varphi(l)}{l} \Big) \prod\limits_{q \mid L}\Big(1- \chi \varphi^{-1}(q)\Big)  \frac{B_{1}(\chi^{-1}\varphi^{-1})}{2} \frac{B_{1}(\chi \varphi^{-1})}{2},$
if $\chi \in \x_{N}^{-\varphi(-1)}$.
\end{enumerate}
\end{proposition}

\begin{proof}
From Proposition \ref{lvalue} we have the following.
\begin{equation*}
\label{Lchi}
L(E_{\varphi,M,L}, \chi, s) = \Big( \prod_{l \mid M} (1 - \chi(l)\varphi(l)l^{-s}) \Big) \Big( \prod_{q \mid L} (1- \chi(q)\varphi^{-1}(q) q^{1-s}) \Big) L( \chi\varphi, s) L( \chi\varphi^{-1}, s-1).
\end{equation*}
By~\cite[p. 544]{MR800251}, we have
$
L(\varphi\chi, 1) =\frac{ \pi i}{\tau(\varphi^{-1}  \chi^{-1})} L( \varphi^{-1} \chi^{-1}, 0 ), \text{\quad} L( \varphi, 0) = -B_1(\varphi), \text{\quad and \quad} \frac{\tau(\chi^{-1})}{ \tau(\chi^{-1} \varphi^{-1})}=\frac{\chi(f) \cdot  \varphi(m_{\chi})}{ \tau(\varphi^{-1})}.
$
Combining these identities, we obtain
\begin{align*}
 \Lambda(E_{\varphi, M,L}, \chi, 1) 
 &: = \frac{\tau(\chi^{-1})}{2 \pi i} L(E_{\varphi, M,L}, \chi, 1 )\\
 & = \frac{\tau(\chi^{-1})}{2 \pi i} \Big( \prod_{l \mid M} (1 - \chi(l)\varphi(l)l^{-1}) \Big) \Big( \prod_{q \mid L} (1- \chi(q)\varphi^{-1}(q)  \Big) L( \chi\varphi, 1) L( \chi\varphi^{-1}, 0)\\
    &=  \frac{\tau(\chi^{-1})}{2 \tau(\varphi^{-1} \chi^{-1})}  \prod_{l \mid M } \Big(1- \frac{\chi(l)\varphi(l)}{l} \Big)\prod_{q \mid L}(1- \chi(q) \varphi^{-1}(q))
 B_1(\varphi^{-1} \chi^{-1}) B_1(\chi\varphi^{-1})\\
   & =  \frac{\varphi(m_{\chi}) \chi(f)}{2 \tau(\varphi^{-1})}   \prod_{l \mid M } \Big(1- \frac{\chi(l)\varphi(l)}{l} \Big) \prod_{q \mid L}(1- \chi(q) \varphi^{-1}(q)) B_{1}(\varphi^{-1}\chi^{-1}) B_{1}(\chi \varphi^{-1}).
 \end{align*}
 This completes part (1) of the proposition. If $\chi \in \x_{N}^{-\varphi(-1)}$, then $m$ is a prime with $m \equiv 3 \pmod{4}$. In particular, we get $\chi_m(-1)=-1$ and $\varphi \chi(-1)=-1$. This implies that $\varphi^{-1} \chi^{-1} \chi_m^{-1}$ 
 is an even character, and hence $B_1(\varphi^{-1} \chi^{-1} \chi_m^{-1})=0$. Part (1) of the proposition implies that $\Lambda(E_{\varphi, M,L}, \chi\chi_m, 1)=0$. Since $\varphi(-1) \tau(\varphi) \tau(\varphi^{-1}) = f$, part (2) of the proposition follows.
\end{proof}
\subsection{Proof of Theorem~\ref{ordercuspidal}}
Now, we are ready to give a proof of Theorem~\ref{ordercuspidal}, and the argument is similar to ~\cite[p 542]{MR800251}.
Note that $K_\varphi \cong \Q(\zeta_{p^r}, \zeta_{p-1}^a)$ for some integer $a$ and consequently, one sees that $\mathcal{O}_{K_\varphi} \cong \Z[\zeta_f, \varphi]$. Recall our assumption about $N$ from Theorem \ref{ordercuspidal}. With the notation as above and by Theorem~\ref{boundary}, we have
$$ \delta_{\Ga_0(N)} (E_{\varphi, M,L}) = \beta_{\Ga_0(N),\varphi, M,L} D_{N,M,L}(\varphi), $$
where

$$
\beta_{\Ga_0(N), \varphi,M,L}  := -\frac{p^{k-r}M \phi(L)}{4n} \frac{ \tau(\varphi^{-1})}{\tau(\xi^{-1})} B_2(\xi^{-1}) \prod_{l \mid pM}  \Big( 1 - \frac{\xi(l)}{l^2} \Big) \in \Q(\zeta_{f}, \varphi):=K_{\varphi}
$$

and 
 $D_{N,M,L}(\varphi) \in  \Div^0(X_0(N)(\Q(\zeta_{f}, \varphi)) \otimes_\Z \sO_{K_{\varphi}}$. Define $T:=Mp^r$ and $\tilde{\beta}_{\Ga_0(N),\varphi,M,L} := T \beta_{\Ga_0(N), \varphi,M,L}$. Although $\tilde{\beta}_{\Ga_0(N), \varphi,M,L}$ does not need to be an element of $\Z[\zeta_f, \varphi]$, the following lemma shows that 
 \[
 \tilde{\beta}_{\Ga_0(N), \varphi,M,L} \in \Z[\frac{1}{6N}, \zeta_f, \varphi] \subset \mathcal{O}_{K_\varphi}[\frac{1}{6N}]:=R \subset K_\varphi \subset \CC.
 \]
As a consequence of Lemma~\ref{notsobad}, we work with $\Z[\frac{1}{6N}]$-coefficients. 
We define $\tilde{\beta}_{\Ga_0(N), \varphi,M,L}$ to detect the integrality statement as in Lemma~\ref{notsobad}, but over $R$,  $\tilde{\beta}_{\Ga_0(N), \varphi,M,L}$ and 
 $\beta_{\Ga_0(N), \varphi,M,L}$ generate the same ideal. 
Notice that $R$ is a localization of a Dedekind domain and hence a Dedekind domain contained in $\CC$. 
\begin{lemma}\label{notsobad}
 $24n\tilde{\beta}_{\Ga_0(N), \varphi,M,L} \in \Z[\zeta_f, \varphi].$ In-particular, $\tilde{\beta}_{\Ga_0(N), \varphi,M,L} \in \Z[\frac{1}{6p}, \zeta_f, \varphi]$
\end{lemma}
\begin{proof}
\begin{align*}
\tilde{\beta}_{\Ga_0(N), \varphi,M,L} & = -\frac{p^{k}M^2 \phi(L)}{4n} \frac{ \tau(\varphi^{-1})}{\tau(\xi^{-1})} B_2(\xi^{-1}) \prod_{l \mid Mp}  \Big( 1 - \frac{\xi(l)}{l^2} \Big)   \\
&=  -\frac{p^{k} \phi(L)}{4n^2}  \tau(\varphi^{-1})\tau(\xi) B_2(\xi^{-1}) \Big( 1 - \frac{\xi(p)}{p^2} \Big)  \prod_{l \mid M} \Big( l^2 - \xi(l)\Big).
\end{align*}
First, observe that for $l \mid M$, we have $\varphi(l) = \pm 1$ and hence $\xi(l)=1$. Note that the Gauss sums are elements of $\Z[\zeta_f, \varphi]$. Finally, observing that if $p \mid n$, then $\xi(p)=0$, we split the proof into two cases.

\begin{enumerate}
 \item[Case 1:] First, we consider the case where $\varphi$ is a quadratic character. 
Then $\xi= 1$ and $B_2(\xi^{-1}) = \frac{1}{6}$. 
\begin{align*}
24n\tilde{\beta}_{\Ga_0(N), \varphi,M,L} & = \tau(\varphi^{-1})  p^{k} \phi(L) \Big( 1 - \frac{1}{p^2} \Big)   \prod\limits_{l \mid M} ( l^2 - 1) \\
& =  p^{k-2} \phi(L)  \tau(\varphi^{-1})   \prod\limits_{l \mid pM} ( l^2 - 1) \Big) \in \Z[\zeta_f, \varphi].
\end{align*}

\item[Case 2:] We now consider case $\varphi^2 \neq 1$. In this situation, $\xi$ is a primitive Dirichlet character of conductor $n=p^t$, with $1 \le t \le r$. Since
$B_2(\xi^{-1}) = \frac{1}{6n} \sum\limits_{i=1}^{n-1} \xi^{-1}(i) (6i^2 - 6in +n^2)$, we re-write $24n\tilde{\beta}_{\Ga_0(N), \varphi,M,L}$ as $c_1 \cdot c_2$ with 
\begin{align*}
  c_1:=\tau(\varphi^{-1})\tau(\xi) \Big (\sum\limits_{i=1}^{n-1} \xi^{-1}(i) [ 6i^2 - 6in +n^2] \Big) \in \Z[\zeta_f, \varphi]\\
c_2:=  - \frac{p^{k} \phi(L)}{n^2}   \prod\limits_{l \mid M} ( l^2 - 1 )
\end{align*}
Now, $\nu_p(c_2)  \ge k -2t \ge k-2r \ge 0$ and hence $c_2 \in \Z$. Consequently, it follows that $24n \tilde{\beta}_{\Ga_0(N), \varphi,M,L} \in \Z[\zeta_f, \varphi]$. \qedhere
\end{enumerate}
\end{proof}
Recall that $R_{N}(E_{\varphi, M,L})$ is the $\Z$-submodule of $\mathbb{C}$ generated by the coefficients of $\delta_{\Ga_{0}(N)}(E_{\varphi, M,L}) $. 
Hence $R_{N}(E_{\varphi, M,L})\otimes_\Z R$ is an $R$-module contained in $\CC$. Observe that 
\[
\tilde{\beta}_{\Ga_0(N), \varphi,M,L} R=\beta_{\Ga_0(N), \varphi,M,L} R :=(\beta_{\Ga_0(N), \varphi,M,L}) \subset R  \subset \CC
\]  
\begin{lemma}
We have an equality of two $R$-modules $R_{N}(E_{\varphi, M,L})\otimes_\Z R=(\beta_{\Ga_0(N), \varphi,M,L})$. As a consequence, 
$R_{N}(E_{\varphi, M,L})\otimes_\Z R$ is a proper ideal of $R$. Hence, it is a free  $R$-module of rank $1$. 
\end{lemma}
\begin{proof}
By Theorem~\ref{boundary}, we write
$$ \delta_{\Ga_0(N)} (E_{\varphi, M,L}) = \beta_{\Ga_0(N),\varphi, M,L} D_{N,M,L}(\varphi). $$
The coefficients of $D_{N,M,L}(\varphi) \in  \Div^0(\partial(X_0(N)))\otimes K_\varphi$ are in $R$. Let the coefficients of $D_{N,M,L}(\varphi)$ generate an ideal $I$ of $R$. By the explicit expression of $I$, we have $I=R$.
 Hence, we have $
  R_{N}(E_{\varphi, M,L}) \otimes_\Z R \simeq \beta_{\Ga_0(N), \varphi,M,L} I  \simeq (\beta_{\Ga_0(N), \varphi,M,L})$.
 In turn, this implies $R_{N}(E_{\varphi, M,L}) \otimes_\Z R = (\beta_{\Ga_0(N), \varphi,M,L})$.
\end{proof}
For any $R$-module $N$, we denote by $N_\gP$ the localization at the prime ideal $\gP$. 
\begin{lemma}
We have an equality of two $R$-modules:
\[
 P_{N}(E_{\varphi, M,L}) \otimes_\Z R=
R + R_{N}(E_{\varphi, M,L}) \otimes_\Z R.
\]
\end{lemma}
\begin{proof}
As a first step, we now show that 
\[
\Lambda_{\pm}(E_{\varphi, M,L}, \chi, 1) \in R[\chi,\frac{1}{m_{\chi}}].
\]
By Proposition~\ref{specialvalues} for $\chi \in \x_{N}^{-\varphi(-1)}$, we have
\begin{align*}
T\Lambda_{\pm}(E_{\varphi, M,L}, \chi, 1) &=&
 \varphi(-m_{\chi}) \tau(\varphi) \chi(f)  \prod_{l \mid M }  (l- \chi\varphi(l) )  \prod_{q \mid L}(1- \chi \varphi^{-1}(q))  \frac{B_{1}(\chi^{-1}\varphi^{-1})}{2} \frac{B_{1}(\chi \varphi^{-1})}{2}\\
 &=&  \tau(\varphi) \cdot U;
\end{align*}
with 
\[
U=\varphi(-m_{\chi})  \chi(f)  \prod_{l \mid M }  (l- \chi\varphi(l) )  \prod_{q \mid L}(1- \chi \varphi^{-1}(q))  \frac{B_{1}(\chi^{-1}\varphi^{-1})}{2} \frac{B_{1}(\chi \varphi^{-1})}{2}.
\]
We first show that $U\in R[\frac{1}{m_\chi}]$. We have the following:
\begin{itemize}
\item 
$ \frac{ B_{1}(\varphi^{-1}\chi) }{2}, \frac{ B_{1}(\varphi^{-1}\chi^{-1}) }{2} \in \Z[\varphi, \chi, \frac{1}{m_{\chi}}] \subset R[\frac{1}{m_\chi}]$ (cf. ~\cite[Theorem 4.2(b)]{MR800251}).
\item 
$\varphi(-m_{\chi})  \chi(f)  \prod_{l \mid M }  (l- \chi\varphi(l) )  \prod_{q \mid L}(1- \chi \varphi^{-1}(q)) \in R$ by definition.
\end{itemize}
\item 
Since $\tau(\varphi) \in R$ (by the definition of Gauss sum), we have 
 \begin{equation*}
\Lambda_{\pm}(E_{\varphi, M,L}, \chi, 1) \in \frac{\tau(\varphi) }{T}\mathcal{O}_{K_\varphi}[\chi,\frac{1}{m_{\chi}}] \subset R[\chi,\frac{1}{m_{\chi}}].
 \end{equation*}
From Theorem~\ref{main}, we obtain the following:
\begin{equation}
\label{Pinclusion}
P:= P_{N}(E_{\varphi, M,L}) \otimes_\Z R\subseteq 
\widetilde{M}:= R + R_{N}(E_{\varphi, M,L})\otimes_\Z R.
\end{equation}
Now, $\widetilde{M}$ is a torsion-free $R$-module. 
We now show that this inclusion of the $R$-module is an equality.
Let $\gP$ be a prime ideal of $R$. 
We show that $P_{\gP}=\widetilde{M}_{\gP}$.  

Definitely $P_{\gP} \subset \widetilde{M}_{\gP}$. Notice that $ \widetilde{M}_{\gP}$ is a finitely generated $R_{\gP}$ module. 
Since $R$ is a Dedekind domain, $R_{\gP}$ is a discrete valuation ring, hence a principal ideal domain.
By the structure theorem for modules over principal ideal domains, 
we may assume that $\widetilde{M}_{\gP} \simeq R_{\gP}^s$. Since $P_{\gP}$ is a submodule of $\widetilde{M}_{\gP}$, we can write $P_{\gP}=R_{\gP}^{s'}$. From the inclusion~\ref{Pinclusion}, we have $s' \leq s$. We now show that $s' \geq s$. 

Definitely, $R_{N}(E_{\varphi, M,L}) \otimes_\Z R \subset P_{N}(E_{\varphi, M,L})\otimes_\Z R$. 
We first prove that there exists a character $\chi$ of conductor $m$ with $m \nmid T$
such that 
\[
R_{\gP}[\chi, \frac{1}{m_{\chi}}] \subset    P_{N}(E_{\varphi, M,L})\otimes_\Z R_{\gP}[\chi, \frac{1}{m_{\chi}}].
\]
Now, given a prime $\gP$ of $R$, we choose a nonquadratic character $\chi \in \x_N$ in a suitable way as follows:  

 \begin{enumerate}
 \item
For any choice of $\pm$, we can choose infinitely many Dirichlet characters $\chi \in \x_N^{\pm}$ such that $\gP \nmid (l-\chi(l) \varphi(l))$ for all primes $l \mid M$ and $\gP \nmid (1-\chi(q) \varphi(q))$ for all primes $q \mid L$. This essentially follows as in the proof of~\cite[Theorem 4.2(a)]{MR800251}. For the convenience of the reader, we give a sketch of the proof.   We choose a prime $m \in S_N$ large enough so that $\gP \nmid m$. We now specify the meaning of ``large".  

Since the $m^{th}$ roots of unity are distinct modulo $\gP$, it follows that for each $l$, $\gP \mid (l - \zeta_m^{a(l)} \varphi(l))$ for at most one choice of $a(l)\in \Z/m\Z$. Similarly, for each $q$, $\gP \mid (1- \zeta_m^{b(q)} \varphi(q))$ for at most one choice of $b(q) \in \Z/m\Z$. Let $a(1), \dots, a(s)$ be such that $\gP \mid (l - \zeta_m^{a(l)} \varphi(l))$ and $b(1), \dots, b(t)$ be such that $\gP \mid (1- \zeta_m^{b(j)} \varphi(q))$. Then there is a Dirichlet character $\chi \in \x_N^{\pm}$ of conductor $m$ such that $\chi(l) \neq \zeta_m^{a(i)}$ for $i \in \{ 1, \dots, s \}$ and $\chi(q) \neq \zeta_m^{b(j)}$ for $j \in \{ 1, \dots, t \}$. The result follows as there are infinitely many choices of $m \in S_N$.
  \item By~\cite[Theorem 4.2(c)]{MR800251}, we see that for all but finitely many choices of $\chi \in \x_N^{-\varphi(-1)}$, both $\frac{1}{2} B_{1}(\varphi^{-1}\chi)$ and $\frac{1}{2} B_{1}(\varphi^{-1}\chi^{-1})$ are 
 $\gP$-units provided $\gP \nmid \mrm{cond}(\chi)$.
 \end{enumerate}
Therefore, for any prime $\gP$, there exists a character $\chi \in \x_N^{-\varphi(-1)}$ such that
\begin{itemize}
\item 
$\gP \nmid (l_i-\chi(l) \varphi(l))$ for all $l \mid M$,
\item 
 $\gP \nmid (1-\chi(q) \varphi(q))$ for all $q \mid L$, 
\item 
$\frac{1}{2} B_{1}(\varphi^{-1}\chi)$ and $\frac{1}{2} B_{1}(\varphi^{-1}\chi^{-1})$ are $\gP$-units.
\end{itemize}
This implies that $\Lambda_{\pm}(E_{\varphi, M,L}, \chi, 1)$ is a $\gP$ unit.  By  ~\cite[p. 527, Lemma 2.2]{MR800251}, 
$\Lambda_{\pm}(E_{\varphi, M,L}, \chi, 1) \in  P_{N}(E_{\varphi, M, L})[\chi, \frac{1}{m_{\chi}}]$. 
Hence, we obtain the inclusion:
\[
R_{\gP}[\chi, \frac{1}{m_{\chi}}] \subset    P_{N}(E_{\varphi, M,L})\otimes_\Z R_{\gP}[\chi, \frac{1}{m_{\chi}}].
\]
Tensoring with $R_{\gP}[\chi, \frac{1}{m_{\chi}}]$ both sides of  $P_{\gP}=R_{\gP}^{s'} $,  we conclude that  
\[
R_{\gP}^{s'}[\chi, \frac{1}{m_{\chi}}]=P_{\gP} [\chi, \frac{1}{m_{\chi}}] \supset R_{N}(E_{\varphi, M,L}) \otimes_\Z R_\gP [\chi, \frac{1}{m_{\chi}}]+R_{\gP}[\chi, \frac{1}{m_{\chi}}].
\]
 Hence, we conclude $s' \geq s$. Since $P_\gP=\widetilde{M}_\gP$ for all prime ideals $\gP$ of $R$, it follows that $P_{N}(E_{\varphi, M, L})  \otimes_\Z R= \widetilde{M}$.

\end{proof}
\begin{proof}[Proof of Theorem \ref{ordercuspidal}]

We are interested in understanding the structure of the $R$-module $ A_{N}(E_{\varphi,M, L})$. We have 
\begin{align*}
 A_{N}(E_{\varphi,M, L})  \otimes_\Z R & =& \left(P_{N}(E_{\varphi, M,L})/R_{N}(E_{\varphi, M,L}) \right)  \otimes_\Z R  \\
                    & = & P_{N}(E_{\varphi, M,L}) \otimes_\Z R/\left(R_{N}(E_{\varphi, M,L}) \otimes_\Z R \right) \\
                       & = & \left( R + R_{N}(E_{\varphi, M,L}) \otimes_\Z R\right)/ \left( R_{N}(E_{\varphi, M,L}) \otimes_\Z R \right)\\
                        & \cong &  R /\left( R_{N}(E_{\varphi, M,L}) \otimes_\Z R \right)\\
                          & \cong &  R/(\beta_{\Ga_0(N), \varphi,M,L}).\\
                  \end{align*}
Let $J':=(\beta_{\Ga_0(N), \varphi,M,L})$ be an ideal of $R$.
We then obtain the isomorphisms 
\[
R/J' \simeq A_{N}(E_{\varphi, M, L}) \otimes_\Z R \simeq \bigoplus\limits_{i=1}^n R/p^{n_i}R.
\]
Thus, we obtain isomorphisms $C_{\Ga_0(N)}(E_{\varphi, M,L}) \otimes_\Z R \simeq \bigoplus\limits_{i=1}^n R/p^{n_i}R \simeq R/J'$.
\end{proof}
We now prove Corollary~\ref{support}. 
\begin{proof}
A prime $l$  divides the order of a finite abelian group $C$ if and only if $l$ divides the order $|C \otimes_\Z R|$ of the $R$-module $C \otimes_\Z R$. 
\end{proof}
\begin{remark}
\label{Shimuracuspidal}
In \cite[p.~538]{MR800251}, the computation of cuspidal subgroups is carried out by Stevens using Shimura subgroups. For general $N$, one can compute the order of the Shimura subgroup~\cite[Corollary~1]{MR1141458}. However, we do not know of any clean way to compute the kernel in the short exact sequence \cite[\S 4.3]{MR800251}. This is the reason we compute the order of the cuspidal subgroup directly in this paper. To complete this computation, we need to perform a careful and detailed analysis of the Eisenstein series and the associated cuspidal subgroup. 

\end{remark}
\section{Ribet's conjecture for non-rational Eisenstein series}
\label{Eisenstein ideals}
Let $\mathcal{T}(N)$ be the Hecke algebra acting on $M_2(\Ga_0(N);\CC)$. For an Eisenstein series $E \in E_2(\Ga_0(N);\CC)$, let $\mathrm{Ann}_{\mathcal{T}(N)}(E)$ be the annihilator of the Eisenstein series $E$. 
 The Eisenstein ideal is the image of $\mathrm{Ann}_{\mathcal{T}(N)}(E)$ in the Hecke algebra $\T(N)$ (the Hecke algebra acting on the space of cusp forms).  We now recall our Theorem~\ref{classification} for the convenience of the readers. Given an integer $N$, write $N=N_1N_2$ as in equation \ref{decomp}. 
 Let $\m$ be a non-rational Eisenstein maximal ideal of $\T(N)$ whose residual characteristic $\ell$ is coprime to $6N_1$.
The ideal $\m$ determines a non-trivial Dirichlet character
$\overline{\epsilon}_\m$ of the conductor $f$ (with $f^2 \mid N_1$). Moreover, there exist square-free natural numbers $t$ and $M$ such that
  $$T_r \equiv \overline{\epsilon}_\m(r) + r \cdot (\overline{\epsilon}_\m)^{-1}(r) \pmod{\m} \text{ \qquad\qquad for primes } r \nmid N,$$ 
    $$U_p \equiv 0 \pmod \m \text{ \qquad\qquad for primes } p \mid ft \text{ with }  t \mid N_1 \text{ and }  (t,f)=1,$$ 
    $$U_s \equiv  s \cdot (\overline{\epsilon}_\m)^{-1}(s) \pmod{\m} \text{ \qquad\qquad for primes } s \mid M \text{ with } M \mid N \text{ and } (M,ft)=1,$$
        $$U_q \equiv  \overline{\epsilon}_\m(q) \pmod{\m} \text{ \qquad\qquad for all other primes } q.$$ 
        In particular, the maximal ideal $\m$ 
with $\kappa(\m) \cong \T/\m\cong \F_\ell[\overline{\epsilon}_\m]$
contains an ideal of the form
$$I_{\overline{\epsilon}_\m,M, t}(N)= \Big(I_{\overline{\epsilon}_\m}(N), U_p \text{ for primes } p \mid ft ,  g_s(U_s) \text{ for primes } s \mid M, h_q(U_q)   \text{ for all other primes } \Big). $$ 
Hence, we have $\m = (\ell, I_{\overline{\epsilon}_\m,M, t}(N))$. We are interested in finding the possible $\ell$ away from $6N$. To do the same, we wish to explicitly write down an Eisenstein series associated with $\m$.

 Consider the  ring $\T(N)/ I_{\overline{\epsilon}_\m,M,t}(N)$
 \cite[p. 704]{MR4358260}. We write for every prime $q$ (see  loc.\ cit.):
\[
\left(\T(N)/ I_{\overline{\epsilon}_\m,M,t}(N)\right) \otimes \Z_q \simeq \prod_{\substack{\m}} \T(N)_{\m}/I_{\overline{\epsilon}_\m,M,t}(N) 
\]
with maximal ideals $\m$ containing $(q,I_{\overline{\epsilon}_\m,M,t}(N))$.

By Theorem ~\ref{classification} as recalled above, it follows that ideals of the form $(q, I_{\overline{\epsilon}_\m,M,t}(N))$ are maximal, and hence for a prime $q$ we have 
$| \big(\T(N)/ I_{\overline{\epsilon}_\m,M,t}(N)\big)[q^\infty] |\neq 0$ if and only if $\m = (q, I_{\overline{\epsilon}_\m,M,t}(N))$  is an Eisenstein maximal ideal, provided $q \nmid 6N_1$.
We understand the support of the index of the Eisenstein ideal $I_{\overline{\epsilon}_\m, M, t}$  and study it by relating this index to the orders of the cuspidal subgroup associated with an Eisenstein series $E_\m$.  We describe the process in the following.

We now choose a Teichm\"uller lift  $\epsilon_\m$ of 
$\overline{\epsilon}_\m$ in a controlled way. 
 Assume that the non-trivial character $\overline{\epsilon}_\m : (\Z/f \Z)^\times \to \kappa(\m)^\times$ is given by $\overline{\epsilon}_\m (a_i) = \alpha_i$, where $a_1, \dots, a_r$ is a set of generators of $(\Z/f\Z)^\times$.  If $k \ge 2$ is of the order $\overline{\epsilon}_\m$, then $k \mid (|\kappa(m)|-1)$. In particular, we have $(k, \ell)=1$. Let $\zeta_k$ be a complex number that is a fixed primitive $k$ -th root of unity; it follows that there exists a prime $\cL$ in $\Z[\zeta_k]$ lying above $\ell$ and an isomorphism $\Theta: \Z[\zeta_k]/(\cL) \simeq \F_\ell[\overline{\epsilon}_\m]$. Assume that $\Theta(\zeta_k)= \alpha \in \F_\ell[\overline{\epsilon}_\m]$ (say). Since $\alpha_i^k=1$ for all $i$, it follows that $\alpha_i = \alpha^{b_i}$ for some integer $b_i$.

 Now, the Teichm\"uller lift $\epsilon_\m$ is the character $\epsilon_\m : (\Z/f\Z)^\times \to \CC^\times$ given by $\epsilon_\m(a_i)= \zeta_k^{b_i}$. Then, it follows that $\epsilon_\m$ is a Dirichlet character of conductor $f$ such that for all primes $r \nmid f$, $\epsilon_\m(\Frob_r) \equiv \overline{\epsilon}_\m(\Frob_r) \pmod{\cL}$. The Teichm\"uller lift $\epsilon_\m$  is unramified at the prime $\ell$. The next definition aims to uniquely determine an Eisenstein series associated with an Eisenstein maximal ideal. 
 
\begin{definition}[Eisenstein series associated to Eisenstein maximal ideal]\label{Em}
 Let $N=N_1N_2$ be an integer as in Equation~\ref{decomp} and $\m$ be a non-rational Eisenstein maximal ideal of $\T(N)$ (with the associated non-trivial Dirichlet character $\overline{\epsilon}_\m$ of the conductor $f$) with the residual characteristic $\ell$ coprime to $6N_1$. Let $\epsilon_\m: (\Z/f\Z)^\times \to \CC^\times$ be the lift of $\overline{\epsilon}_\m$ as described above. By Remark~\ref{classificationremark}, we have $\m =( \ell, I_{\overline{\epsilon}_\m, M, t}(N))$. By the Eisenstein series associated with $\m$, we mean the non-rational Eisenstein series $E_\m :=E_{\epsilon_\m, Mt, Lt} \in E_2(\Ga_0(N), \CC)$, where $L= \prod\limits_{ \substack{q \mid N,\\  q \nmid Mft}} q$. 
\end{definition}
Recall that this Eisenstein series $E_{\varphi, Mt , Lt}$ (cf. the description just before Proposition~\ref{lvalue}) is an eigenform for all Hecke operators, and the eigenvalues are given by
\begin{align*}
T_{r}(E_{\varphi, Mt, Lt}) &=  (\varphi(r)+ r \varphi^{-1}(r)) E_{\varphi, Mt,Lt} \quad  \text{for}\ r \nmid N, \\
U_{l}(E_{\varphi, Mt, Lt}) &= l \varphi^{-1}(l) E_{\varphi, Mt, Lt} 
\qquad \qquad \quad \text{for}\  l \mid M, \\
U_{q}(E_{\varphi, Mt, Lt}) &=  \varphi(q) E_{\varphi, Mt, Lt}  
\qquad \qquad \qquad  \text{for}\  q \mid L ,  \\
U_{p}(E_{\varphi, Mt, Lt}) &=  U_{p}(E_\varphi) =0  
\qquad \qquad \qquad \quad \text{for}\   p \mid ft.  
\end{align*}

Let $K=\Q(\epsilon_\m)$, $S= \Z[\frac{1}{6N}, \epsilon_\m]$ and $\T_S(N):= \T(N) \otimes_{\Z} S$.
We are interested in the Eisenstein ideals $\widetilde{I_{\epsilon_\m,M,t}}(N)$, that is, the image in $ \T_S(N)$ of $\mathrm{Ann}_{\mathcal{T}_S(N)}(E_{\m})$.   

Consider the ideal $\widetilde{I'_{\epsilon_\m,M,t}}(N)$ of $\T_S(N)$ generated by the following elements:
\begin{enumerate}
\item $T_r - \epsilon_\m(r) - r \epsilon_\m^{-1}(r)$  for primes  $r \nmid N$,
\item $U_p$ for primes $p \mid ft$,
\item $U_s - s \epsilon_\m^{-1}(s)$ for primes $s \mid M$,
\item $U_q - \epsilon_\m(q)$ for primes $q \mid L$.
\end{enumerate}
Definitely $T_S(N) \supset \widetilde{I_{\epsilon_\m,M,t}}(N) \supset \widetilde{I'_{\epsilon_\m,M,t}}(N)$. 
\begin{lemma}\label{faliureofRamanujan2}
There exists a non-zero ideal $J$ of $\Z[\epsilon_\m]$ such that
$$ \T_S(N)/ \widetilde{I_{\epsilon_\m,M,t}}(N) \cong S/J.$$
In particular,  the index of the Eisenstein ideal $\T_S(N)/\widetilde{I_{\epsilon_\m,M,t}}(N) $ (up to the power of $2, 3$ and the primes dividing $N$) is given by the norm $N(J)$ of the ideal $J$.
\end{lemma}
\begin{proof}
We have an inclusion $S  \subset 
\T_S(N)$. This induces a surjective ring homomorphism
$$\kappa: S = \Z[\frac{1}{6N}, \epsilon_\m] \to \T_S(N)/\left(\widetilde{I_{\epsilon_\m,M,t}}(N) \right) \text{ given by } z \mapsto z \pmod{\widetilde{I_{\epsilon_\m,M,t}}(N)}.$$
Recall that this map is surjective since, modulo $\widetilde{I_{\epsilon_\m,M,t}}(N)$ ( and hence modulo $\widetilde{I'_{\epsilon_\m,M,t}}(N)$), every $T_r$ comes from an element of $S$. 

 Suppose $\T_S(N)/\left(\widetilde{I_{\epsilon_\m,M,t}}(N) \right) \cong S$. 
 From the description of $\widetilde{I'_{\epsilon_\m,M,t}}(N)$ 
 and again from the non-degenerate bilinear pairing 
 \cite[(10), p. 465]{MR1047143}, there exists a weight $2$ Hecke eigenform that is also a cusp form  over $S \subset \mathbb{C}$ with $T_r$ eigenvalue $a_r(f)=\epsilon_\m(r) + r \epsilon_\m^{-1}(r)$ for all primes $r \nmid N$ and hence, in particular, for some prime $r \ge 7$. 
 Note that the Deligne (Ramanujan) bound is given by $|a_r(f)|^2 \le 4r$. Writing $\epsilon_\m(r)= x+iy$ with $x^2+y^2=1$, the Ramanujan bound gives us $$ 1+r^2+ 2r(2x^2-1) \le 4r \text{ or equivalently } (1-r)^2 \le 4r(1-x^2) \le 4r.$$
 Note that as $r \ge 7$, $(1-r)^2 > 4r$, hence the cusp form does not satisfy the Ramanujan bound. 
 As a consequence, we obtain $\T_S(N)/\widetilde{I_{\epsilon_\m,M,t}}(N) \cong S/J$ for a non-zero ideal $J$.
\end{proof}

 For the remainder of this section, we assume that $N$ satisfies the hypothesis of Theorem~\ref{ordercuspidal}. Assume that $N=p^kN_2$, with $k \ge 2$ and $N_2$ a square-free natural number such that all prime factors $\widetilde{q}$ of $N_2$ satisfy $\widetilde{q} \equiv \pm 1 \pmod{p^{\lfloor 
\frac{k}{2} \rfloor}}$.

Recall that the cuspidal subgroup associated with $E_\m$ is the finite subgroup of $J_0(N)(\mathbb{Q}(\zeta_f,\varphi))$, and we have calculated the order of $C_{\Ga_0(N)}(E_\m)\otimes R$ in Theorem \ref{ordercuspidal}.  We now relate $|C_{\Ga_0(N)}(E_\m)\otimes R|$ to the cardinality of a suitable quotient of $\T_R(N) := \T(N) \otimes_\Z R$. 
Since $t=1$ by our assumption on $N$, we denote $\widetilde{I_{\epsilon_\m,M}}:=\widetilde{I_{\epsilon_\m,M,1}}$ (respectively $I_{\epsilon_\m,M}:=I_{\epsilon_\m,M,1}$).

For a prime $q$, we define $\alpha(q), \beta(q) \in \N$ as follows:
$$ \nu_q(|\T_R(N)/\left(\widetilde{I_{\epsilon_\m,M}}\T_R(N) \right)|)  = \alpha(q) \text{ and } \nu_q(|C_{\Ga_0(N)}(E_\m)\otimes R|)=\beta(q).$$

We show that $\alpha(q) = \beta(q)$ if $q \nmid 6N $. As we know $\beta(q)$, this will also give $\alpha(q)$ (the primes that produce Ramanujan-like congruences). 

Recall that $C^{Ren}_{\Ga_0(N)}(E)$ is the cuspidal subgroup associated with $E$ as defined by Ren
(cf.\S~\ref{Cuspidal}).
  
\begin{lemma}\label{annihilates}
The cuspidal subgroup $C^{Ren}_{\Ga_0(N)}(E_{\m})  \otimes R$ is a $\T_R(N)/\left(\widetilde{I_{\epsilon_\m,M}}(N)\right)$ module. 
 \end{lemma}
 \begin{proof}
Note that $T_l \circ w_N=w_N \circ T_l^{t}$ (cf. Ribet. \cite{MR1047143}).  It follows $\delta_{\Ga_0(N)}(T_l E_{\m}) = T^{t}_l \delta_{\Ga_0(N)}(E_{\m})$ (cf. \cite[p. 110]{MR670070}, 
see also~\cite[Lemma 1.11]{MR1376558}). As a consequence, $\widetilde{I_{\epsilon_\m, M}}(N)$ annihilates $C^{Ren}_{\Ga_0(N)}(E_{\m}) \otimes R$ from the very definition of $\widetilde{I_{\epsilon_\m, M}}(N)$.
\end{proof}
By our main Theorem~\ref{ordercuspidal}, we have an isomorphism of $R$-modules:
\[
C_{\Ga_0(N)}(E_\m) \otimes R \simeq R/(\beta_{\Ga_0(N), \varphi,M,L}).
\]
As an $R$-module, $C_{\Ga_0(N)}(E_\m)\otimes R$ is a quotient of $R$ by the ideal $J':= (\beta_{\Ga_0(N), \varphi,M,L})$. As a consequence, we have the $R$-module isomorphism
\[
\operatorname{End}_R \left(C_{\Ga_0(N)}(E_\m)\otimes R\right) \simeq C_{\Ga_0(N)}(E_\m)\otimes R \cong R/J'. 
\]

The lemma~\ref{annihilates} shows that $\widetilde{I_{\epsilon_\m,M}}\T_R(N)$ annihilates $C^{Ren}_{\Ga_0(N)}(E_\m)  \otimes R$.
We see that $C^{Ren}_{\Ga_0(N)}(E_\m)  \otimes R$ is a $\T_R(N)/\left(\widetilde{I_{\epsilon_\m,M}}\T_R(N)\right) \cong R/JR$ module. This implies $J (C^{Ren}_{\Ga_0(N)}(E_\m) \otimes R )=0$.  

Hence, there exists a surjection of finite quotient rings:
\[
\T_R(N)/\left(\widetilde{I_{\epsilon_\m,M}}\T_R(N)\right) \twoheadrightarrow  \mathrm{End}_R(C^{Ren}_{\Ga_0(N)}(E_\m)  \otimes R).
\]
Note that, as an $R$-module, $ C^{Ren}_{\Ga_0(N)}(E_\m)\otimes R$ is cyclic and $|C_{\Ga_0(N)}(E_\m)\otimes R|=| C^{Ren}_{\Ga_0(N)}(E_\m)\otimes R|$. We deduce that
$$ |C_{\Ga_0(N)}(E_\m)  \otimes R|= |C^{Ren}_{\Ga_0(N)}(E_\m)  \otimes R|=|\mathrm{End}_R(C^{Ren}_{\Ga_0(N)}(E_\m)  \otimes R)| \leq    |\T_R(N)/\left(\widetilde{I_{\epsilon_\m,M} \T_R(N)}\right)|= |R/\left(JR\right)|.$$
We conclude that $\alpha(q) \ge \beta(q)$ for all $q \nmid 6N$.

Denote by $\cX_0(N)$ the Deligne-Rapoport smooth model of $X_0(N)$ over $\Z[\frac{1}{N}]$  \cite{MR0337993} (see also \cite[\S 1.1]{MR3304683}). 
For any ring $V$, let $M_2^{A}(\Ga_0(N), V)$ be the space of algebraic modular forms in the sense of Katz. Let $M_2^B(\Ga_0(N),V)$ (respectively $S_2^B(\Ga_0(N),V)$  be the space of $p$-adic modular forms (respectively cusp forms) in the sense of Serre with coefficients in a ring $V$ (cf.~\S~\ref{maximalideals}).  We use this crucially in the next theorem. For brevity, we use the same ``$q$" for the primes and $q$-expansions. However, it will be clear from the context.  

 Recall the following result from \cite[Proposition 1.4.9]{MR3304683}. 
For any $\Z[\frac{1}{N}]$ algebra $V$, we have the following inclusion:
\begin{equation}
\label{comparemodular}
M_2^B(\Ga_0(N),V) \xhookrightarrow{} M_2^A(\Ga_0(N),V).
\end{equation}

\begin{theorem}\label{index=order}
 Let $\m$ be a non-rational Eisenstein maximal ideal of $\T(N)$ whose residual characteristic $\ell$ is coprime to $6N$, and let $R$ be the ring as in the introduction.
  For any prime $q \nmid 6N$, we have  $$  \alpha(q)=\beta(q).$$

\end{theorem}

\begin{proof}
To show $\alpha(q)= \beta(q)$, it is enough to show $\alpha(q) \le \beta(q)$ for $q \nmid 6N $. This is obvious if $\alpha(q)=0$. We now assume that $\alpha(q) \ge 1$ or, equivalently
, $q \mid |\T_R(N)/\widetilde{I_{\epsilon_\m,M}} \T_R(N)| = |R/JR|$.

Write $JR = J_{q}J_{\backslash q}$, where the prime factorization of the ideal $J_q$ (resp. $J_{\backslash q}$) only contains primes of $R$ that lie above $q$ (resp. do not lie over $q$).

Consider the ideal $\widetilde{I} := (J_q, \widetilde{I_{\epsilon_\m,M}}) \subset \T_R(N)$. Since 
\[
\T_R(N)/\widetilde{I} \cong \left(R/JR\right)/ J_q (R/JR)\cong R/J_qR,
\]
we see  
\[
\nu_q(|R/J_qR|) =\nu_q( |\T_R(N)/ \widetilde{I}|)=\nu_q( |\T_R(N)/  (J_q, \widetilde{I_{\epsilon_\m,M}})|) =\nu_q( | \T_R(N)/\left(\widetilde{I_{\epsilon_\m,M}} \T_R(N)\right)|) = \alpha(q).
\]
By  \cite[(10), p. 465]{MR1047143}, 
we have a non-zero cusp form $F \in  S^B_2(\Ga_0(N),\T_R(N)/ \widetilde{I}) = S_2^B(\Ga_0(N), R/\left(J_q R)\right) $  whose Fourier expansion at the cusp $i\infty = 
\bstwomat{1}{N}_{\Ga_0(N)}$ is given by
$$F= \sum_{n \ge 1} (T_n \mod{\widetilde{ I}})\ q^n.$$
This is a modular form with coefficients at $\infty$ in the ring $R/(J_q R)$. 

Note that the Eisenstein series $E_\m$ associated with our Eisenstein maximal ideal $\m$ (cf. Definition~\ref{Em}) is an element of $M_2^B(\Ga_0(N),R)$.  Thus, we have a  modular form 
       $$G= F - \overline{E_\m}  \in M^B_2(\Ga_0(N), R/(J_q R)),$$ 
 where $\overline{ E_\m }$ denotes the reduction of $E_\m$ modulo the ideal $J_qR$ of $R$.

     Since $q \nmid 6N$, we have $R/(J_qR)$ 
as a $\Z[\frac{1}{N}]$ algebra. The special fiber $\cX_0(N)/\overline{\F_q}$ is connected since $q \nmid N$ \cite[p. 227, Corollary 5.6]{MR0337993}. Using the inclusion~\ref{comparemodular}, we consider $G$ as an element in the space $M_2^A(\Ga_0(N), R/(J_qR))$ of algebraic modular forms in the sense of Katz. We now apply the $q$-expansion principle to this $G$.  The map 
\[
M_2^{A}(\Ga_0(N),R) \rightarrow R \otimes \Z[[q]];
\]
Obtained by evaluating $G$ at the Tate curve is injective by Katz~\cite[Theorem 1.6.1, \S 1.6]{MR0447119} (See also Brian Conrad's notes on $q$ expansion \cite[Lemma 3.2]{Conrad}  that works fine since we use only the fact that the special fiber of the moduli stack is connected as $q \nmid 6N$). Since the $q$ expansion of $G$ at the cusp $i\infty$ is $0$, the modular form $G$ is identically $0$.  Hence, the constant term in the Fourier expansions of $E_\m$ at any cusp is also $0 \in R/(J_qR)$.

We now wish to compute the Fourier coefficients of the Eisenstein series $E_\m$ at the cusp $\bstwomat{1}{f}_{\Ga_0(N)}$. 
Note that
\[
\delta_{\Ga_0(N)}(E_\m) = \displaystyle\sum_{y \in 
 \partial(X_0(N)) } e_{\Ga_0(N)}(y) \cdot a_0(E\mid_{\theta_y}) \{y\} = \beta_{\Ga_0(N), \epsilon_\m, M, L} D_{N, M, L}(\epsilon_\m).
 \]
 We denote $a_0(E\mid_{\theta_y})$ by $a_0(E\mid_{y})$, for simplicity.
 Comparing the coefficients on both sides at the cusp 
 $\bstwomat{1}{f}_{\Ga_0(N)}$, we get 
$$e_{\Ga_0(N)}(\bstwomat{1}{f}_{\Ga_0(N)}) a_0(E_{\epsilon_\m, M, L}\mid \bstwomat{1}{f}_{\Ga_0(N)})  = \beta_{\Ga_0(N), \epsilon_\m, M, L} \cdot U'. $$
 Here,  $U' \in R^{\times}$ is a unit by our assumption on $N$ (cf. Theorem~\ref{boundary}).  
 We again use the fact that $e_{\Ga_0(N)}(\bstwomat{1}{f}_{\Ga_0(N)}) = \frac{N}{f^2}$. Since $N=p^k N_2$ and $f=p^r$, after simplification, we obtain $$ a_0 (E_ {\epsilon_\m, M, L}\mid \bstwomat{1}{f}_{\Ga_0(N)}) = \frac{\beta_{\Ga_0(N), \epsilon_\m, M, L} \cdot  U'}{p^{k-2r}N_2}.$$    
 
   Since the Fourier coefficient of $E_\m \pmod{J_q}$ at cusps  $\bstwomat{1}{f}_{\Ga_0(N)}$ is $0$, we see that 
 \[
 \mrm{Num}(a_0(E_\m \mid \bstwomat{1}{f}_{\Ga_0(N)})) := (a_0(E_\m \mid \bstwomat{1}{f}_{\Ga_0(N)}) \cap R \subset J_qR.
 \]
 As a consequence, we obtain,
 $$ \mrm{Num}( \beta_{\Ga_0(N), \epsilon_\m, M, L})  \subset \mrm{Num}(a_0(E_\m \mid \bstwomat{1}{f}_{\Ga_0(N)}) \subset J_qR.$$ 
 Finally, we obtain our desired result $\beta(q) \geq \alpha(q)$ from the 
following inequality:
\[
\beta(q) =\nu_q( |C_{\Ga_0(N)}(E_\m)\otimes R|) =\nu_q( |R/(\Num\Big( \beta_{\Ga_0(N), \epsilon_\m, M, L}\Big))|) \ge \nu_q(|R/(J_qR)|) = \alpha(q).
\]

\end{proof}

Now we are ready to give a proof of Theorem \ref{mainthmRibet}.\\

\begin{proof}[Proof of Theorem \ref{mainthmRibet}]
Let $\ell \nmid 6N$ be the residual characteristic of a non-rational Eisenstein maximal ideal $\m$ of $\T(N)$. Since $N=p^kN_2$ with $N_2$ square-free, by Theorem~\ref{classification}, $\m= (\ell, I_{\overline{\epsilon}_\m,M})$. Now 
\[
 \T_S(N)/(\cL, \widetilde{I_{\epsilon_\m,M}}) \cong \left(\Z[\epsilon_\m]/J\right)  \otimes \Z[\frac{1}{6N}]/\left(\cL(\frac{\Z[\epsilon_\m]}{J}  \otimes \Z[\frac{1}{6N}])\right).
\]
The last isomorphism follows from Lemma~\ref{faliureofRamanujan2}.
As a consequence, we have $\ell \mid |\Z[\epsilon_\m]/J|$ and therefore $\ell \mid  |\T_R(N)/\left(\widetilde{I_{\epsilon_\m,M}}\T_R(N)\right)|$.

By Theorem~\ref{index=order}, it follows that $\ell \mid |C_{\Ga_0(N)}(E_\m)\otimes R |$ and hence $\ell \mid |C_{\Ga_0(N)}(E_\m) |$. In other words, $C_{\Ga_0(N)}(E_\m)[\ell] \neq \{0\}$. 
Now $C_{\Ga_0(N)}[\m] = C_{\Ga_0(N)}[(\ell, I_{\overline{\epsilon}_\m,M,t})] \supseteq C_{\Ga_0(N)}(E_\m)[\ell] \neq \{0\}$ as the cardinality of the last group is divisible by $\ell$.
\end{proof}

Given an integer $N=p^kN_2$, it follows that the residual characteristic of a non-rational Eisenstein ideal $\m$ of $\T(N)$ must divide $6N \cdot |C_{\Ga_0(N)}(E_\m)|$. 
Thus, determining the prime divisors of $|C_{\Ga_0(N)}(E_\m)|$ will identify which primes may appear as a residual characteristic of some non-rational Eisenstein ideal of $\T(N)$. To do this, we define the following two subsets of rational primes as follows.
\[
S_1(N) := \{ r\in \N:  r \mid (q^2-1) \text{ for some prime } q \mid N \}, S_2(N) = \bigcup_\varphi S_2(N, \varphi).
\]
Here, the union is over the set of nontrivial Dirichlet characters $\varphi$ of conductor $f$ with $f^2 \mid p^k$,
 $$S_2(N, \varphi) := \{s \in \N:  s \mid N_{\Q(\varphi)/\Q}\big(6nB_2(\xi^{-1}) \},$$
where $\xi$ is the primitive Dirichlet character associated with $\varphi^2$. Note that $6nB_2(\xi^{-1}) \in \Z[\varphi]$.
 \begin{proposition}
 \label{possiblecharacteristic}
 Let $N=p^k N_2$ be an integer as in Theorem \ref{theorem1}  and $\m$ be a non-rational Eisenstein maximal ideal of $\T(N)$. The residual characteristic $\ell$
 of $\m$ belongs to the set 
 \[
 D_R:=\{ q \mid 6N \} \cup S_1(N) \cup S_2(N) \subset \N .
 \]
 \end{proposition}
 
 \begin{proof}
 If $\ell \mid 6N$, then there is nothing to prove.  We now assume that $\m$ is a non-rational Eisenstein ideal of $\T(N)$ whose
 residual characteristic $\ell$ is coprime to $6N $. 
 In this case, $\ell$ must divide $|C_{\Ga_0(N)}(E_\m) \otimes R|= |R/(\tilde{\beta}_{\Ga_0(N), \epsilon_\m, M, L})|$. Recall that
 
 \begin{equation*}
\tilde{\beta}_{\Ga_0(N), \epsilon_\m,M,L} =
\begin{cases}
 -\tau(\epsilon_\m^{-1}) 
  \frac{p^{k-2} \phi(L)}{24} \prod\limits_{q \mid Mp} ( q^2 - 1)&  \text{if $\epsilon_\m$ is the quadratic character,} \\
 -\tau(\epsilon_\m^{-1})\tau(\xi) \frac{p^{k} \phi(L)}{24n^3} (6n B_2(\xi^{-1}))  \prod\limits_{q \mid M} ( q^2 - 1 ) , &  \text{Otherwise. } \\
\end{cases}
\end{equation*}
 Here, $\xi$ is the primitive Dirichlet character of conductor $n$ corresponding to $\epsilon_\m^2$.

      Case $1$: If $\epsilon_\m$ is a quadratic character, then
$\ell \mid |R/\left((\tilde{\beta}_{\Ga_0(N), \epsilon_\m, M, L})\right)|$ if and only if
$\ell$ divides $|R/\left( \phi(L) \prod\limits_{q \mid Mp} ( q^2 - 1 ) \right)|$, which implies
that $\ell \in S_1(N)$. 

Case $2$:-
Similarly, if $\epsilon_\m$ is not quadratic, then $ \ell \mid  |R/\left((\tilde{\beta}_{\Ga_0(N), \epsilon_\m, M, L})\right)|$ if and only if
$$
\ell \mid |R/\left(  \phi(L)  6nB_2(\xi^{-1}) 
\prod\limits_{q \mid M} (q^2 - 1 ) \right) |. 
$$ 
Since $L$ is square-free, we have $\ell \in S_1(N) \cup S_2(N)$.
\end{proof}

 \section{Explicit congruences between non-rational Eisenstein series and cusp forms }
\label{numerical examples}
 Let $N=p^kN_2$ be an integer as in Theorem~\ref{ordercuspidal}. By Proposition \ref{possiblecharacteristic}, the possible residual characteristics $\ell$
 of $\m$ belong to $\{ p \mid 6N \} \cup S_1(N) \cup S_2(N)$.  Theorems~\ref{mainthmRibet} and Theorem ~\ref{index=order} show that if there exists a non-rational Eisenstein ideal $\m$ with residual character $\ell$ coprime to $6N$, then there exists an Eisenstein series $E_{\epsilon_\m,M,L}$ associated to $\m$ and there exists a cusp form $f \in S_2(\Ga_0(N))$ such that 
 $$ E_{\epsilon_\m,M,L}  \equiv f \pmod{ \cL} $$ 
 for some prime $\cL \mid \ell$.

 For a prime $q \mid | C_{\Ga_0(N)}(E_{\varphi, M, L})| $ with 
 $q \nmid 6N$, we give examples of explicit congruence between $E_{\varphi,M,L}$ and $f \in S_2(\Ga_0(N))$ above a prime over $q$. 
 We compute these examples using SAGE \cite{sagemath} and LMFDB \cite{lmfdb}.

 \begin{example}
 Take $N=121=11^2$ and $\varphi$ is the quadratic character of conductor $11$ (given by $\varphi(2)=-1$). Then $\tau(\varphi^{-1}) = \sqrt{-11}$, $\xi = 1$, $n=1$, $\tau(1)=1$, $B_2(1) = \frac{1}{6}$. We have
 $$\tilde{\beta}_{\Ga_0(N),\varphi} = -5 \sqrt{-11} .$$
 So, $5$ is a possible residual characteristic of an Eisenstein ideal $\m$ whose associated Eisenstein series is $E_\varphi$. The $q$-expansion of $E_\varphi =\sum\limits_{ n \ge 1} a_n q^n$, where  $a_n = \sum\sum\limits_{bc=n} \varphi(c) \varphi^{-1}(b)b$, we get
 \begin{align*}
 E_\varphi = \sum\limits_{n=1}^\infty \varphi(n) (\sum\limits_{d \mid n} d) q^n\\
  = q - 3q^2 + 4 q^3 + 7 q^4 + 6 q^5 - 12 q^6 - 8 q^7 -15 q^8 + 13 q^9 -18 q^{10}+ 28q^{12} +O(q^{13}).
 \end{align*}
 Let $ f \in S_2^{\mrm{new}}(121)$  be a newform with LMFDB label $121.2.a.d$ and its $q$-expansion is given by
 $$q + 2q^{2} - q^{3} + 2q^{4} + q^{5} - 2q^{6} + 2q^{7} - 2q^{9} + 2q^{10} - 2q^{12}+ O(q^{13}).$$
We obtain 
$$E_\varphi(z) \equiv f(z) \pmod{5}.$$
The corresponding non-rational Eisenstein maximal ideal $\m$ of characteristic $5$, with $\T/\m \cong \F_5$, is given by
 $$\m= \langle 5, U_{11}, \{ T_r -1 -r \}_{\text{ primes } r \equiv 1, 3 , 4 , 5 , 9 \pmod{11}}, \{ T_s +1+s \}_{\text{ primes } s \equiv 2, 6 , 7 , 8 , 10 \pmod{11}} \rangle.$$
\end{example}

\begin{example}
Take $N=121= 11^2$ and $\varphi$ is the Dirichlet character of conductor $11$, of order $10$, given by $ \varphi(2) = \zeta_{10}$.  Then, $\xi= \varphi^2$ is a primitive Dirichlet character of conductor $11$ with $\xi^{-1} =\varphi^8$. We deduce
$B_2(\xi^{-1}) = \frac{1}{33} (61 - 59 \zeta_5 - 23 \zeta_5^2 - 47 \zeta_5^3 + 13 \zeta_5^4)$. Using SAGE, one can check $5 \mid N_{\Q(\zeta_5)/\Q} (B_2(\xi^{-1}))$. Thus, $ 5 | N(\tilde{\beta}_{\Ga_0(N),\varphi} )$ and therefore $5$ is a possible residual characteristic of an Eisenstein ideal $\m$ whose associated Eisenstein series is $E_\varphi$. 
Since the order of $\varphi$ is $10$, the character $\overline{\epsilon}_\m$ is the quadratic character $(\Z/11\Z)^\times \to \F_5^\times$. 
The Fourier expansion of $E_\varphi = \sum\limits_{n=1}^\infty \varphi(n) (\sum\limits_{d \mid n} \varphi^8(d) d) q^n$ is given by
$$q + \zeta_{10}( 1 + 2 \zeta_5^4) q^2 - \zeta_{10}^3(1 + 3 \zeta_5^2) q^3 + \zeta_5(1+ 2 \zeta_5^4 + 4 \zeta_5^3) q^4 + \zeta_5^2( 1 + 5 \zeta_5) q^5 - \zeta_5^2 ( 1 + 2 \zeta_5^4 + 3 \zeta_5^2 + 6 \zeta_5) q^6 + O(q^7).	$$
Let $f(z)$ denote the newform of $\Ga_0(121)$ whose LMFDB label is $121.2.a.d$. Then, we have 
$$E_\varphi(z) \equiv f(z) \pmod{\lambda}, $$ where 
$\lambda = (1-\zeta_5)$ be the prime above $5$ in $\Q(\zeta_5) = \Q(\zeta_{10})$.  We conclude that the corresponding non-rational Eisenstein maximal ideal $\m$ of characteristic $5$, with $\T/\m \cong \F_5$, is given by
 $$\m= \langle 5, U_{11}, \{ T_r -1 -r \}_{\text{ primes } r \equiv 1, 3 , 4 , 5 , 9 \pmod{11}}, \{ T_s +1+s \}_{\text{ primes } s \equiv 2, 6 , 7 , 8 , 10 \pmod{11}} \rangle.$$

\end{example}
 \begin{example}
Take $N=725=5^2 \cdot 29$. There are $3$ nontrivial Dirichlet characters of conductor $5$, namely $\{ \varphi, \varphi^2, \varphi^3 \}$, where $\varphi$ is given by $\varphi(2)= i$. Obviously, $\varphi^{-1} = \varphi^3$ and $\varphi^2$ are the quadratic Dirichlet characters of conductor $5$. Observe that,  $\tau(\varphi^2) = \sqrt{5}$, $B_2(\varphi^2) = \frac{4}{5}$, $\tau(\varphi) = i \sqrt[4]{-15+20i}$ and $\tau(\varphi^3)= i \sqrt[4]{-15-20i}$. Note that, in $\Q(i, \zeta_5)$, the only prime divisors of $\tau(\varphi)$ and $\tau(\varphi^3)$ are primes lying above $5$.
The space of non-rational Eisenstein series is $6$-dimensional and given by $E_{\varphi}^{\mrm{ord}_{29}}=E_{\varphi,1,29}, E_{\varphi}^{\crit_{29}}=E_{\varphi, 29,1}, E_{\varphi^2}^{\ord_{29}}=E_{\varphi^2,1,29}, E_{\varphi^2}^{\crit_{29}}=E_{\varphi^2, 29,1}, E_{\varphi^3}^{\ord_{29}}=E_{\varphi^3,1,29}$, and $E_{\varphi^3}^{\crit_{29}}=E_{\varphi^3, 29,1}$.
Hence, compute $\tilde{\beta}_{\Ga_0(N), \varphi^i, 29, 1}$ and see that $7$ (is the only prime other than $2,3$ and $5$ ) that divides $N(\tilde{\beta}_{\Ga_0(N), \varphi^i, 29, 1})$ for $i=1,2,3$. 

The Fourier expansion of $E_{\varphi^2}$ is given by $E_{\varphi^2}(z) = \sum\limits_{n=1}^\infty \varphi^2(n) (\sum\limits_{d \mid n} d) q^n  = q -3q^2 - 4q^3+ 7q^4+ 12q^6-8q^7-15q^8+13q^9+12q^{11} -28q^{12} -14 q^{13}+ 24q^{14}+O(q^{16}).$
Consider the newform $f(z)$ of level $\Ga_0(725)$ with LMFDB label $725.2.a.b$ and the coefficient field $\Q(\beta)$, where $\beta =
\sqrt{2}$. The Fourier expansion of $f$ is given by
 \begin{align*}
 q + ( 1 + \beta ) q^{2} + ( -1 - \beta ) q^{3} + ( 1 + 2 \beta ) q^{4} + ( -3 - 2 \beta ) q^{6} + 2 \beta q^{7} + ( 3 + \beta ) q^{8} + \\
 2 \beta q^{9} + ( 1 - \beta ) q^{11} + ( -5 - 3 \beta ) q^{12} + ( 1 + 2 \beta ) q^{13} + ( 4 + 2 \beta ) q^{14} +O(q^{16}).
 \end{align*}
Let $\wp = (3-\sqrt{2})$ be the prime in $\Z[\sqrt{2}]$ lying above $7$. We obtain
$$ E_{\varphi^2}^{\crit_{29}}(z)  \equiv E_{\varphi^2}^{\ord_{29}}(z) \equiv f(z) \pmod{\wp}.$$
We conclude that the non-rational Eisenstein maximal ideal of characteristic $7$, corresponding to  $E_{\varphi^2}^{\crit_{29}}$, with $\T/\m \cong \F_7$, is given by
$$\m = \langle 7, U_{5}, U_{29}-1, \{ T_r - 1 - r \}_{\text{ primes } r \equiv 1, 4 \pmod{5}}, \{ T_s + 1 + s \}_{\text{ primes } s \equiv 2, 3 \pmod{5}} \rangle. $$

We now consider the newform $g(z)$ of level $\Ga_0(725)$ with LMFDB label $725.2.a.l$, with the coefficient field $K$ determined by  $m(x)=x^6-13x^4+41x^2-1$. The Fourier expansion of $g$ is given by
\begin{align*}
g(q)= q + \beta_{1} q^{2} + ( \beta_{1} - \beta_{4} ) q^{3} + ( 2 + \beta_{2} ) q^{4} \\
+ ( 2 - \beta_{3} ) q^{6} + ( -\beta_{1} + \beta_{4} + \beta_{5} ) q^{7} + ( 2 \beta_{1} + \beta_{4} + \beta_{5} ) q^{8} + ( 2 - \beta_{2} - \beta_{3} ) q^{9} +O(q^{10}).
\end{align*}
Here,  $\beta_1$ is a root of $m(x)$, $\beta_2 = \beta_1^2 - 4 $, 
$\beta_{3}=( \beta_1^{4} - 8 \beta_1^{2} + 5)/2$, 
$\beta_{4}=( \beta_1^{5} - 12 \beta_1^{3} + 35 \beta_1)/2$ and 
$\beta_{5}=( -\beta_1^{5} + 14 \beta_1^{3} - 47 \beta_1)/2$.
The Fourier expansions of $E_\varphi$ and $E_{\varphi^3}$ are given by
$$E_\varphi(z) = \sum\limits_{n=1}^\infty \varphi(n) (\sum_{d \mid n} \varphi^2(d) d)q^n, E_{\varphi^3}(z) = \sum_{n=1}^\infty \varphi^3(n) (\sum_{d \mid n} \varphi^2(\frac{1}{d}) d)q^n.$$
Note that $E_\varphi$ and $E_{\varphi^3}$ are defined over $\Q(i)$, and they are Galois $\Gal(\Q(i)/\Q)$-conjugate. Using SAGE, we see that 
the prime factorization of $(7)$ in $K(i)$ is given  by $\p_1^2\p_2^2\p_3\p_4$, where $\p_1 = ( (\frac{1}{2} \beta_1^5 + \frac{1}{4} \beta_1^4 - 6 \beta_1^3 - 2 \beta_1^2 + \frac{35}{2} \beta_1 + \frac{13}{4}) i - (\frac{1}{4} \beta_1^5 - \frac{7}{2} \beta_1^3 - \frac{1}{2} \beta_1^2 + \frac{47}{4} \beta_1 + \frac{5}{2}))$ and $\p_2 = ( (\frac{1}{2} \beta_1^5 + \frac{1}{4} \beta_1^4 - 6 \beta_1^3 - 2 \beta_1^2 + \frac{35}{2} \beta_1 + \frac{13}{4}) i + (\frac{1}{4} \beta_1^5 - \frac{7}{2} \beta_1^3 - \frac{1}{2} \beta_1^2 + \frac{47}{4} \beta_1 + \frac{5}{2}))$. Then,
$$ E_{\varphi}^{\crit_{29}} \equiv E_{\varphi}^{\ord_{29}} \equiv g(z) \pmod{\p_1}\quad \mrm{and}\quad E_{\varphi^3}^{\crit_{29}} \equiv E_{\varphi^3}^{\ord_{29}} \equiv g(z) \pmod{\p_2}. $$
\end{example}
\begin{example}
Take $N=234=3^2 \cdot 2 \cdot 13$ and $\varphi$, which are the non-trivial Dirichlet characters of the conductor $3$ given by $\varphi(2)=-1$. 
There are $4$ non-rational Eisenstein series (Eisenstein series whose Fourier coefficients are not rational) of weight $2$ and level $\Ga_0(N)$, namely $E_{\varphi}^{\ord_2, \ord_{13}}, E_{\varphi}^{\ord_2, \crit_{13}}, E_{\varphi}^{\crit_2, \ord_{13}}$ and $E_{\varphi}^{\crit_2, \crit_{13}}$. If
$E_\varphi(z)= \sum\limits_{n=1}^\infty \varphi(n) (\sum\limits_{d \mid n} d) q^n = q-3q^2+7q^4-6q^5+8q^7-15q^8 + 18q^{10} -12 q^{11}+ 14q^{13} - 24q^{14} + O(q^{16}),$
then
$$E_\varphi^{\crit_{13}} = E_\varphi(z) - E_\varphi(13z) = q-3q^2+7q^4-6q^5+8q^7-15q^8 + 18q^{10} -12 q^{11}+ 13q^{13} - 24q^{14} + O(q^{16}),$$
$$E_\varphi^{\ord_{13}} = E_\varphi(z) - 13 E_\varphi(13z) = q-3q^2+7q^4-6q^5+8q^7-15q^8 + 18q^{10} -12 q^{11}+ q^{13} - 24q^{14} + O(q^{16}).$$
$$E_{\varphi}^{\ord_2, \crit_{13}} = E_\varphi^{\crit_{13}}(z) + 2 E_\varphi^{\crit_{13}}(2z) = q-q^2+q^4-6q^5+8q^7-q^8 + 6q^{10} -12 q^{11}+ 13q^{13} - 8q^{14} + O(q^{16}),$$
$$E_{\varphi}^{\crit_2, \crit_{13}} = E_\varphi^{\crit_{13}}(z) +  E_\varphi^{\crit_{13}}(2z) = q-2q^2+4q^4-6q^5+8q^7-8q^8 + 12q^{10} -12 q^{11}+ 13q^{13} - 16q^{14} + O(q^{16}),$$
$$E_{\varphi}^{\ord_2, \ord_{13}} = E_\varphi^{\ord_{13}}(z) + 2 E_\varphi^{\ord_{13}}(2z) = q-q^2+q^4-6q^5+8q^7-q^8 + 6q^{10} -12 q^{11}+ q^{13} - 8q^{14} + O(q^{16}),$$
$$E_{\varphi}^{\crit_2, \ord_{13}} = E_\varphi^{\ord_{13}}(z) +  E_\varphi^{\ord_{13}}(2z) = q-2q^2+4q^4-6q^5+8q^7-8q^8 + 12q^{10} -12 q^{11}+ q^{13} - 16q^{14} + O(q^{16}).$$
Proceeding as before, we can calculate
$$\tilde{\beta}_{\Ga_0(N), \varphi,1, 26} = -4 \sqrt{-3}, \tilde{\beta}_{\Ga_0(N), \varphi, 2, 13} = -12 \sqrt{-3}, \tilde{\beta}_{\Ga_0(N), \varphi, 26, 1} = -7 \cdot 24 \sqrt{-3}, \tilde{\beta}_{\Ga_0(N), \varphi, 13, 2} = -7 \cdot 8 \sqrt{-3}.$$
Thus, the prime $7$ cannot be a residual characteristic of the Eisenstein maximal ideal whose associated Eisenstein series is either $E_{\varphi}^{\ord_2, \ord_{13}}$ or $E_{\varphi}^{\crit_2, \ord_{13}}$. However,  $7$ can be the residual characteristic of the Eisenstein maximal ideal whose associated Eisenstein series is either $E_{\varphi}^{\ord_2, \crit_{13}}$ or $E_{\varphi}^{\crit_2, \crit_{13}}$.

Note that there are $5$ new forms of weight $2$, level $\Ga_0(234)$ (that is, with LMFDB label $234.2.a.a$ to $234.2.a.e$), all of which have rational Fourier coefficients.   One can easily check $E_{\varphi}^{\ord_2, \crit_{13}} \equiv 234.2.a.b \pmod 7$.  We conclude that the corresponding non-rational Eisenstein maximal ideal of characteristic $7$, with $\T/\m \cong \F_7$, is given by
 $$\m= \langle 7, U_{3}, U_2 +1 , U_{13} -1, \{ T_r -1 -r \}_{\text{ primes } r \equiv 1 \pmod{3}}, \{ T_s +1+s \}_{\text{ primes } s \equiv 2 \pmod{3}} \rangle.$$

By the Fourier expansion, we see $E_{\varphi}^{\crit_2, \crit_{13}}$ is not congruent to any of the above newforms modulo $7$. There exists no maximal Eisenstein ideal $\m$ of the residual characteristic $7$ that corresponds to $E_{\varphi}^{\crit_2, \crit_{13}}$.

\end{example}



\bibliographystyle{crelle2}
\bibliography{Eisensteinquestion1.bib}
\end{document}